%% file: IlluminateCINT_arXiv.tex
\documentclass[final]{siamltex}
\usepackage{xcolor}
\usepackage{amsmath}
\usepackage{amssymb}
\usepackage{mathrsfs}
\usepackage{graphicx}
\usepackage{hyperref}
\usepackage{bm}
\usepackage{ucs}
\usepackage[utf8x]{inputenc}
\usepackage{wrapfig}
\usepackage{color}
\usepackage{float}
\usepackage{epstopdf}

\usepackage{amsmath,amstext,amssymb,graphicx,hyperref}
\usepackage{verbatim}
\usepackage{color}
\usepackage{algorithm,algorithmic}

\usepackage{mathtools}

\usepackage{epstopdf}
\usepackage{color}
\usepackage{graphicx}
\usepackage{float}
\usepackage[section]{placeins}

\usepackage{tikz}
\usetikzlibrary{arrows}

\input{mathmacros.tex}

\def\dsp{\displaystyle}

\def\Gc{{\cal G}}

\def\Zc{{\cal Z}}
\def\norm#1{\|#1\|}

\def\vect#1{\mbox{\boldmath{$#1$}}}

\def\vx{\vec{\vect x}}
\def\vy{\vec{\vect y}}
\newcommand{\bx}{\vect x}

\newcommand{\cC}{\mathcal{C}}

\def\mE{\mathbb{E}}

\def\mC{\mathbb{C}}

\def\we{\widehat{e}}

\newcommand{\Pm}{\vect P}
\newcommand{\Mm}{\vect M}

\newcommand{\eps}{\varepsilon}
\renewcommand{\epsilon}{\eps}

\def\bfrho{\mbox{\boldmath$\rho$}}

\DeclarePairedDelimiter{\abs}{\lvert}{\rvert}

\begin{document}

\title{Multifrequency interferometric imaging with intensity-only measurements}
\author{Miguel Moscoso, Alexei Novikov, George Papanicolaou and Chrysoula Tsogka}
\maketitle
\begin{abstract}
We propose an illumination strategy for interferometric imaging 
that allows for robust depth recovery from intensity-only measurements. 
For an array with colocated sources and receivers, we show that all the possible interferometric data  
for multiple sources, receivers and frequencies can be recovered from intensity-only measurements 
provided that we have sufficient source location and frequency illumination diversity. 
There is no need for phase reconstruction in this approach.
Using interferometric imaging methods 
we show that in homogeneous media there is no loss of resolution when imaging with intensities-only.
If in these imaging methods we reduce incoherence by restricting the multifrequency interferometric data 
to nearby array elements and nearby frequencies we obtain robust images in weakly
inhomogeneous background media with a somewhat reduced resolution.

\end{abstract}
\begin{keywords}
array imaging, phase retrieval.
\end{keywords}

\section{Introduction}

Coherent array imaging when the phases of the signals received at the array cannot be measured is a difficult problem because much of
the information about the sought image is contained in the lost phases. Imaging without phases arises in many applications such as
crystallography \cite{Harrison93, Shechtman_2015}, ptychography \cite{Rodenburg08} 
and optics \cite{Nugent10, Goodman05, Walther63, Poon2006} 
where images are formed from the spectral intensities. In most of these cases, the media through which the probing signals propagate are assumed to be homogeneous. 

The  earliest and most widely used methods for imaging with intensity-only measurements
are alternating projection algorithms \cite{Fienup82, GS72}. The basic idea is to project the iterates on the intensity data sequentially in both the real and 
the Fourier spaces. Although these algorithms are very efficient for reconstructing
the missing phases in the data, and performance is often good in practice, they do not always converge to the true, missing phases. This is especially so if strong constrains or prior information about the object to be imaged, such as spatial support and non-negativity, are not reliably available \cite{Fienup13, Shechtman_2015}.
We do not use phase retrieval methods here. Instead, we exploit illumination diversity to recover all missing phase information and then image holograpically.
We assume that the missing phase information is largely coherent, that is, it is not so corrupted by medium inhomogeneities and measurement noise so that even when recovered it will not be useful in coherent imaging. We address in detail coherence issues in this paper.


\subsection*{The array imaging problem}
\begin{figure}[htbp]
\begin{center}
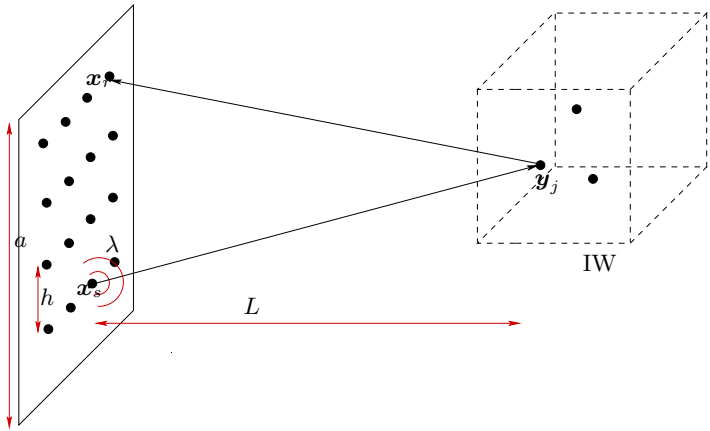
\caption{General setup of an array imaging problem. The transducer at $\vect x_s$  emits a probing signal and the reflected intensities are recorded at $\vect x_r$. The scatterers located at ${\vect y}_{j}$, $j=1,\dots,M$ 
are at distance $L$ from the array and inside the image window IW. }
\label{cs_3d_illustration}
\end{center}
\end{figure}

We consider an active array of size $a$ consisting of $N$ transducers separated by a distance $h$ which is of the order of the central wavelength $\lambda_0$ 
of the probing signals. The transducers emit probing signals of different frequencies $\omega_l$, $l=1,\dots,S$, from positions ${\vect x}_s$ and record the reflected intensities at positions ${\vect x}_r$, $s,r = 1,2,\ldots,N$ (see Figure~\ref{cs_3d_illustration}).  

The goal
is to determine the positions ${\vect y}_{j}$ and reflectivities $\alpha_j\in\mC$, $j=1,\dots,M$, of a set of $M$ point-scatterers within a region of interest, called the image window (IW), which is at a distance $L$ form the array, as shown in Figure \ref{cs_3d_illustration}.

\subsection*{Holographic Imaging}   
A holographic imaging approach with intensity-only measurements is presented in ~\cite{Novikov14,Moscoso16}. The main
idea is to exploit illumination diversity by
designing illumination strategies that recover the missing phase information from intensity-only measurements.
It was shown in~\cite{Novikov14,Moscoso16} how 
 by using an appropriate protocol of illuminations and the polarization identity, the single frequency 
matrix $\Mm(\omega) = {\Pm(\omega)}^*  \Pm(\omega) $ can be determined from intensity-only measurements at that frequency.
Here $\Pm(\omega)= [ P({\vect x}_r,{\vect x}_s;\omega)]_{r,s=1}^{N}$ is the full array response matrix of the imaging system, including phases, with
${\vect x}_r,{\vect x}_s$ the receiver and source locations, and $\omega$ the radian frequency.
The matrix $\Mm(\omega)$ is called the time reversal matrix as it arises in
ultrasonic time reversal experiments \cite{ Fink_2001, Fink_1992} and has been studied extensively there \cite{BPT_2003, BPZ_2002}. 
We will refer to $\Mm(\omega)$ as the single frequency interferometric data matrix. 
Once we have this data matrix we can image with the DORT method \cite{Prada96} which uses the eigenvectors of $\Mm(\omega)$, or MUSIC \cite{Schmidt86, Gruber04}, which also uses the eigenvectors of $\Mm(\omega)$. Here DORT and MUSIC are the acronyms:
Decomposition de l'Operateur de Retournement Temporel (Decomposition of the Time Reversal Operator), 
and Multiple Signal Classification, respectively. 
These are phase-sensitive imaging methods that involve only phase differences contained in $\Mm(\omega)$ and, therefore, 
they provide interferometric information. The illumination strategies in
~\cite{Novikov14,Moscoso16} are a form of digital holography \cite{Wolf69, Schotland2009, Goodman05, Poon2006} since the resulting image does have phase information. As already noted, we do not
reconstruct phases from intensity measurements, but rather we recover the missing phase information using illumination
diversity.

Imaging with $\Mm(\omega)$ 
at a single frequency $\omega$  is not robust relative to small perturbations in the unknown phases unless the array is very large \cite{Chai16}. The perturbations can come from 
medium inhomogeneities or from the discretization of the image window. Having $\Mm(\omega_l)$ at multiple frequencies $\omega_l,~l=1,2,\ldots,S$  still may not provide 
robustness with respect to depth in the image. Methods that use the eigenvectors frequency by frequency, as in MUSIC, are not robust.

\subsection*{Interferometric robust imaging}   
It is known \cite{BPT-05, BPT-ADA} that we can get image robustly if we have interferometric data 
\begin{equation}
\label{eq:data_full}
d((\vect x_r, \vect x_{r'}),(\vect x_s, \vect x_{s'}),(\omega,\omega'))=\overline{P(\vect x_r, \vect x_s, \omega)}  P(\vect x_{r'}, \vect x_{s'}, \omega')
\end{equation}
at multiple frequency pairs $(\omega,\omega')$, receiver location pairs $(\vect x_r, \vect x_{r'})$ and source location pairs $(\vect x_s, \vect x_{s'})$. 
The main result of this paper is that we can recover such data $d((\vect x_r, \vect x_{r'}),(\vect x_s, \vect x_{s'}),(\omega,\omega'))$ 
for pairs of arguments from intensity-only measurements.
Here receivers and transmitters are colocated in the same array. When the imaging system has separate transmitting and receiving arrays then we can recover only 
single receiver elements, one receiver at a time, 
\begin{equation}
\label{eq:data_single}
d((\vect x_r,\vect x_r),(\vect x_s, \vect x_{s'}),(\omega,\omega')) = \overline{P(\vect x_r, \vect x_s, \omega)}  P(\vect x_{r}, \vect x_{s'}, \omega')
\end{equation}
for all pairs of frequencies, and source locations from intensity-only measurements. 

In a homogeneous medium, imaging with  $d((\vect x_r, \vect x_{r'}),(\vect x_s, \vect x_{s'}),(\omega,\omega'))$
can be done by 
\begin{equation}
\label{eq:interf_full}
\begin{array}{ll}
\dsp  {\cal I}^{Interf}({\vect y}^s)  &\dsp = \sum_{\vect x_s, \vect x_{s'}} \sum_{\vect x_r, \vect x_{r'}} \sum_{\omega_l, \omega_{l'}}
d((\vect x_r, \vect x_{r'}),(\vect x_s, \vect x_{s'}),(\omega_l,\omega_{l'})) \\
&\times G_0(\vect x_r, \vect y^s, \omega_l ) G_0(\vect x_s, \vect y^s, \omega_l )\overline{G_0(\vect x_{r'}, \vect y^{s}, \omega_{l'} )}
\overline{G_0(\vect x_{s'}, \vect y^{s}, \omega_{l'} )}
\end{array}
\end{equation}
with $G_0(\vect x_r, \vect y^s, \omega_l )$ the free space Green's function for the Helmholtz equation (see expression~\eqref{homo_green} below), and $\vect y^s$ a point in the image window IW. 
Replacing the data by its expression (\ref{eq:data_full}) we note that ${\cal I}^{Interf}({\vect y}^s) $ equals the square of the Kirchhoff Migration  imaging
function
\begin{equation}
\begin{array}{ll}
\dsp  {\cal I}^{Interf}({\vect y}^s) &\dsp = \left| \sum_{\vect x_s} \sum_{\vect x_r} \sum_{\omega_l} \overline{P(\vect x_r, \vect x_s, \omega_l )}  G_0(\vect x_r, \vect y^s, \omega_l ) G_0(\vect x_s, \vect y^s, \omega_l )\right|^2 \\[16pt]
&\dsp = 
\left| {\cal I}^{KM}({\vect y}^s) \right|^2.
\end{array}
\label{eq:INTERF}
\end{equation}
Here, the Kirchhoff migration functional 
\begin{equation}
{\cal I}^{KM}({\vect y}^s)= \sum_{\vect x_s} \sum_{\vect x_r} \sum_{\omega_l} \overline{P(\vect x_r, \vect x_s, \omega_l )}  
G_0(\vect x_r, \vect y^s, \omega_l ) G_0(\vect x_s, \vect y^s, \omega_l )
\label{eq:KM}
\end{equation}
is simply the back propagation of the array response matrix in a homogeneous medium, both for source and receiver points.
Note that it is the {\it square} of the Kirchhoff migration functional that we obtain with intensity-only measurements. 

The main result of this paper can now be restated as follows. For colocated source and receivers on a single array, we can obtain full-phase, holographic images from intensity-only measurements by 
exploiting illumination and frequency diversity. That is, in a homogeneous medium there is no loss of resolution when imaging only with intensities if we have sufficient source and frequency illumination diversity.

\subsection*{Inhomogeneous background medium and CINT} 
In a randomly inhomogeneous medium it is well known \cite{BPT-05,BGPT-rtt} that Kirchhoff migration does not work well even if we have the full array 
response matrix, phases included. The interferometric functional (\ref{eq:interf_full}) does not work either, since it is just the square of the 
Kirchhoff migration functional.
For weakly inhomogeneous random media we can image with the coherent interferometric (CINT) functional, which has the form
\begin{equation}
\label{eq:CINT}
\begin{array}{ll}
\dsp  {\cal I}^{CINT}({\vect y}^s)  &\dsp = \hspace*{-0.7cm}
\sum_{\begin{array}{c} 
\scriptsize \vect x_s, \vect x_{s'} \\
\scriptsize |\vect x_s - \vect x_{s'}| \le X_d \end{array}} 
\hspace*{-0.3cm}
\sum_{
\begin{array}{c} 
\scriptsize \vect x_r, \vect x_{r'}
\\
\scriptsize  |\vect x_r - \vect x_{r'}| \le X_d \end{array}} 
\hspace*{-0.3cm}
\sum_{
\begin{array}{c} 
{\scriptsize \omega_l, \omega_{l'}} \\
{\scriptsize | \omega_l -\omega_{l'} | \le \Omega_d} 
\end{array}}
\hspace*{-0.5cm}
\overline{P(\vect x_r, \vect x_s, \omega_l )}  P(\vect x_{r'}, \vect x_{s'}, \omega_{l'}) \\[24pt]
&\times G_0(\vect x_r, \vect y^s, \omega_l ) G_0(\vect x_s, \vect y^s, \omega_l )\overline{G_0(\vect x_{r'}, \vect y^{s}, \omega_{l'} )}
\overline{G_0(\vect x_{s'}, \vect y^{s}, \omega_{l'} )},
\end{array}
\end{equation}
assuming full phase information is available by the array response matrix $P(\vect x_r, \vect x_s, \omega_l )$.

We note that the CINT imaging functional uses the same data 
that we get from intensity-only measurements, $d((\vect x_r, \vect x_{r'}),(\vect x_s, \vect x_{s'}),(\omega_l,\omega_{l'}))$,
but now the sums are restricted to nearby pairs of sources, receivers and frequencies. 
The main idea in CINT is that since the waves propagate in fluctuating media their phases are distorted and the multifrequency interferometric data
$d((\vect x_r, \vect x_{r'}),(\vect x_s, \vect x_{s'}),(\omega_l,\omega_{l'}))$, given by (\ref{eq:data_full}), 
remain coherent only over small frequency and space offsets. We call decoherence distance, $X_d$, and decoherence frequency, $\Omega_d$, the 
largest spatial and frequency intervals, respectively, over which the multifrequency interferometric data $\overline{P(\vect x_r, \vect x_s, \omega_l )}  P(\vect x_{r'}, \vect x_{s'}, \omega_{l'})$ remain coherent. We mean by coherent that the distortion of phases by inhomogeneities is weak.
These decoherence parameters $X_d$ and $\Omega_d$ are a-priori unknown. They can be estimated directly from the data using statistical techniques like the
 variogramm \cite{ripley}.  However, optimal imaging results are obtained when $X_d$ and $\Omega_d$ are estimated {\em on the fly} during the image formation process as in 
 adaptive CINT \cite{BPT-ADA}. 

Most importantly, what matters in imaging is that statistical stability is gained by this appropriate restriction of the multifrequency interferometric data. 
The statistical stability of CINT is shown in \cite{BGPT-rtt}. Specifically, this means that the variance of the image is small compared to its mean square, 
with respect to the realizations of the fluctuating medium. Therefore, the imaging results do not depend on any particular realization of the random medium. 
However, statistical stability comes at the cost of loss in resolution: cross-range resolution now becomes $ \lambda_0 L/X_d$ instead of $ \lambda_0 L/a$,
and resolution in range or depth is $c_0/\Omega_d$ instead of $c_0/B$, with $a$ being the array size and $B$ the bandwidth. Typically, 
we have $\Omega_d < B$ and $X_d < a$, and often
$\Omega_d \ll B$ and $X_d \ll a$.  

\subsection*{Single receiver multifrequency interferometric imaging (SRINT)} 
Restricting the data to intensity-only measurements at a single receiver, we obtain $d((\vect x_r, \vect x_{r}),(\vect x_s, \vect x_{s'}),(\omega_l,\omega_{l'}))$.
Using only data from a single receiver, we introduce 
the following single receiver coherent interferometric imaging (SRINT) functional 
\begin{equation}
\begin{array}{ll}
\dsp  {\cal I}^{SRINT}({\vect y}^s)  &\dsp = \hspace*{-0.7cm}
\sum_{\begin{array}{c} 
\scriptsize \vect x_s, \vect x_{s'} \\
\scriptsize |\vect x_s - \vect x_{s'}| \le X_d \end{array}} 
\hspace*{-0.3cm}
\sum_{
\begin{array}{c} 
{\scriptsize \omega_l, \omega_{l'}} \\
{\scriptsize | \omega_l -\omega_{l'} | \le \Omega_d} 
\end{array}}
\hspace*{-0.5cm}
d((\vect x_r, \vect x_{r}),(\vect x_s, \vect x_{s'}),(\omega_l,\omega_{l'})) \\[24pt]
&\times G_0(\vect x_r, \vect y^s, \omega_l ) G_0(\vect x_s, \vect y^s, \omega_l )\overline{G_0(\vect x_{r}, \vect y^{s}, \omega_{l'} )}
\overline{G_0(\vect x_{s'}, \vect y^{s}, \omega_{l'} )}\, .
\end{array}
\label{eq:SRINT}
\end{equation}
Note that there is no sum over receivers here. We only have one receiver at $\vect x_r$.

We use in this paper the SRINT imaging functional (\ref{eq:SRINT}) in a computationally efficient form involving only matrix multiplications (see~\eqref{eq:AMA} below). 
The key idea is the introduction of a mask, {\em i.e.}, a matrix that 
is composed of zeros and ones only, depending on the spacing between the indices in the matrix so as to restrict the data $d((\vect x_r, \vect x_{r}),(\vect x_s, \vect x_{s'}),(\omega_l,\omega_{l'}))$ to the ones satisfying the constraints $|\vect x_s - \vect x_{s'}| \le X_d$ and  $ | \omega_l -\omega_{l'} | \le \Omega_d$ that should be used in SRINT imaging.

The performance of the proposed interferometric method is explored with numerical simulations in an optical (digital) microscopy regime. We observe in the simulations that 
in homogeneous media we can image with the same resolution as if phases where recorded and the method is robust with respect to the discretization of the image window. When 
the ambient medium is weakly inhomogeneous the interferometric approach removes some of the uncertainty in the data due to the fluctuating phases, which tends to stabilize the 
images and this is seen clearly in the simulations. We also compare the performance of the interferometric approach with  MUSIC which is shown to be sensitive to phase errors 
and does not provide robust results unless the illuminating and receiving arrays are large \cite{Chai16}. The fact that the SRINT imaging functional, which uses data obtained with 
intensity-only measurements, gives images that are robust to weak fluctuations in the ambient medium is another main result in this paper. It is surprising that such robust, 
holographic imaging can be obtained with intensity-only measurements.

The paper is organized as follows.
In Section~\ref{sec:model} we formulate our data model for intensity-only
measurements. 
In Section~\ref{sec:single} we formulate our single frequency data model for intensity-only
measurements.
In Section~\ref{sec:multi} we formulate our multi-frequency data model for intensity-only
measurements.
In Section~\ref{sec:alexei} we describe our illumination strategy for holographic imaging,
and in Section~\ref{sec:imaginghomo} we describe the imaging algorithms of single receiver interferometry (SRINT) and multiple signal classification (MUSIC).  
In Section~\ref{sec:numerics}, we explore with numerical simulations the robustness of the imaging methods in an optical (digital) microscopy regime. 
In Section~\ref{sec:imagingrandom} we discuss aspects of imaging in inhomogeneous background media.
Section~\ref{sec:conclusions} contains our conclusions.

\section{Single frequency data models and imaging with phases}
\label{sec:model}
We consider the array  imaging configuration of Figure~\ref{cs_3d_illustration}, where an array consisting of $N$ transducers is used to probe the image window (IW). Our goal is to  
reconstruct a sparse scene consisting of $M$ point-scatterers. The unknowns are both the locations 
${\vect y}_{j}$ and the reflectivities $\alpha_j\in\mC$, $j=1,\dots,M$, of the scatterrers. 

For imaging purposes the IW is discretized using a uniform grid of $K$ points ${\vect y}_k$, $k=1,\ldots,K$. We assume that $K>N$ and often we have $K\gg N$.
By point-like scatterers we mean very small scatterers compared to the central wavelength.
We also assume that the scatterers are far apart or are weak, so multiple scattering between them is negligible. We refer to \cite{Chai14} for array imaging problems with multiple scattering.

If the scatterers are far apart or the reflectivities are small, the interaction between them is weak and the Born approximation is applicable. In this case, the response at ${\vect x}_r$ (including phases) due to a pulse of angular frequency $\omega_l$ sent from ${\vect x}_s$, and reflected by the $M$ scatterers,  is given by 
\begin{equation}
\label{response}
P({\vect x}_r,{\vect x}_s;\omega_l)=\sum_{j=1}^M\alpha_j G ({\vect x}_r,{\vect y}_j; \omega_l) \, ,
G ({\vect y}_{j},{\vect x}_s; \omega_l)\, ,
\end{equation}
where $G ({\vect x},{\vect y}; \omega_l)$ denotes the Green's function in a general, inhomogeneous medium that characterizes the propagation of a signal of angular frequency $\omega_l$ from point ${\vect x}$ to point ${\vect y}$. In a homogeneous medium we denote the Green's function by $G_0 ({\vect x},{\vect y}; \omega_l)$
and it is given by
\begin{equation}
\label{homo_green}
 G_0(\vec{\vect x},\vec{\vect y}; \omega_l) = \frac{\exp\{i \omega_l \abs{\vec{\vect x} - \vec{\vect y}}/c_0\}}{4\pi \abs{\vec{\vect x} - \vec{\vect y}}} \, ,
\end{equation}
where $c_0$ is the speed of propagation.
We note that the number of pixels $K$ in the image window is typically large compared to the number of reflectors $M$, $K>M$,  and we also assume that $M<N$, the number of array transducers, so that $K>N>M$. Even though the data model (\ref{response}) that we use here is simple and somewhat stylized, it is quite flexible and can deal with complex reflectors, not just widely spaced point reflectors, especially when multiple scattering effects are included.

\begin{figure}[htbp]
\begin{center}
\includegraphics[scale=0.9]{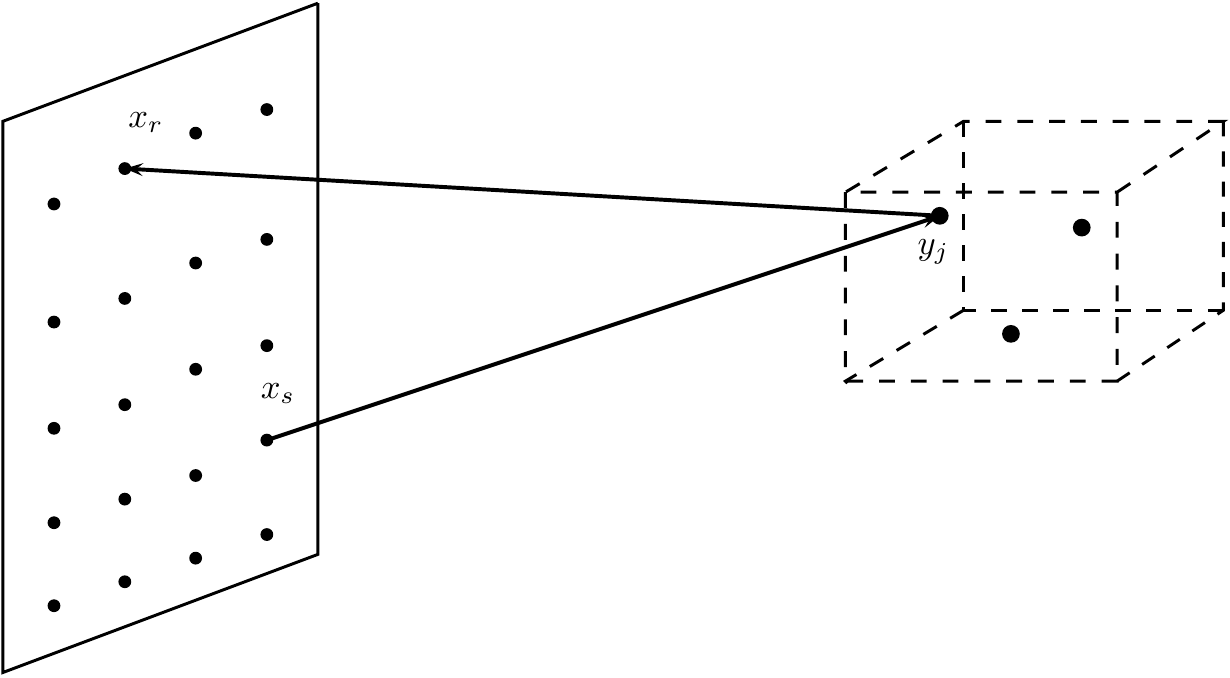}
\caption{Schematic of the components of the {\em array Green's function vector}  $\vect g ({\vect y};\omega_l)$ that represent the signals received on the array elements when a point source located at $\vect y$ sends a unit amplitude at frequency $\omega_l$.}
\label{fig:vectg}
\end{center}
\end{figure}

To write the single frequency data received on the array in a more compact form, we define the {\em array Green's function vector} 
$\vect g ({\vect y};\omega_l)$ at location ${\vect y}$ and with frequency $\omega_l$ as 
\begin{equation}
\label{GreenFuncVec}
\vect g ({\vect y};\omega_l)=[ G({\vect x}_{1},{\vect y};\omega_l),\cdots,G({\vect x}_{N},{\vect y};\omega_l)]^t\, ,
\end{equation}
where $.^t$ denotes  transpose. This vector represents the signal of frequency $\omega_l$ received at the array due to a point source at ${\vect y}$ as 
illustrated in Figure \ref{fig:vectg}. It can also be interpreted as the illumination vector  at the array that targets or beamforms to the position 
${\vect y}$ in the image window. 
 
We also introduce the generic {\it reflectivity vector} $\vect\rho=[\rho_{1},\ldots,\rho_{K}]^t\in\mC^K$, whose entries are  given by the possible reflectivities on the IW grid.
The generic reflectivities $\{\rho_j\}$ are zero except where the true reflectivity is not zero and then they equal $\{\alpha_j\}$. 
With this notation, we write the full response matrix  corresponding to a single frequency $\omega_l$ as a sum of rank-one matrices, so
\begin{equation}
\label{responsematrix}
\Pm(\omega_l)=
[P({\vect x}_r,{\vect x}_s;\omega_l)]
=\sum_{j=1}^M\alpha_j\vect g(\vect y_{j};\omega_l)\vect g^t(\vect y_{j};\omega_l)
=\sum_{k=1}^K\rho_{k}\vect g (\vect y_{k};\omega_l)\vect g^t(\vect y_{k};\omega_l).
\end{equation}
Using \eqref{GreenFuncVec}, we also define the $N\times K$ single frequency sensing matrix 
\begin{equation}\label{sensingmatrix}
\vect{\Gc}(\omega_l)=[\vect g(\vect y_1;\omega_l)\,\cdots\,\vect g (\vect y_K;\omega_l)] \,,
\end{equation}
whose column vectors are the signals received at the array due to point sources of frequency $\omega_l$ at the grid points $\vect y_k$, $k=1,\dots,K$. 
$\vect{\Gc}(\omega_l)$ maps a distribution of sources of frequency $\omega_l$ in the IW to the data of the same frequency received on the array. 
Using \eqref{sensingmatrix}, we write \eqref{responsematrix} in matrix form as
\begin{equation}\label{responsematrix1}
\Pm(\omega_l) =\vect\Gc(\omega_l) \diag(\bfrho )\vect\Gc^t(\omega_l).
\end{equation}
For a fixed frequency $\omega_l$, this matrix is  a linear transformation from the illumination space $\mC^N$ to the data space $\mC^N$. 
Indeed, for an illumination vector of frequency $\omega_l$
\begin{equation}
\vect f(\omega_l) =[f_{1}(\omega_l),\ldots, f_{N}(\omega_l)]^t\, ,  
\end{equation}
whose components are the signals $f_{1}(\omega_l),\ldots,f_{N}(\omega_l)$  sent from the transducers in the array, $\vect\Gc^t(\omega_l)\vect f(\omega_l)$ is the vector of size $K$ of signals of the same frequency at each grid point of the IW. These signals are reflected by the scatterers on the 
grid, with reflectivities given by the vector $\bfrho$, and then they are propagated back to the array by the matrix $\vect\Gc(\omega_l)$.

Summarizing, for a given illumination vector $\vect f^{(j)}(\omega_l)$, the single frequency, full phase $N-$vector array data  model 
$\vect b^{(j)}(\omega_l)$ is given by
\begin{equation}
\label{eq:sfimaging}
\vect b^{(j)}(\omega_l)=\vect P(\omega_l)\vect f^{(j)}(\omega_l)\,\, , \quad j=1,2,\dots, \aleph \, ,
\end{equation} 
whose components are
\begin{equation}
\label{eq:datasf}
b_r^{(j)}(\omega_l) 
= \sum_{k=1}^K \rho_k G({\vect
x}_r,{\vect y}_{k};\omega_l) \sum_{s=1}^N  G({\vect y}_{k},{\vect 
x}_s;\omega_l) f^{(j)}_s(\omega_l),\ r=1,2,\ldots,N \,.
\end{equation}
Here, $\sum_{s=1}^N  G({\vect y}_{k},{\vect 
x}_s;\omega_l) f^{(j)}_s(\omega_l)$ is the total illumination received at pixel $\vect y_k$ which is multiplied by the reflectivity $\rho_k$ and then propagated to the receiver $\vect x_r$ on the 
array with the Green's function 
$G({\vect x}_r,{\vect y}_{k};\omega_l) $. Note that in this model \eqref{eq:datasf}, the data is a linear function of both reflectivities $\bfrho$ and illuminations $\vect f^{(j)} (\omega_l) $.
Introducing the operator ${\mathcal A}_{f^{(j)}}(\omega_l)$ that transforms reflectivities $\bfrho$ to data ${\vect b}^{(j)}(\omega_l)$ we can also write the data in the following form 
\begin{equation}
\label{eq:datasf2}
{\vect b}^{(j)}(\omega_l) = {\mathcal A}_{f^{(j)}} (\omega_l){\vect \rho} \, ,
\end{equation}
where ${\mathcal A}_{f^{(j)}}(\omega_l)$ is defined as
\begin{equation}
\label{eq:defA}
\left[ {\mathcal A}_{f^{(j)}} (\omega_l) \right]_{rk} = G({\vect
x}_r,{\vect y}_{k};\omega_l)\sum_{s=1}^N  G({\vect y}_{k},{\vect 
x}_s;\omega_l) f^{(j)}_s(\omega_l) \, ,
\end{equation}
for  $r=1,2,\ldots,N$ and $k=1,2,\ldots,K$.

\subsection*{Single frequency imaging with phases}
The single frequency, full phase imaging problem can be stated as follows: Given a set of illuminations $\{\vect f^{(j)}(\omega_l)\}_{j=1,2,\dots,\aleph}$, 
determine the location and reflectivities of the scatterers from  the data \eqref{eq:sfimaging}. 
The representation (\ref{eq:datasf2}) of the data allows us to write the imaging problem as a linear system
\begin{equation}
\label{eq:linear}
{\mathcal A}^0_{f^{(j)}} (\omega_l){\vect \rho} = {\vect b}^{(j)}(\omega_l).
\end{equation}
Here, ${\mathcal A}^0_{f^{(j)}}(\omega_l)$ is our model for the operator that transforms reflectivities to data 
in homogeneous media and, therefore, 
we use in its definition the Green's function in a homogeneous, reference medium $G_0({\vect
x}_r,{\vect y}_{k};\omega)$ given by \eqref{homo_green}, i.e.,
\begin{equation}
\label{eq:defA0}
\left[ {\mathcal A}^0_{f^{(j)}} (\omega_l) \right]_{rk} = G_0({\vect
x}_r,{\vect y}_{k};\omega_l)\sum_{s=1}^N  G_0({\vect y}_{k},{\vect 
x}_s;\omega_l) f^{(j)}_s(\omega_l) \, ,
\end{equation}
for  $r=1,2,\ldots,N$ and $k=1,2,\ldots,K$. 

Several approaches can be used for computing the solution of the linear system \eqref{eq:linear}.  
Kirchhoff Migration consists in estimating the reflectivity by applying $\left({\mathcal A}^0_{f^{(j)}}(\omega_l)\right)^*$ 
to the data ${\vect b}^{(j)}(\omega_l)$, where the superscript $*$ denotes conjugate transpose. That is,
$$ \bfrho^{KM} =  \left({\mathcal A}^0_{f^{(j)}}(\omega_l)\right)^* {\vect b}^{(j)}(\omega_l) =  
\left({\mathcal A}^0_{f^{(j)}}(\omega_l)\right)^* {\mathcal A}_{f^{(j)}}(\omega_l) \bfrho \, .$$
We expect $ \bfrho^{KM}$ to be a good estimate of the true reflectivity $\bfrho$ when the model ${\mathcal A}^0_{f^{(j)}}(\omega_l)$ is close to the true operator ${\mathcal A}_{f^{(j)}}(\omega_l)$. It is known that
$ \left({\mathcal A}^0_{f^{(j)}}(\omega_l)\right)^* {\mathcal A}^0_{f^{(j)}}(\omega_l)$ is close to a diagonal matrix when the discretization of the image window conforms to the physical resolution limits which 
are $\lambda_0 L/a$ in cross-range and $\lambda_0 (L/a)^2$ in range (depth) in the single frequency case. A better estimate of the reflectivity is the least square solution of the linear system \eqref{eq:linear}.
This is an $\ell_2$ approach which is robust to additive uncorrelated noise and gives good results when the system \eqref{eq:linear} is overdetermined, that is for $N > K$. When the scene is sparse 
and the system  \eqref{eq:linear} is underdetermined,  the reflectivity can be estimated accurately and efficiently using the singular value decomposition of $\Pm(\omega_l)$ and an $\ell_1$ minimization approach as in \cite{CMP13}.

\section{Single frequency intensity-only data and imaging}\label{sec:single}
If only the single frequency intensities $\beta_i(\omega_l) = \abs{ b_i(\omega_l) }^2$ are recorded at the array, $i=1,\dots,N$,  
the imaging problem is to determine the location and reflectivities of the scatterers from the absolute values of each component in \eqref{eq:datasf2}, i.e., 
from 
\begin{equation}\label{dataI-1}
\vect \beta^{(j)}(\omega_l)=\diag \left(  \left( {\mathcal A}_{f^{(j)}} (\omega_l){\vect \rho} \right) \left({\mathcal A}_{f^{(j)}} (\omega_l){\vect \rho} \right)^* \right)
=\diag \left( {\mathcal A}_{f^{(j)}} (\omega_l){\vect \rho} {\vect \rho}^* {\mathcal A}^*_{f^{(j)}} (\omega_l) \right)
\end{equation}
for $j=1,2,\dots \, ,\aleph$. 
The corresponding imaging problem consists in seeking a $\bfrho$ solution of the system
 \begin{equation}\label{eq:non-linear}
\diag \left( {\mathcal A}^0_{f^{(j)}} (\omega_l){\vect \rho} {\vect \rho}^* \left( {\mathcal A}^0_{f^{(j)}} (\omega_l)\right)^* \right) = \vect \beta^{(j)}(\omega_l), \ \  j=1,2,\dots \, ,\aleph. 
\end{equation}

This single frequency imaging problem is now nonlinear in the unknown reflectivities $\bfrho$.  There are several ways in which this intensities-only imaging problem can be addressed. One is by convexification.
Because coherent imaging without phases is a non-convex, nonlinear problem, an alternative convex approach has been considered when the signals propagate 
in a homogeneous medium~\cite{CMP11,Candes13}.  In this approach, the original vector  
intensity-only, non-linear imaging problem is reformulated as a low-rank matrix linear imaging problem, which can be solved by using nuclear norm  minimization. 
This makes the intensity-only imaging problem convex over the appropriate matrix vector space and, therefore, the unique true solution can be found in the noise-free case \cite{Candes09, Recht10}.  However, because the original vector of unknown reflectivities $\bfrho$ is replaced by the rank one matrix $\bfrho \bfrho^*$, the size of the resulting optimization problem increases quadratically with the number of unknowns $K$. The computational cost of this approach is prohibitively high except in very special cases where the a-priori support and overall location of the reflectivities is known so that the window size $K$ can be reduced.  

Another approach to address this nonlinear imaging problem is to to use alternating projections \cite{Fienup13}, provided we are in the Fraunhofer regime 
or close to it, which means that the data is the discrete Fourier transform of the 
reflectivities (up to scaling and an overall phase) and the image is the discrete inverse Fourier transform of the data. With this approach 
we reconstruct the missing phases with acceptable accuracy if there is enough prior information 
such as support, positivity, symmetries, etc. 
This is the preferred approach when (i) the ambient medium is inhomogeneous and even if we do have the phases, they are 
randomized and cannot be used in (coherent) imaging, and (ii) we do not have the possibility of illumination diversity 
so that a holographic method can be used, assuming that the missing phases are coherent.
It is the holographic approach that we address here, when there is coherence 
in the phases and illumination diversity is available.

\subsection*{Imaging with MUSIC}
\label{sec:music}
We discuss briefly
the MUltiple SIgnal Classification method (MUSIC), which is a subspace projection algorithm that uses the 
singular value decomposition (SVD) of the full data array response matrix $\Pm(\omega)$ to form the images. 
It is an algorithm that is widely used to image the locations of $M<N$ point scatterers in a region of interest, restricted to an image window IW. 
Once the locations are known, their reflectivities can be found from the recorded intensities using convex optimization (see, e.g.~\cite{Chai16}).
Note that MUSIC can also be applied when the single frequency interferometric data matrix $\Mm(\omega)$ is available instead of $\Pm(\omega)$. 
Therefore, MUSIC can be used with intensity-only measurements \cite{Novikov14,Moscoso16}. 

We first consider MUSIC for a single frequency using the time reversal matrix $\Mm(\omega)$. 
This matrix maps illumination vectors to themselves and its eigenvectors 
$ V_j(\omega)$, $j=1\dots,M$, corresponding to non-zero eigenvalues are illumination vectors that beamform to the scatterers. They form the signal subspace. 
The remaining eigenvectors 
${V}_j(\omega)$, $j=M+1,\ldots,N$, span a subspace called the noise subspace. 
The beamforming vectors $\vect g_0(\vect y^s,\omega)$, defined by (\ref{GreenFuncVec}) with the homogeneous Green's function $G_0(\vect x_r,\vect y^s, \omega)$, will be approximately orthogonal to the noise subspace
only when $\vect y^s$  is close to a scatterer location $\vect y_j$. In this case, $\sum_{j=M+1}^{N} |\vect g_0^T(\vect y^s,\omega) V_j(\omega)|$ is close to zero and 
it follows that the scatterer locations must correspond to the peaks of the functional
\begin{equation}
\label{MUSIC_0}
\dsp \mathcal{I}(\vect y^s)=\frac{1}{\sum_{j=M+1}^{N} |\vect g_0^T(\vect y^s,\omega) V_j(\omega) |^2 }.  
\end{equation}
Often the number of scatterers is small so that the dimension of the signal subspace is much smaller than that of the noise subspace. We, therefore,  use the (normalized) equivalent functional
\begin{equation}
\label{MUSIC}
\mathcal{I}^{MUSIC}(\vect y^s)=\frac{\min_{1\le j\le K}\|\mathcal{P_N}\vect g_0(\vect y_j,\omega)\|_{\ell_2}}{\|\mathcal{P_N}\vect g_0(\vect y^s,\omega)\|_{\ell_2}},\, 
\end{equation}
with the
projection onto the noise subspace defined as
\begin{equation}
\label{proyection}
\mathcal{P_N}\vect g_0(\vect y,\omega)=\vect g_0(\vect y,\omega) -\sum_{j=1}^M (\vect g_0^T(\vect y,\omega) V_j(\omega))
 V_j(\omega).
\end{equation}
We can also define the following imaging functional 
\begin{equation}
\label{SIGNAL}
\mathcal{I}^{SIGNAL}(\vect y^s)=\frac{\|\mathcal{P_S}\vect g_0(\vect y^s,\omega)\|_{\ell_2}}{\max_{1\le j\le K}{\|\mathcal{P_S}\vect g_0(\vect y_j,\omega)\|_{\ell_2}}}   ,\,
\end{equation}
with the
projection onto the signal subspace defined as
\begin{equation}
\label{projection_signal}
\mathcal{P_S}\vect g_0(\vect y,\omega)=\sum_{j=1}^M (\vect g_0^T(\vect y,\omega) V_j(\omega))
 V_j(\omega).
\end{equation}
We note that \eqref{MUSIC}  is not robust to ambient medium inhomogeneities, unless 
the array is very large~\cite{Chai16}. Generalizations of MUSIC for multiple scattering and extended scatterers have also 
been developed (see, for example, \cite{Gruber04} and \cite{hou06}).

\section{Multifrequency data models}\label{sec:multi}
Now consider the case in which signals of different frequencies can be used to probe the medium. We introduce the composite column vector of all $S$ illuminations at the different frequencies $\omega_l$, $l=1,\dots,S$, 
\begin{equation}
\vect f = [\vect f(\omega_1)^t,  \vect f(\omega_2)^t,\dots,\vect f(\omega_S)^t]^t\,,
\end{equation}
whose dimension is $N\cdot S$,
and the full response matrix for multiple frequencies 
\begin{equation}\label{Pmultiple}
\Pm = [\Pm(\omega_1), \Pm(\omega_2),\dots,\Pm(\omega_S)]\, ,
\end{equation}
whose dimension is $N\times (N\cdot S)$.
With this notation, given a set of, say, $\aleph$ composite vector illuminations $\{\vect f^{(j)} \}_{j=1,2,\dots, \aleph}$ at multiple frequencies, the corresponding imaging problem is to determine the location and reflectivities of the scatterers from  the multifrequency array data, vectors of dimension $N$,
\begin{equation}\label{data}
\vect b^{(j)}=\Pm \vect f^{(j)} \,\, , \quad j=1,2,\dots, \,  \aleph,
\end{equation}
recorded at the array, including phases. All the  information for imaging, including phases, is contained in the full multifrequency response model matrix $\Pm$. 
As noted, the size of this matrix is $N\times (N\cdot S)$. Assuming that reflectivities do not depend on frequency the rank of this matrix is $M\cdot S$, the number of scatterers times the number of frequencies, as can be seen from (\ref{responsematrix}). This is so for scatterers in a general configuration in the discretized image window, and with distinct frequencies. We note, however, that there are special configurations where the rank can be smaller, although this does not influence the resolution theory that assumes a generic configuration. 

 If only the multifrequency intensities $\beta^{(j)}_i = \abs{ b^{(j)}_i }^2$ are recorded at the array, $i=1,\dots,N$,  $j=1,\dots, \,  \aleph$, the imaging problem is to determine the location and reflectivities of the scatterers from the absolute values of each component in \eqref{data}, i.e., 
\begin{equation}\label{dataI}
\vect \beta^{(j)}=\diag \left( \left(\vect P \vect f^{(j)} \right) ^* \left( \vect \Pm \vect f^{(j)} \right) \right)
=\diag \left( \left(\vect f^{(j)} \right)^* \vect \Pm^* \vect \Pm \vect f^{(j)} \right)
\,\, , \quad j=1,2,\dots \, ,\aleph.
\end{equation}
Because the {\em multifrequency interferometric data matrix (MFIDM)} 
\begin{equation}
\label{MFIDM}
\Mm= \Pm^* \Pm
\end{equation}
is involved in \eqref{dataI}  we will use it directly for imaging as in \cite{Novikov14}. 
Here the size of $\Mm$ is $N\cdot S \times N\cdot S$, and $\Pm$ is given by (\ref{Pmultiple}).

The main result of this paper 
is that the matrix $\Mm$  can be obtained from intensity-only measurements using an appropriate illumination strategy and the polarization identity
\begin{equation}
\label{eq:polariden}
2 < \vect x, \vect y> = \|\vect x + \vect y\|^2 - \|\vect x\|^2 - \|\vect y\|^2  + i (\|\vect x - i \vect y\|^2  - \|\vect x\|^2 -  \|\vect y\|^2 )\,\, .
\end{equation}
 The polarization identity \eqref{eq:polariden} allows us to find the inner product between two signals, and hence its phase differences, from linear combinations of the 
 magnitudes (squared) of these signals. In the next section we show how phase information can be recovered using illumination diversity and the polarization identity \eqref{eq:polariden}, that is, how to recover $\Mm$ in (\ref{MFIDM}) from intensity-only measurements. We then show how to image with this information, as already outlined in the introduction.

\section{Illumination strategy for holographic imaging}
\label{sec:alexei}
In~\cite{Novikov14,Moscoso16}, it was shown that the single frequency interferometric data matrix
$\Mm(\omega) = {\Pm(\omega)}^*  \Pm(\omega) $, where $\Pm(\omega)= [ P({\vect x}_r,{\vect x}_s;\omega)]_{r,s=1}^{N}$ is the full array response matrix,
can be recovered when signals of 
the same frequency are used for illuminations and only the intensities are measured at the receivers.
This is equivalent to recovering the inverse Fourier transform of $\Mm(\omega)$, which is the data cross-correlation matrix, when only the intensities are measured at the receivers. We therefore recover all phase differences between the elements of the array response matrix $\Pm(\omega)$ using suitable illumination diversity and measuring only intensities.
 Next, we consider a generalization of this methodology and show how 
the multifrequency interferometric data matrix $\Mm$ defined by (\ref{MFIDM}) can be obtained from intensity-only measurements with suitable illuminations.

\subsection*{General case}
We consider first the general case in which sources and receivers are not necessarily colocated, that is, they are not placed at the same positions, and we recover the elements of the MFIDM, $\Mm$, associated with one receiver at location $\vect x_r$, which we denote by $\Mm_r$,
using a suitable illumination strategy. We describe first the structure of $\Mm_r$.
We note that the $r$th row of the multifrequency array response matrix $\Pm$, given by (\ref{Pmultiple}), has the form
\begin{equation}
\label{eq:Pr}
\Pm_r=(p_{r 1}, p_{r 2}, \dots, p_{r N\cdot S})
\end{equation}
where the entry $p_{r j}$, with $j= s+ (l-1)\cdot N$, denotes the received signal at $\vect x_r$ when 
the source at ${\vect x}_s$ sends a  signal with frequency $\omega_l$.
With this notation, $\Mm_r$ is the rank-one matrix 
\begin{equation}
\label{eq:Mr}
\Mm_r=\Pm^*_r \Pm_r\, .
\end{equation}

Let $\vect\we_{i} = [0, 0, . . . , 1, 0, . . . , 0]^T$ with $i=i(s,l)$
be the illumination vector representing a signal of magnitude one and  frequency $\omega_l$ sent from the source ${\vect x}_s$, and where $i= s+ (l-1)\cdot N$. 
The $(i,j)$ entry of $\Mm_r$, which is obtained from intensity measurements at  ${\vect x}_r$ when the illuminations $\vect\we_{i}$ and $\vect\we_{j}$ 
are used is  
\begin{equation}
\label{eq:Mr_entry}
m^r_{ij} = \overline{p}_{r i} p_{r j} = \left(\Pm_r\vect\we_{i}  \right)^*\Pm_r\vect\we_{j}\, .
\end{equation}
The key point here is that  \eqref{eq:Mr_entry} can be obtained from intensity-only measurements using
the polarization identity. Indeed, the polarization identity~\eqref{eq:polariden} gives us
 \begin{equation}\label{polarization_re}
\mbox{Re}(m^r_{ij})= 
\frac{1}{2} \left( \norm{\Pm_r\vect\we_{i+j}}^2  - 
\norm{\Pm_r\vect\we_i}^2 - \norm{\Pm_r \vect\we_j}^2\right)
\end{equation}
\begin{equation}\label{polarization_im}
\mbox{Im}(m^r_{ij})= 
 \frac{1}{2} \left( \norm{\Pm_r\vect\we_{i-{\bf i} j}}^2  -
  \norm{\Pm_r \vect\we_i}^2 - \norm{\Pm_r \vect\we_j}^2\right), 
  \end{equation}
where $\vect\we_{i+j}= \vect\we_{i}+\vect\we_{j},  \quad \vect\we_{i-{\bf i} j}= \vect\we_{i}-{\bf i} \vect\we_{j}$, where ${\bf i}= \sqrt{-1}$. 
In \eqref{polarization_re} and \eqref{polarization_im}, $\mbox{Re}(\cdot)$ and $\mbox{Im}(\cdot)$ denote the real and imaginary parts
of a complex number, respectively. Since all entries on the right-hand side of~\eqref{polarization_re} and~\eqref{polarization_im}  involve intensity-only measurements on the 
$r$th receiver, we can recover all the entries $m^r_{i j}$ in $\Mm_r$.

\subsection*{Symmetric case}
In the general case we recovered, from intensity-only measurements at a single receiver $\vect x_r$, the rank one matrix $\Mm_r$
resulting from signals sent from different sources and with different frequencies. Sources and receivers do not have to be
colocated, but we cannot obtain all the elements of $\Mm$ this way. In other words, we cannot recover MFIDM $\Mm$ from the set of all matrices $\Mm_r$. 
This is so, because each rank one matrix $\Mm_r$ is obtained up to a global phase that will be different for each
receiver location $\vect x_r$. We now  show that when sources and receivers are placed at the same positions, the full MFIDM $\Mm$ can
be recovered from intensity-only measurements.

For colocated sources and receivers, the full array response matrix is symmetric for each frequency-block $\Pm(\omega_l)$.
Let $p_{i j}^l$ represents the signal measured at  receiver $\vect x_i$ due to a signal of frequency $\omega_l$ sent from source $\vect x_j$. Then,
\begin{equation}\label{eq:mig}
p_{i j}^l=p_{j i}^l
\end{equation}
because of wave field reciprocity.
Let  $c_{j} = j+ (l-1)\cdot N$ and $c_{i} = i+ (l-1)\cdot N$, then for the full response matrix $\Pm$, given by (\ref{Pmultiple}), the identity~\eqref{eq:mig} 
becomes
$p_{i c_{j}}=p_{j c_{i}}$. Using these symmetries and the notation (\ref{eq:Mr_entry})
 we can recover all products
 \begin{equation}\label{eq:cross_cor}
\overline{p}_{i k} p_{j n} = \frac{\overline{p}_{i k} p_{i 1} \overline{p}_{j 1} p_{j n}}{  p_{1 c_i} \overline{p}_{1 c_j}} = \frac{m^i_{k 1} m^j_{1 n}}{m^1_{c_i c_j}}
\end{equation}
for different sources, frequencies, and  receivers and, therefore, we can recover the full MFIDM $\Mm$. 

\section{Interferometric imaging}
\label{sec:imaginghomo}
We present in the next section a direct method to form the images from intensity-only measurements at a single receiver with range (depth) and cross range resolution equivalent to those 
that are obtained by migrating the full response matrix including phases.

\subsection*{Single receiver interferometry imaging (SRINT)}

We write the row-vector $\Pm_r \in\mC^{(N\cdot S)}$, defined in \eqref{eq:Pr},  in the form
\begin{equation}
\label{eq:pmr}
\Pm_r =  \bfrho^t \vect\Gc_r \,,
\end{equation}
Here, $\vect\Gc_r^t$ is the $(N\cdot S)\times K$ model matrix in a homogeneous or heterogeneous medium that maps a distribution of scatterers in the 
IW to the data  received at the array, i.e.,
\begin{equation}
\label{eq:G}
 \vect\Gc_r^t = 
\begin{bmatrix}
G(\vect x_r, \vect y_1;\omega_1) \vect g(\vect y_1;\omega_1)\,\cdots\,G(\vect x_r ,\vect y_K;\omega_1)\vect g(\vect y_K;\omega_1) \\
G(\vect x_r, \vect y_1;\omega_2)\vect g(\vect y_1;\omega_2)\,\cdots\,G(\vect x_r ,\vect y_K;\omega_2)\vect g(\vect y_K;\omega_2) \\
\hdotsfor{1} \\
G(\vect x_r, \vect y_1;\omega_S)\vect g(\vect y_1;\omega_S)\,\cdots\,G(\vect x_r, \vect y_K;\omega_S)\vect g(\vect y_K;\omega_S) 
\end{bmatrix}.
\end{equation}
In this expression, $G(\vect x_r ,\vect y_j;\omega_l)$, $j=1,\ldots,K$ and $l=1,\ldots,S$, denotes the Green's function with source at grid point $\vect y_j$, receiver at  $\vect x_r$ and frequency $\omega_l$, and $\vect g(\vect y_i;\omega_l)$ defined by (\ref{GreenFuncVec}) is the vector of illuminations from the array to grid point $\vect y_i$ in the image window.

To form the image given the data $m_{ij}^r$ \eqref{eq:Mr_entry} corresponding to $d((\vect x_r, \vect x_{r}),(\vect x_s, \vect x_{s'}),(\omega_l,\omega_{l'}))$ for all source locations $s,s'=1,\ldots,N$ and frequencies $l,l'=1,\ldots,S$, where $ i = s + (l-1) \cdot N$ and $ j = s' + (l'-1) \cdot  N$,
we compute the imaging functional \eqref{eq:INTERF} which for the single receiver case can be re-written in the following matrix form
\begin{equation}
\label{eq:interfR}
{\cal I}^{Interf}({\vect y}^s) =  \diag( \vect\Gc_{0r} \Mm_r\vect\Gc_{0r}^*)({\vect y}^s) \, ,
\end{equation}
where $\vect\Gc_{0r}$ is the model matrix in a homogeneous medium. The imaging functional \eqref{eq:interfR} can be viewed as migrating the multifrequency interferometric data matrix. In fact, \eqref{eq:interfR} corresponds to coherent interferometric imaging (CINT)  for one receiver. More precisely, since we have not yet introduced any thresholding, \eqref{eq:interfR} is just the usual Kirchhoff migration functional squared
as explained in the introduction (see \eqref{eq:INTERF} and \eqref{eq:KM}). 
\subsection*{Resolution}
To understand the resolution of the imaging functional \eqref{eq:interfR}, we  consider 
scatterers in a homogeneous medium and we substitute \eqref{eq:pmr} into \eqref{eq:interfR} 
\begin{equation}
\label{eq:AMAexp}
{\cal I}^{Interf}({\vect y}^s) = \diag( \vect\Gc_{0r}\vect\Gc_{0r}^* \, \overline{\bfrho} \bfrho^t \, \vect\Gc_{0r} \vect\Gc_{0r}^*) \, ,
\end{equation}
which shows that, in a homogeneous medium, \eqref{eq:interfR} produces sharp images if $\vect\Gc_{0r} \vect\Gc_{0r}^*$, is close to a diagonal matrix, that is,
if the columns of $\vect\Gc_{0r}$ are nearly orthogonal. The near orthogonality of $\vect\Gc_{0r}$ is satisfied when the discretization of the IW is 
compatible with the  resolution provided by the bandwidth and the array size, $c/B$ in range or depth, and  $ \lambda_0 L/a$ in cross-range. 
In theory this is the case when the distance
between two adjacent grid points $\vect y_k$ and $\vect y_{k'}$ is such that
\begin{equation}
 |\vect y_k - \vect y_{k'}| \gg \max\left\{\frac{ \lambda_0 L }{a}, \frac{c}{B}\right\} \, .
\end{equation}
In practice we can have a much finer discretization in the image window if we use an a-posteriori thresholding to select image peaks.

We note that for reasonable discretizations of the IW the imaging functional \eqref{eq:AMAexp} (and, therefore, \eqref{eq:interfR}) produces an 
image of $\diag(\overline{\bfrho} \bfrho^t)$, which is just the complex conjugate of the diagonal entries of the unknown matrix $\bfrho \bfrho^*$ in 
\eqref{eq:non-linear} that is found in \cite{CMP11} through a nuclear norm minimization process.

\subsection*{Thresholding and Masks}
\label{sec:masks}
Thresholding was introduced in interferometric imaging in  \cite{BPT-05,BPT-ADA}  so that frequency-offsets, and source and receiver location offsets are 
restricted to within coherence limits. The resulting coherent interferometric (CINT) functional has the form (\ref{eq:CINT}). There, it is assumed that the full array response matrix 
$\Pm(\omega_l)= [ P({\vect x}_r,{\vect x}_s;\omega_l)]_{r,s=1}^{N}, ~l=1,2,\ldots,S$ is recorded, including phases. The frequency cutoff $\Omega_d <B$ and the source or receiver cutoff $X_d <a$ are not know in advance, as noted in the introduction, but can be determined in the image formation process. This thresholding removes noise and stabilizes the image at the expense of somewhat reduced resolution \cite{BGPT-rtt, BGPT-beamforming}.

A main result in this paper is that the multifrequency interferometric data matrix $\Mm_r$ in (\ref{eq:Mr}) can be obtained from intensity-only measurements at a single receiver, as explained is Section \ref{sec:alexei}, and robust imaging can be done interferometrically with SRINT as given by (\ref{eq:SRINT}).

In the framework of SRINT the thresholding can be easily done by multiplying (in the sense of element-wise multiplication) the matrix $\Mm_r$ in 
\eqref{eq:interfR} by a mask, that is a matrix that has only zeros and ones depending on the spacing between the indices in the matrix so as to restrict 
the data used in imaging. This is an efficient way of implementing CINT for one receiver, that is, SRINT.
Incorporating the mask in \eqref{eq:interfR} leads to the following matrix form of the SRINT imaging functional
\begin{equation}
\label{eq:AMA}
{\cal I}^{SRINT}({\vect y}^s) = diag( \vect\Gc_{0r} \vect{\Zc}\odot\Mm_r\vect\Gc_{0r}^*)({\vect y}^s) \, ,
\end{equation}
where $\vect{\Zc}$ denotes the mask. The product $\vect{\Zc}\odot\Mm_r$ denotes element-wise multiplication. 

We put masks on the data to reduce anticipated decoherence of measurements that arise if sources and/or receivers are far apart. 
Wide bandwidth also can lead to decoherence of measurements, and masks that limit the frequency offset 
are also needed. The need for masks can be understood as follows. If we send two signals of nearby frequencies $\omega_l$ and $\omega_{l'}$ 
from nearby sources at $\vect x_s$ and $\vect x_{s'}$ respectively, then they travel through essentially the same medium and will be affected 
in a similar way by the random inhomogeneities. 
This is quantified by the distance between sensors on the array being smaller than the decoherence distance so $|\vect x_s-\vect x_{s'}| < X_d$, 
and the frequency offset being smaller then the decoherence frequency so $| \omega_l -\omega_{l'}| < \Omega_d$. 
Thus, a mask is a matrix composed by zeros and ones that restrict data only to the coherent nearby source locations and frequencies.

\section{Numerical Simulations}
\label{sec:numerics}

In this section, we present numerical simulations that illustrate the performance of the proposed interferometric method. All length scales are 
measured in units of the central wavelength $ \lambda_0$. 

We consider a regime in optical (digital) microscopy with a central frequency of $f_0=600$ THz which corresponds to $ \lambda_0=500 nm$. We will assume that we can make measurements for multiple 
frequencies covering a total bandwidth of $120$THz. All wavelengths considered here are in the visible spectrum and range from blue to green light. 
The size of the array that we use for imaging is $a=500 \lambda_0$ { and has $N=81$ equispaced transducers}. The distance from the 
array to the image window IW is $L=10000 \lambda_0$. The size of the IW is $160  \lambda_0 \times 80  \lambda_0$ and the pixel size is 
$2 \lambda_0\times  \lambda_0$ in cross-range and range (depth), respectively. In all the figures that follow, the true locations of the point scatterers are 
indicated with white crosses. We note that these specifications do not correspond to any specific device. They are broadly compatible with current spatial light modulator technologies. 

\subsection{Robustness to image window discretization with a full array and full phase information}

We assume here that the phases can be recorded at the array and that the medium between the array and the IW is homogeneous.
In the single frequency case, signals are emitted from all 
the elements in the array, one at a time, and the reflections are recorded at all of them as well. Therefore, the data available for imaging is 
the full, single frequency $N\times N$ array response matrix $\Pm(\omega)$ defined in \eqref{responsematrix}.

To assess the robustness of the different imaging methods with respect to the discretization of the IW we consider two configurations: one with the scatterers placed on the
grid and a second one where the scatterers are displaced with respect to the grid. More precisely, 
the off-grid scatterers are displaced by half the grid size in both directions from a grid point. 

\begin{figure}[htbp]
\begin{center}
\begin{tabular}{ccc}
\includegraphics[scale=0.22]{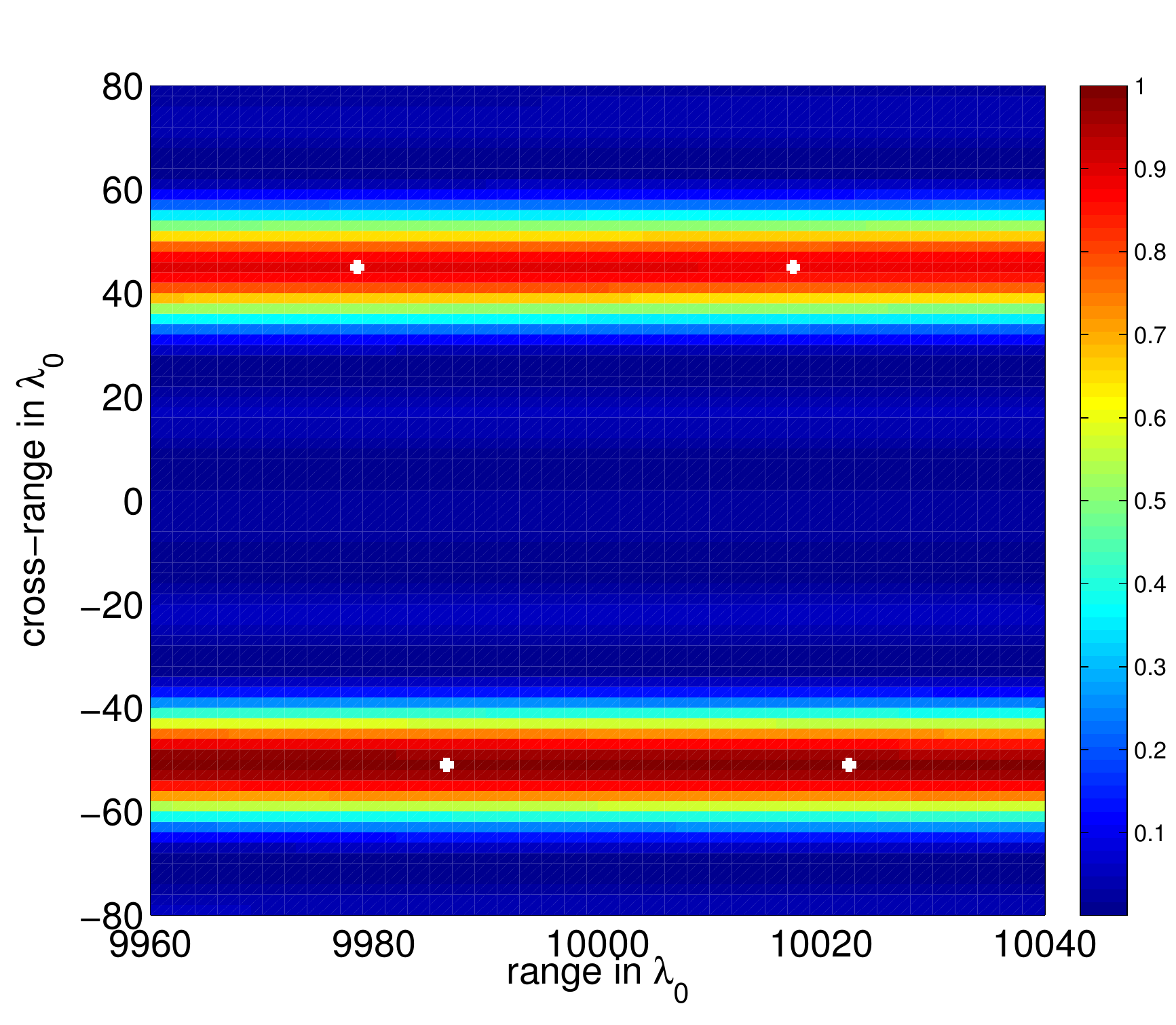}&
\includegraphics[scale=0.22]{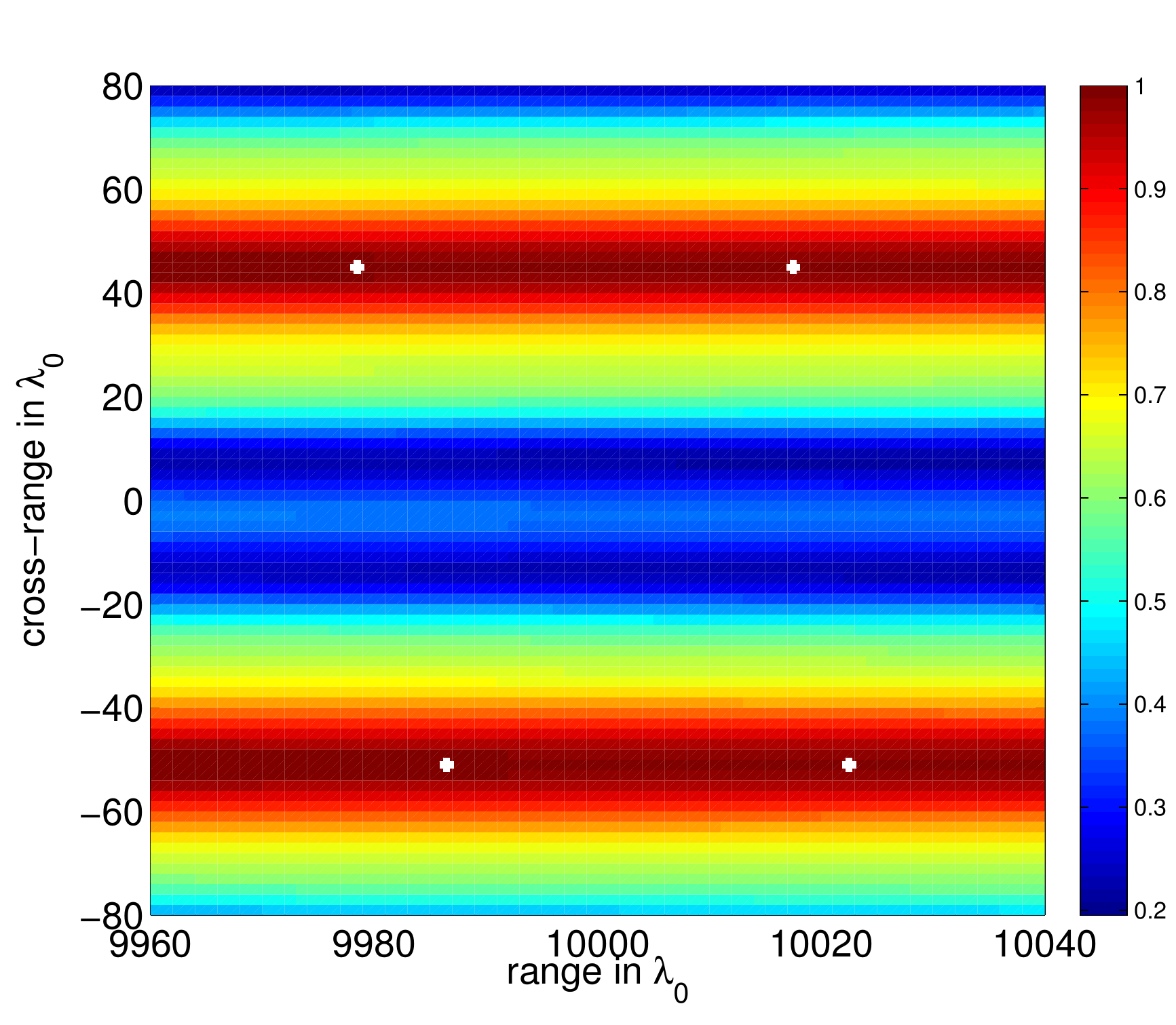}&
\includegraphics[scale=0.22]{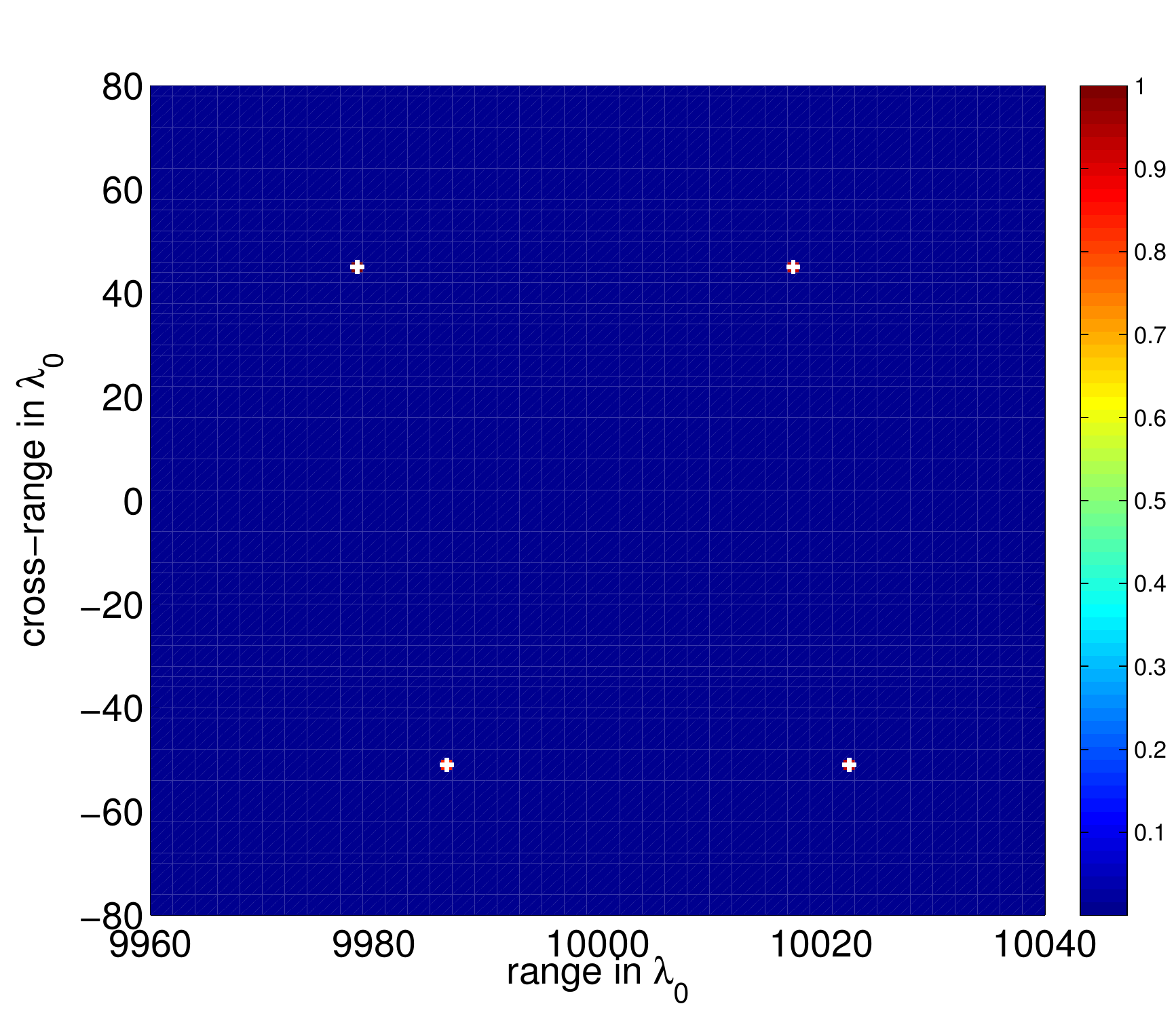}\\
\includegraphics[scale=0.22]{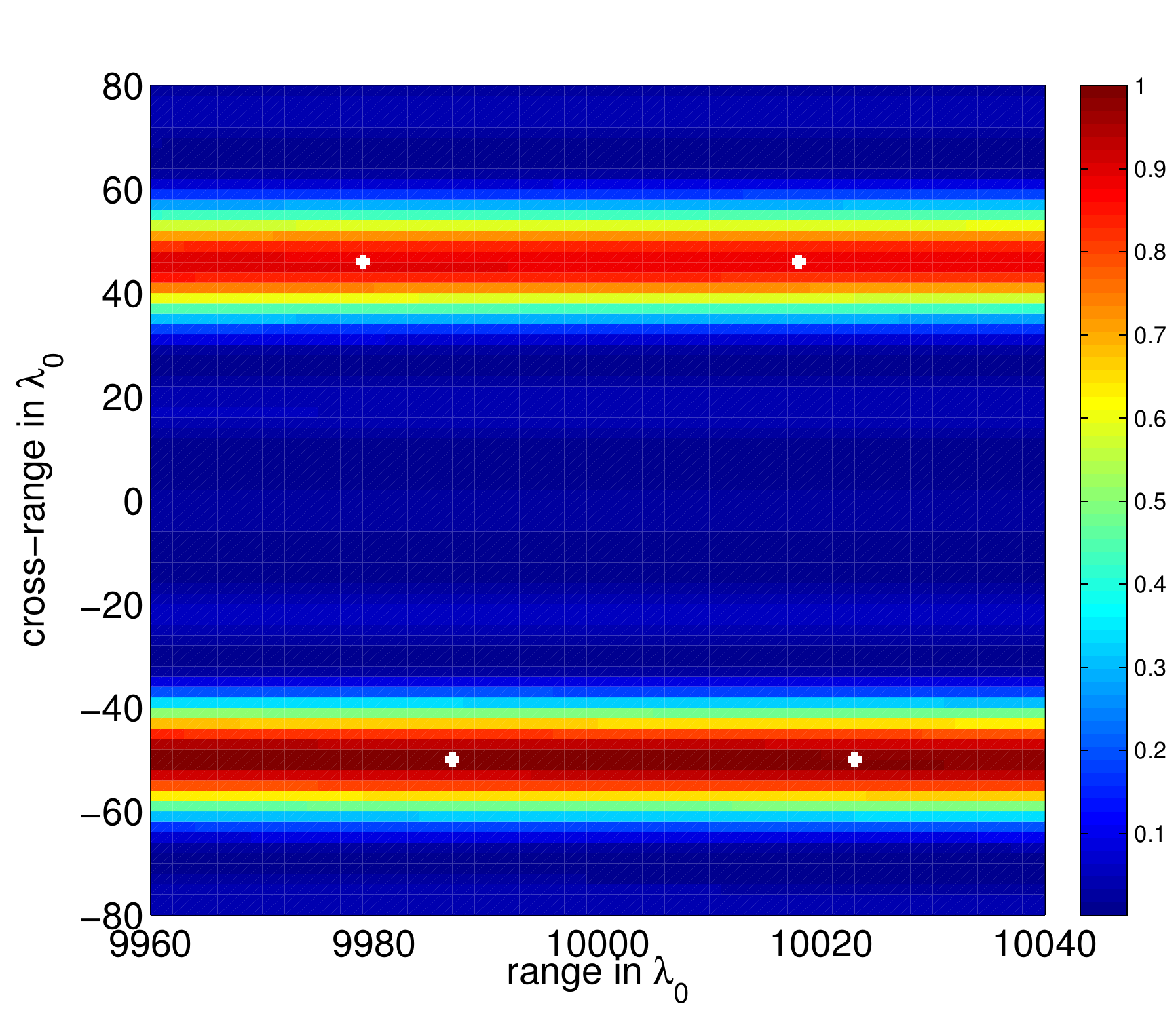}&
\includegraphics[scale=0.22]{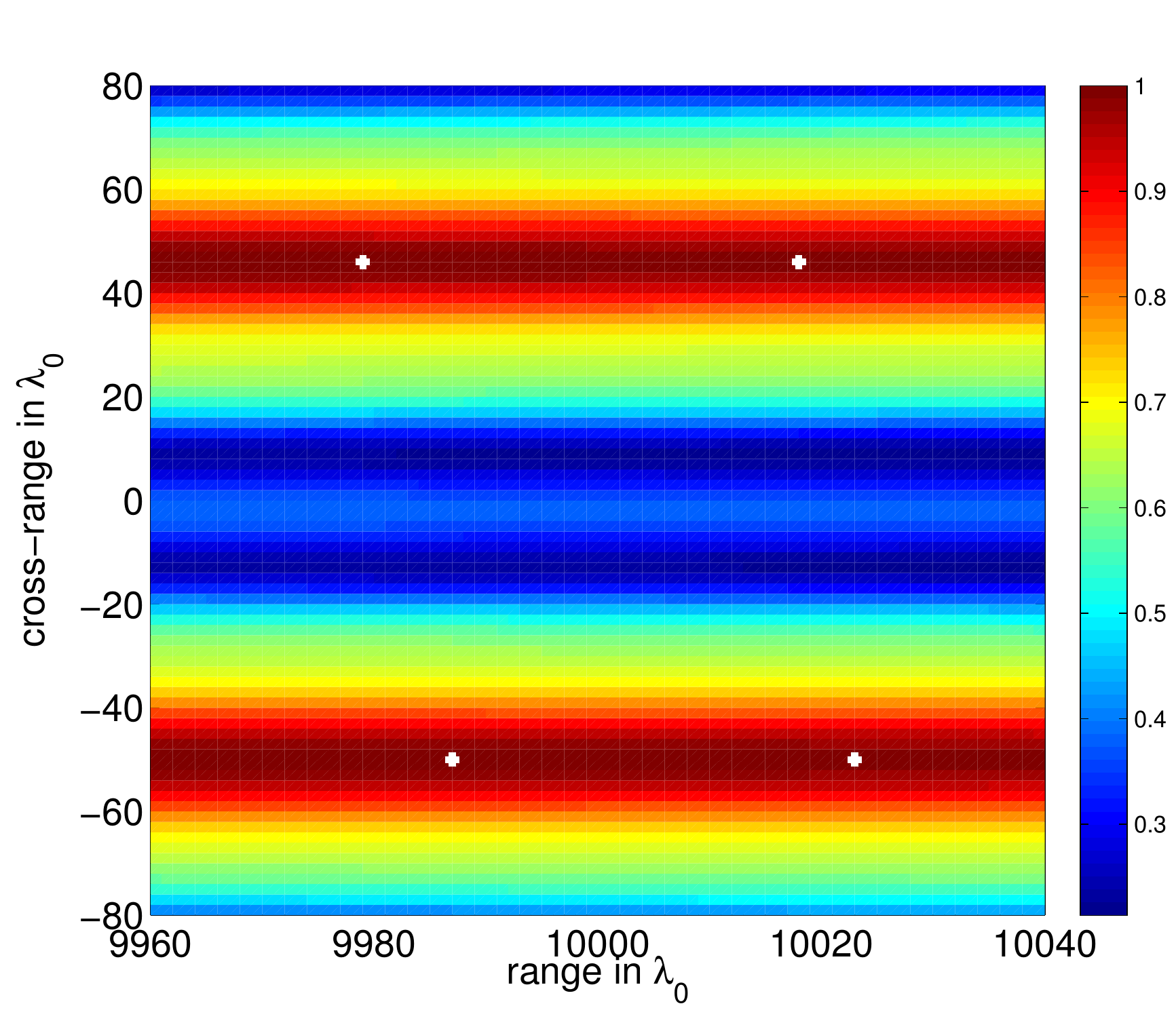}&
\includegraphics[scale=0.22]{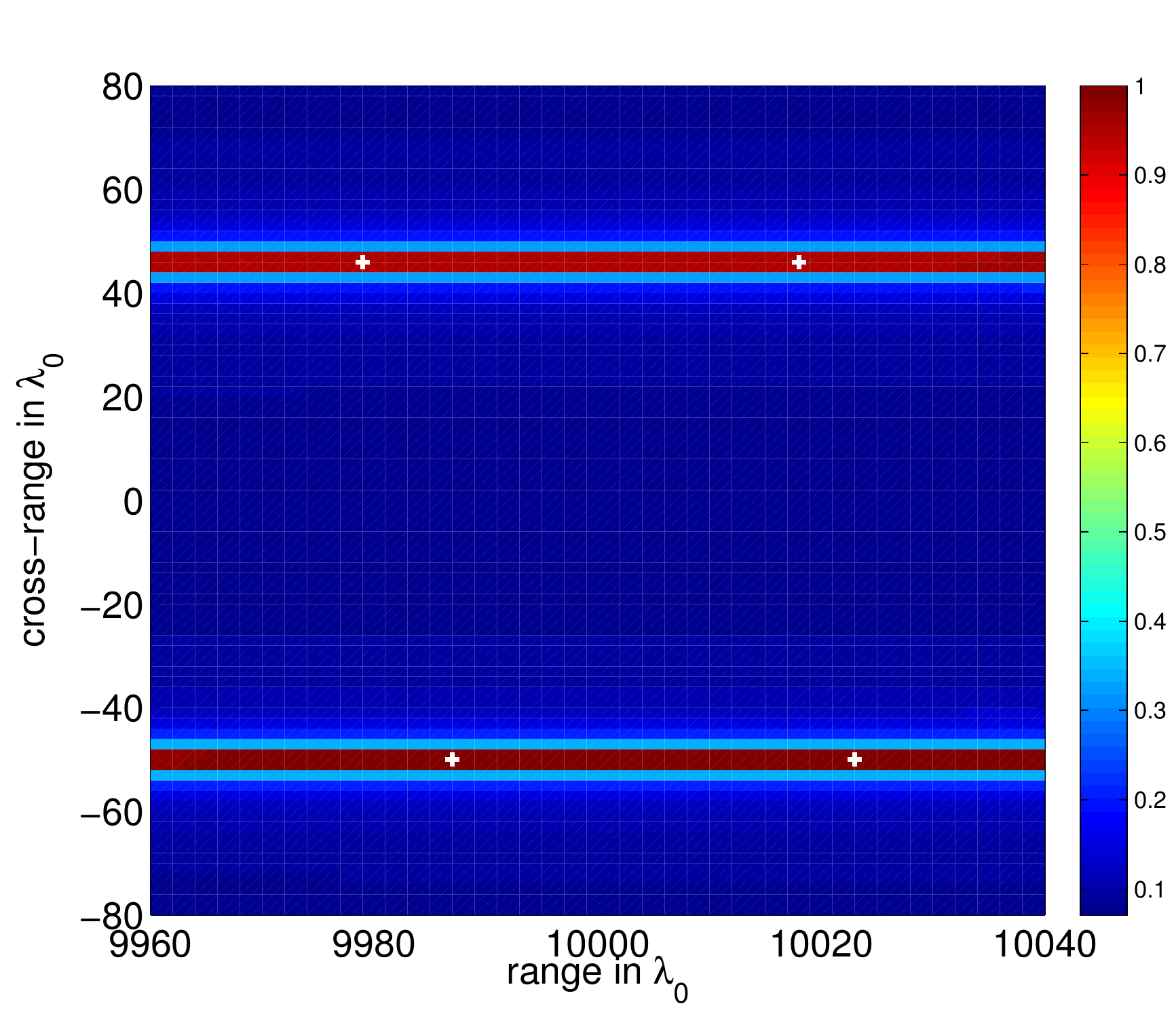}\\
\end{tabular}
\end{center}
\caption{{\bf Single frequency full data, including phases}. Homogeneous medium. Top row: Scatterers {\bf on} the grid. 
Bottom row: Scatterers {\bf off} the grid. From left to right: $\mathcal{I}^{KM}$ as defined in \eqref{eq:KM} with $S=1$ corresponding to the single frequency $f_0$, $\mathcal{I}^{SIGNAL}$ as defined in \eqref{SIGNAL} and $\mathcal{I}^{MUSIC}$ 
as defined in \eqref{MUSIC}. No additive noise in the data. }
\label{fig:h1}
\end{figure}

Figure \ref{fig:h1} shows the images obtained when there is no additive noise in the data. In the top row of this figure
 all the scatterers are placed on the imaging grid, while in the bottom row the scatterers are off-grid. Since the array size $a$ is small with respect 
to the distance $L$ ($L/a=20$) and we only have one frequency, we expect a cross-range 
resolution of $ \lambda_0 L/a=20  \lambda_0$ and a range resolution of $  \lambda_0 (L/a)^2=400  \lambda_0$.
We see this in the images shown in the left and middle columns of Figure~\ref{fig:h1} obtained with 
$\mathcal{I}^{KM}$ (left column), as defined in \eqref{eq:KM} for $S=1$ corresponding to the single frequency $f_0$, and 
$\mathcal{I}^{SIGNAL}$ (middle column), as defined in \eqref{SIGNAL}, no matter whether the scatterers are on or off the grid. 
On the other hand, $\mathcal{I}^{MUSIC}$, as defined in \eqref{MUSIC}, gives
very precise estimates of the scatterer's locations when these are on the grid, as can be seen in the top right plot of this figure. However, 
when the scatterers are off-grid the MUSIC image deteriorates dramatically as it is shown in the bottom-right plot of Figure~\ref{fig:h1}. 
These simulations illustrate clearly the lack of robustness of the MUSIC algorithm with respect to modeling errors such as
off-grid displacements.

\begin{figure}[htbp]
\begin{center}
\begin{tabular}{ccc}
\includegraphics[scale=0.22]{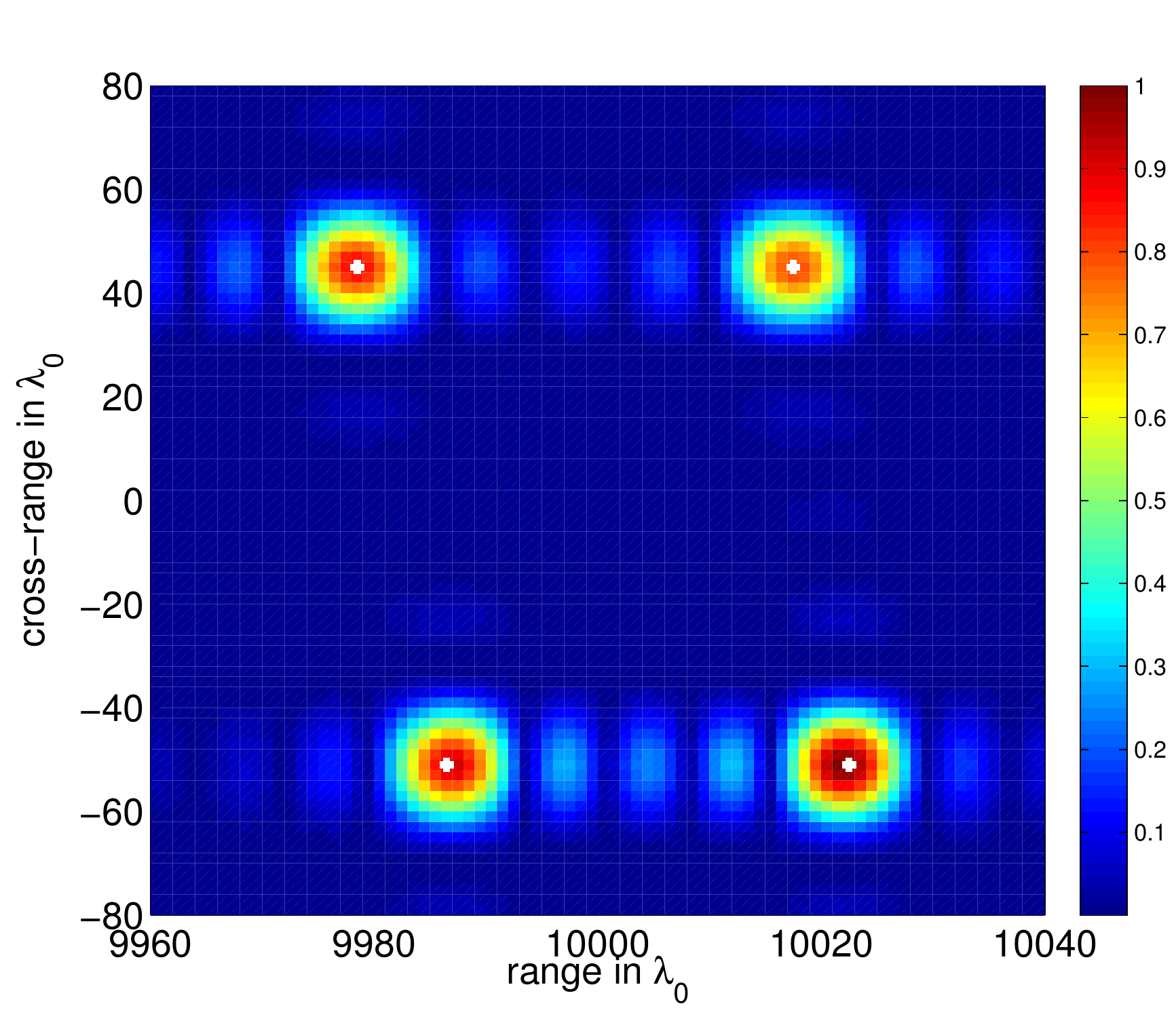}&
\includegraphics[scale=0.22]{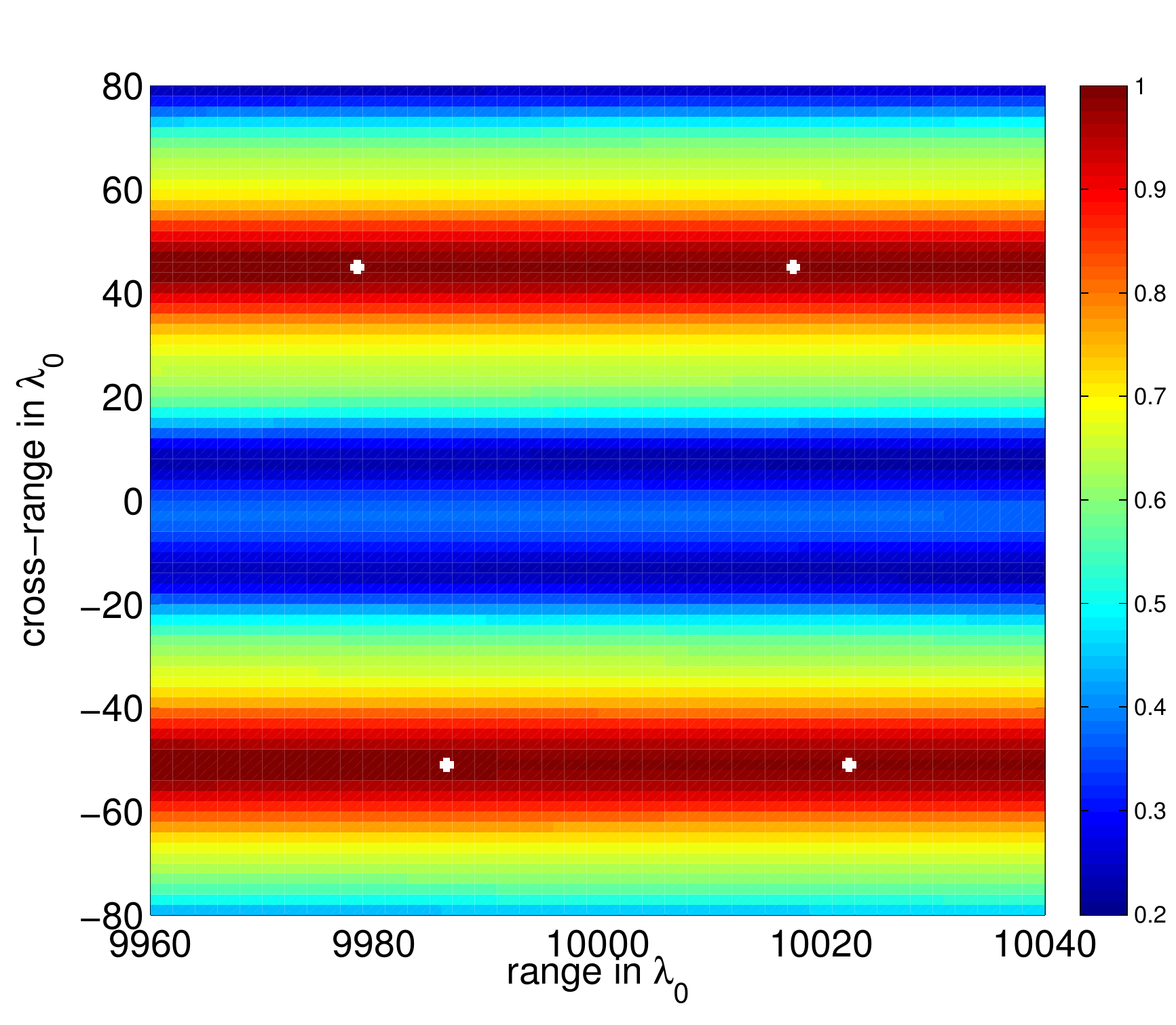}&
\includegraphics[scale=0.22]{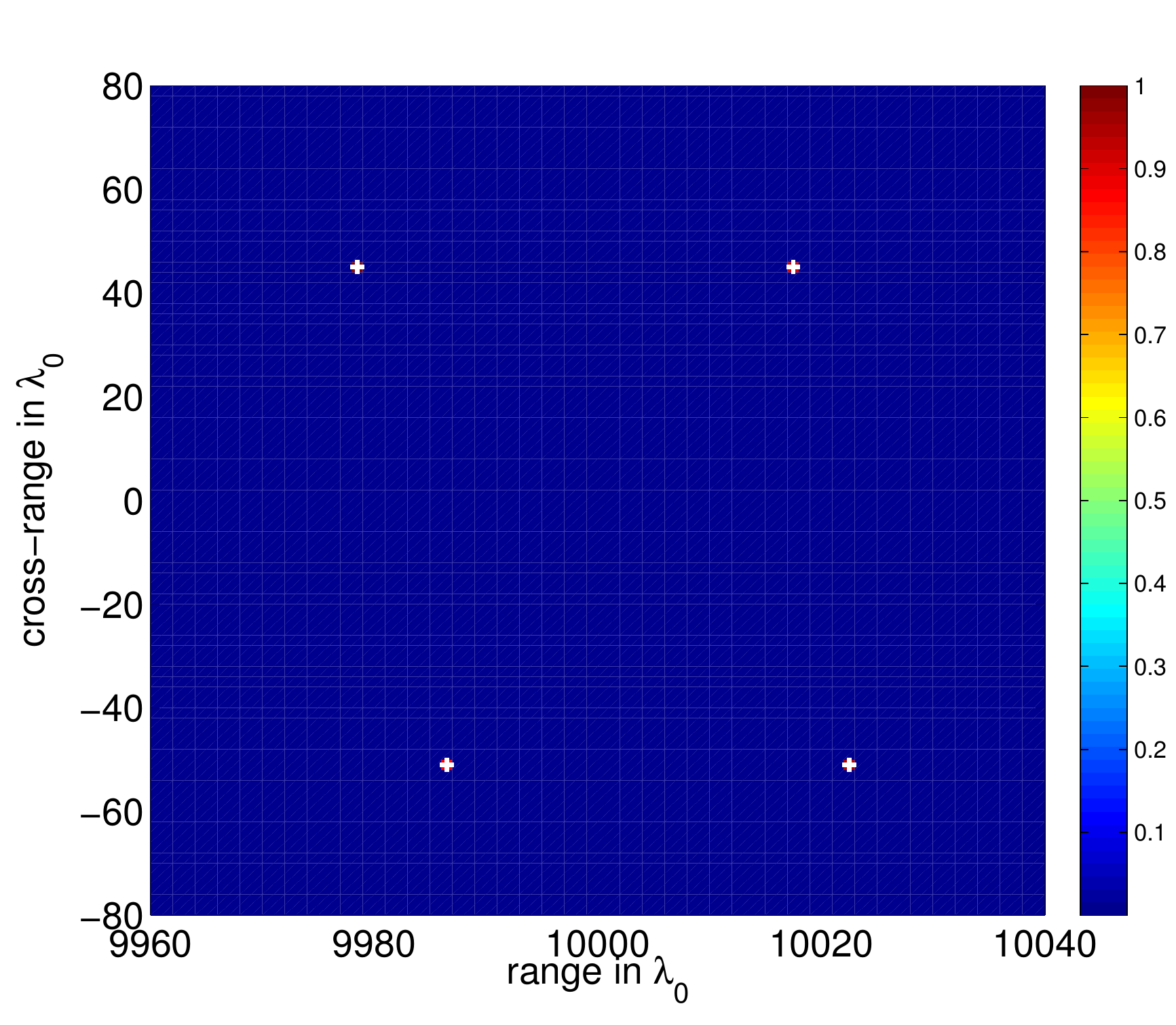}\\
\includegraphics[scale=0.22]{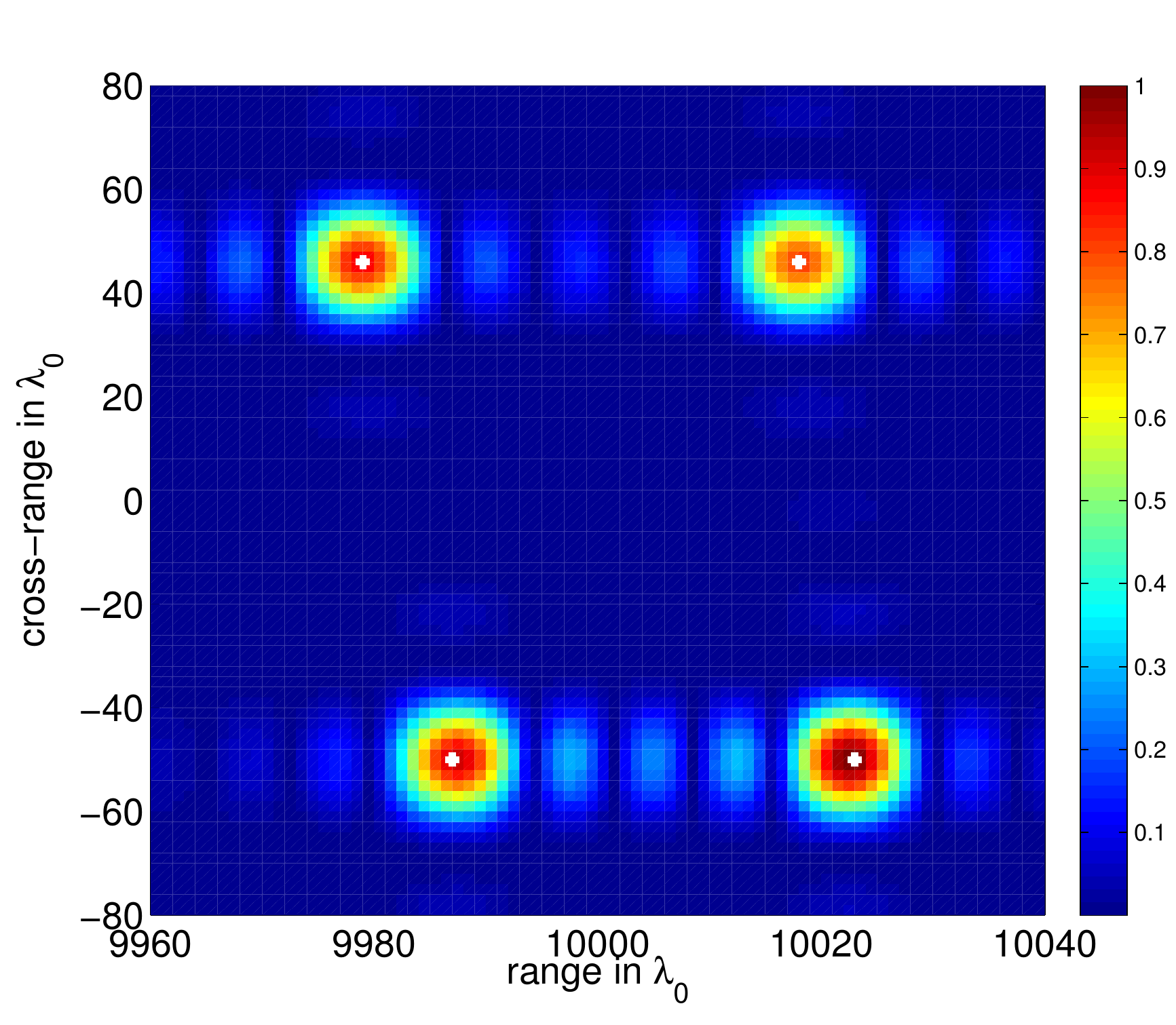}&
\includegraphics[scale=0.22]{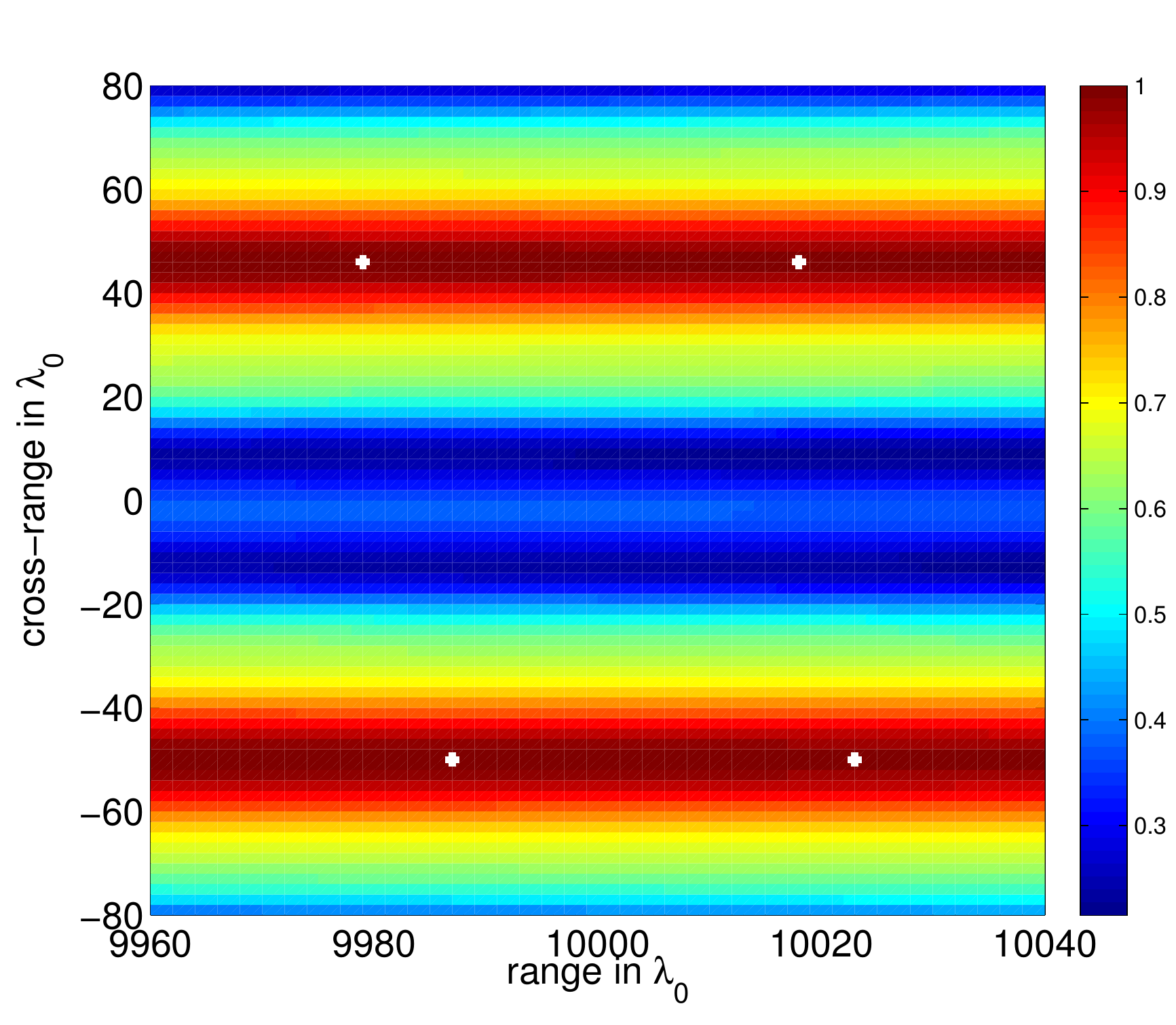}&
\includegraphics[scale=0.22]{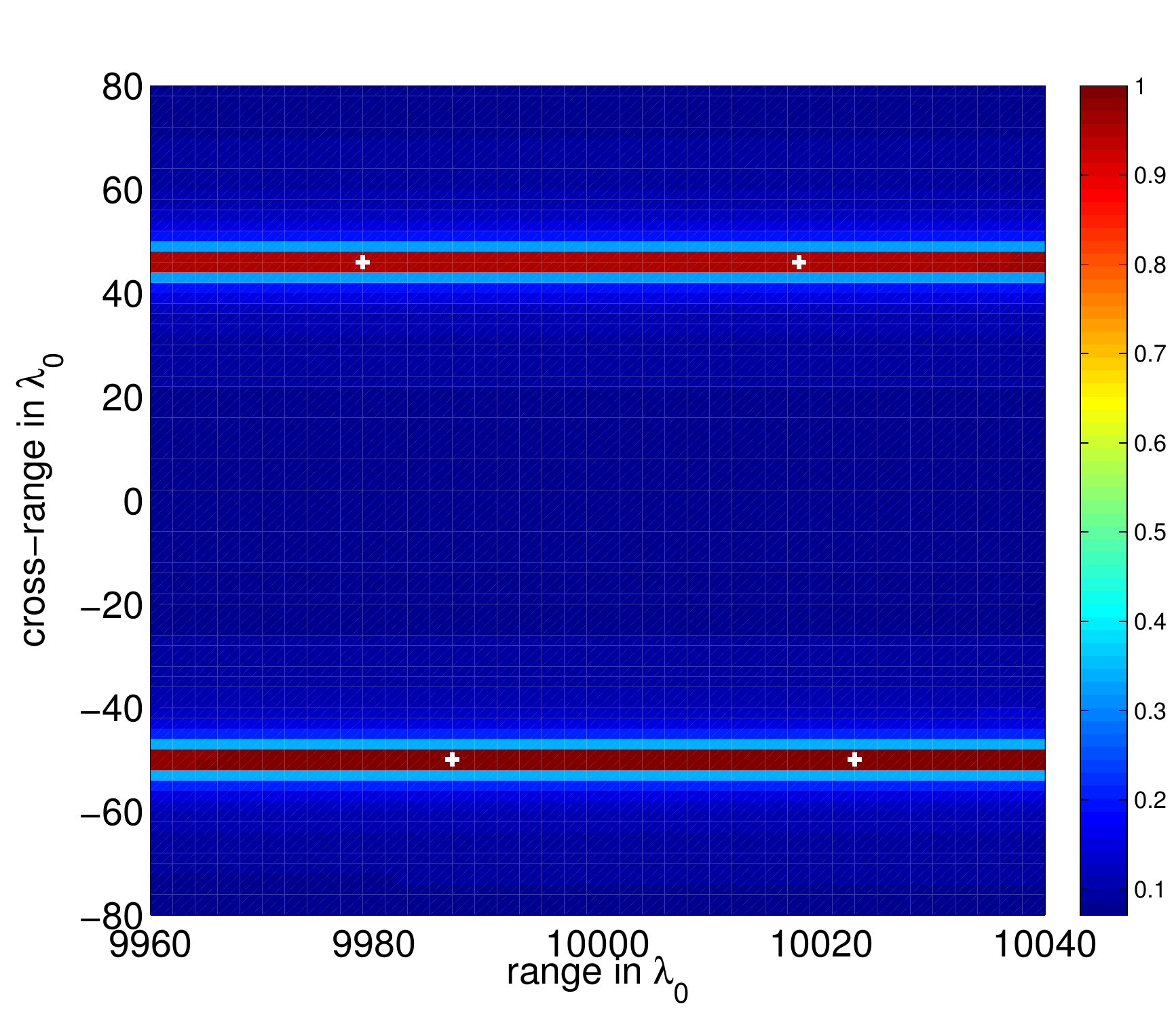}\\
\end{tabular}
\end{center}
\caption{{\bf Mutifrequency full data (including phases)}. Using 16 frequencies equally spaced in the bandwidth $[580,620]$THz. Homogeneous medium. Top row: Scatterers on the grid. Bottom row: Scatterers off the grid. 
From left to right: $\mathcal{I}^{KM}$ as defined in \eqref{eq:KM}, $\mathcal{I}^{SIGNAL}$ as defined in \eqref{SIGNAL} and $\mathcal{I}^{MUSIC}$ as defined in \eqref{MUSIC}. No additive noise in the data. }
\label{fig:h2}
\end{figure}

In a second numerical simulation, we consider probing signals with $S=16$ different frequencies equally spaced in the bandwidth $B=[580,620]$THz.
The data is now
the multiple frequency, $N\times N\cdot S$ response matrix $\Pm$ defined in \eqref{Pmultiple}. Figure~\ref{fig:h2} shows the imaging results.
As expected, MUSIC does not benefit from multiple frequencies since the projection onto the null 
space is performed frequency by frequency. In other words, conventional multiple frequency MUSIC corresponds to adding 
all the single frequency images incoherently over frequencies and, hence, the bottom right MUSIC image in Figure~\ref{fig:h2}, obtained with
16 frequencies, is  not any better than its single frequency counterpart in Figure~\ref{fig:h1}. 
On the other hand, KM performs very well when multiple frequencies are available, as it is shown in the top and bottom plots in the 
left column of Figure~\ref{fig:h2}. We note that these images remain unchanged showing the robustness of KM 
 with respect to off-grid displacements. Indeed, what matters is the point spread function of KM which determines 
 the resolution of the image, that is, $ \lambda_0 L/a=20  \lambda_0$ in cross-range and $C_0/B =  \lambda_0 f_0/B=15 \lambda_0$ 
 in range (we used $f_0/B=15$ in this case).  

While $\mathcal{I}^{SIGNAL}$ and $\mathcal{I}^{MUSIC}$ are both subspace projection algorithms, there is an important difference
between them that has considerable impact in the visualization of the images. Indeed, $\mathcal{I}^{SIGNAL}$ directly displays the norm of 
the data projected onto the signal subspace.  $\mathcal{I}^{MUSIC}$, however, first computes the norm of the data projected onto the noise subspace and then
displays its inverse, that is, it displays one over the norm of the data projected onto the noise subspace. It is the 'one over' that 
makes these two methods so different when there is neither additive noise nor modeling errors due to off-grid placements. In these cases, $\mathcal{I}^{MUSIC}$ gives exactly zero 
in the denominator at the scatterer's locations and, therefore, it acts like a sharp thresholding which, however, does not necessarily work in the 
presence of noise or modeling errors.

\subsection{Imaging with intensity-only measurements}

Now, let assume that we have only the intensities recorded at the array. 
As explained in Section \ref{sec:alexei}, if signals of multiple frequencies are used to probe the medium,
the multifrequency interferometric data matrix $\Mm_r$ can be recovered from intensity measurements at a single receiver ${\vect x}_r$
using an appropriate illumination strategy. Then, images can be formed by using the functional~\eqref{eq:AMA}.
In this subsection and in the rest of the Section, the signals used to
recover the matrix $\Mm_r$ are recorded at the receiver located at the center of the array.
\subsubsection{Robustness and resolution in homogeneous media}
\begin{figure}[htbp]
\begin{center}
\begin{tabular}{cc}
\includegraphics[scale=0.22]{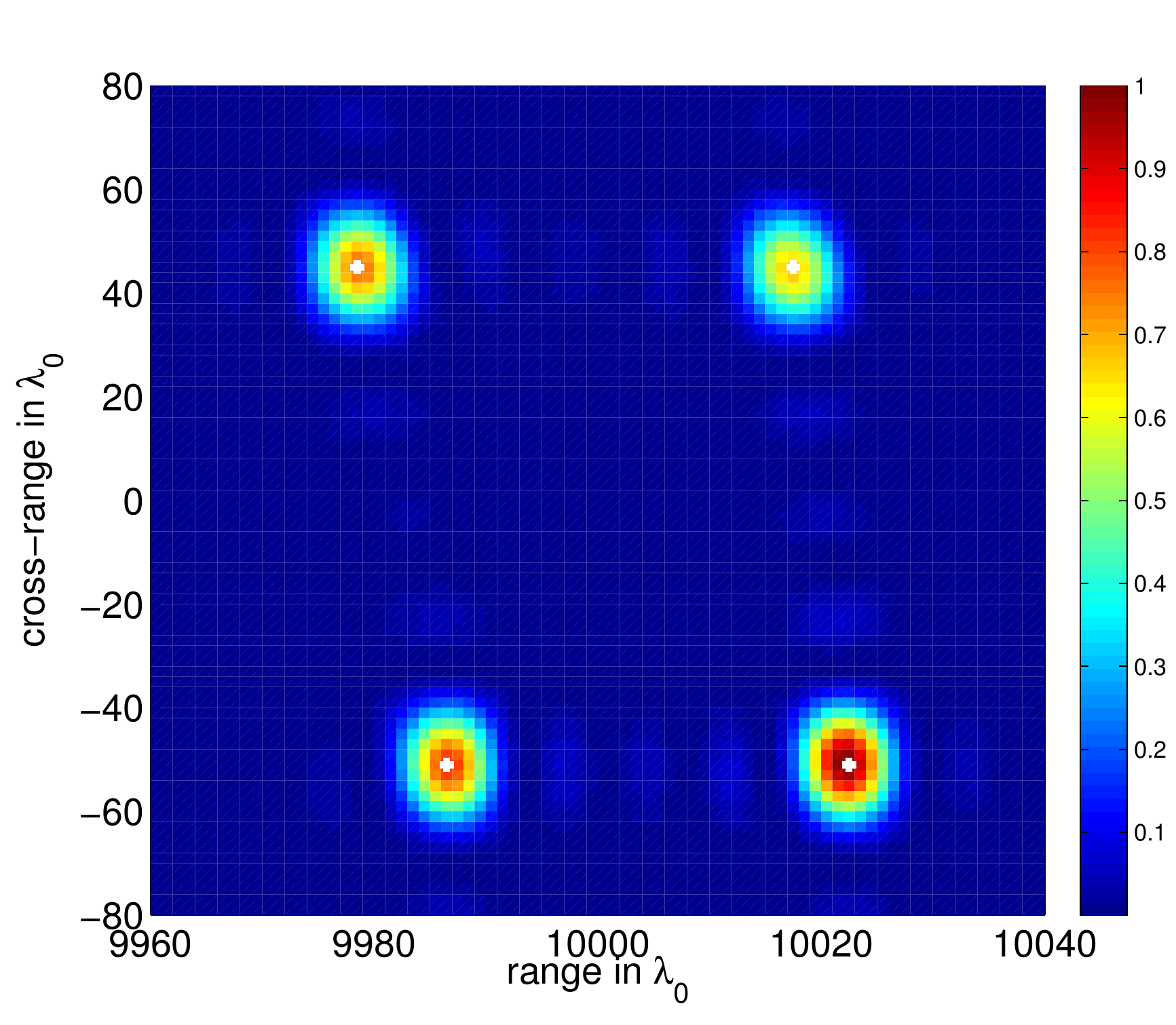}&
\includegraphics[scale=0.22]{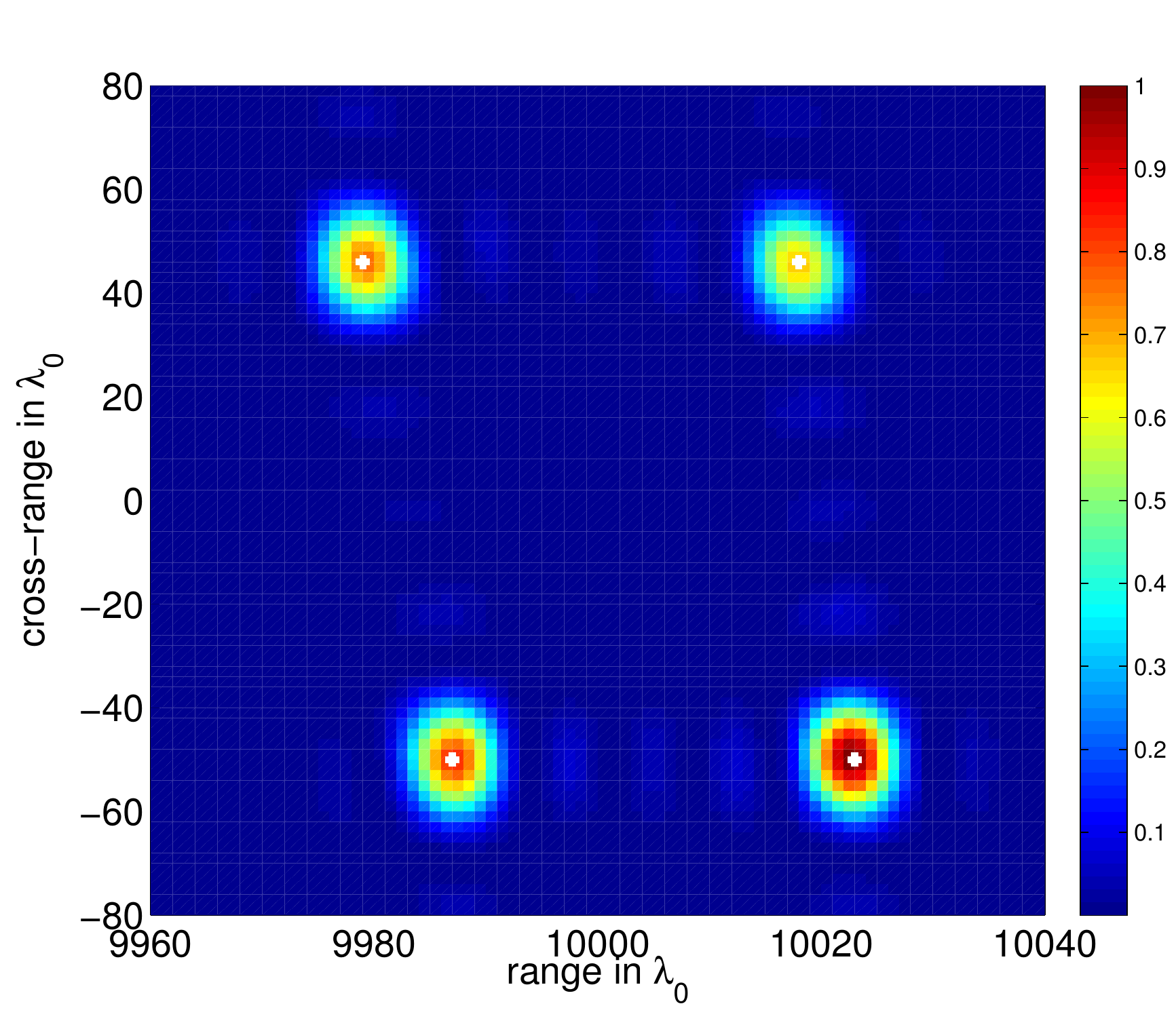}
\end{tabular}
\end{center}
\caption{{\bf Single receiver multifrequency interferometric data}. Homogeneous medium. Imaging with ${\cal I}^{Interf}$ as defined in \eqref{eq:interfR} using 16 frequencies equally spaced in $B=[580,620]$THz. For the left image the scatterers are on the grid while for the right image the scatterers are off the grid.}
\label{fig:h3}
\end{figure}

First, we assume that the medium between the array and the IW is homogeneous. In Figure~\ref{fig:h3}, 
we consider the same setup as in Figure \ref{fig:h2} except that (i) we do not record phases, and (ii) we only use one receiver. 
We see that the images shown in Figure~\ref{fig:h3}, obtained  with ${\cal I}^{Interf}$ as defined in \eqref{eq:interfR}, 
are of the same quality as those in Figure~\ref{fig:h2}, obtained with $\mathcal{I}^{KM}$ as defined in \eqref{eq:KM}, when phases are recorded 
and all source/receiver elements of the array are used for imaging. This means that imaging can be done
just as well without phases if one controls the illuminations! Moreover, the results do not deteriorate when the scatterers are off-grid,
as can be seen in the right plot of Figure~\ref{fig:h3}. This robustness is important and, as we will see, it will persist even when considering imaging 
in inhomogeneous, random media. 

Indeed, there is a similarity between these two types of data uncertainties. 
They both induce errors in the measured (or recovered) phases.
The off-grid case, however, is a systematic error, which is the same for all array elements, while the error 
induced by the random phase model depends on the path that connects the scatterer to each array element (see \eqref{eq:random_green_func} below). 
We note that depending on the correlation length of the random medium, the errors in the  measured (or recovered) phases are more or less 
correlated across the array elements. 

\begin{figure}[htbp]
\begin{center}
\begin{tabular}{ccc}
\includegraphics[scale=0.25]{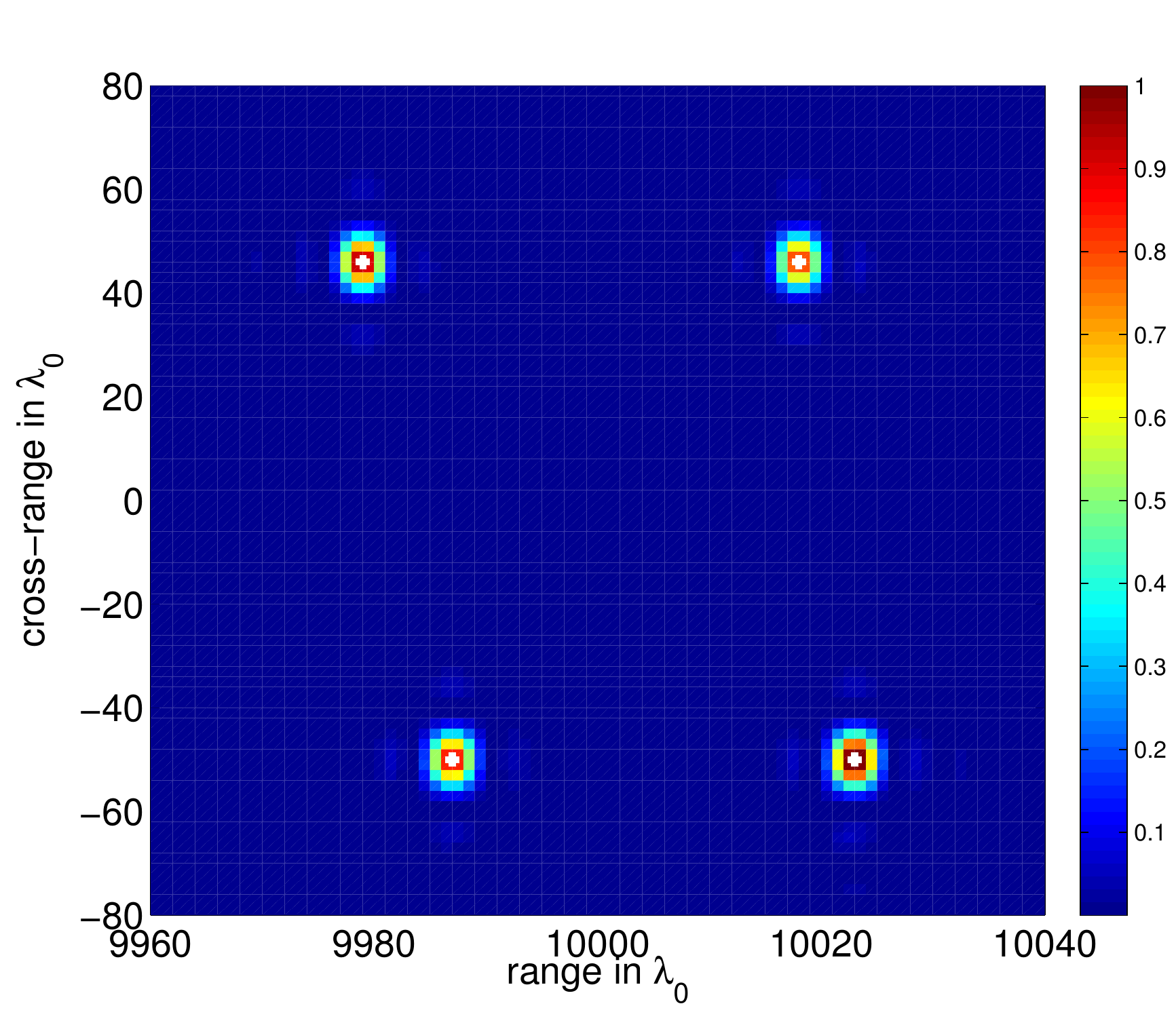}
\end{tabular}
\end{center}
\caption{Same as Figure \ref{fig:h3} but {\bf doubling the array size and the bandwidth}.  Homogeneous medium. The scatterers are off the grid. No additive noise in the data. }
\label{fig:resolution}
\end{figure}

Before considering imaging in random media, we note that the resolution obtained with the imaging functional~\eqref{eq:interfR} 
is, as for conventional imaging, of the order of $ \lambda_0 L/a$ in cross-range, and $C_0/B =  \lambda_0 f_0/B$ in range \cite{BPT-05}. 
This is confirmed in Figure~\ref{fig:resolution}, where we show an image obtained 
after doubling the array size and the bandwidth compared to Figure~\ref{fig:h3}. In Figure~\ref{fig:resolution}, we observe an 
improvement in resolution by a factor of~2 in both the cross-range and range (depth) directions compared to Figure~\ref{fig:h3}.
 
\subsection{Robustness and resolution in inhomogeneous media}
\begin{figure}[htbp]
\begin{center}
\begin{tabular}{ccc}
\includegraphics[scale=0.70]{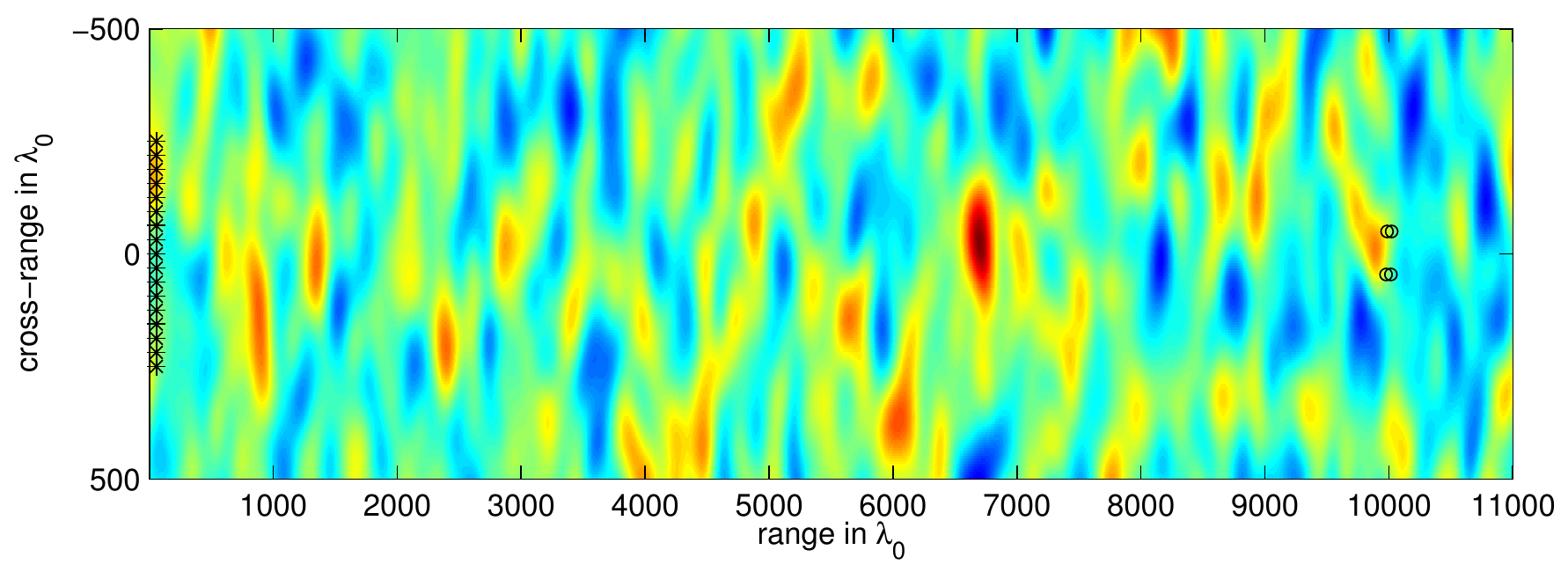} 
\end{tabular}
\end{center}
\caption{One realization of the random medium used in the simulations.} 
\label{fig:rm0}
\end{figure}

To study the performance of the imaging functionals in heterogeneous media, we consider the setup shown in Figure~\ref{fig:rm0}. 
The four scatterers (very close to each other to be seen) on the right are shown with black disks, 
and the array elements on the left are indicated with black stars. The array response matrix is computed using \eqref{responsematrix} 
and the random travel time model for the Green's function \eqref{eq:random_green_func} explained below. 
The medium fluctuations are modeled as in \eqref{eq:random_wave_speed}. The correlation length of the fluctuations is $l=100 \lambda_0$,
and the amplitude of the fluctuations is $\sigma=4 \,10^{-4}$. Here, we increase the bandwidth to the maximum available 
$B=[540,660]$THz, 
and we discretize it using $46$ equally spaced frequencies. The size of the IW is $160 \lambda_0 \times 80 \lambda_0$,
and the pixel size is $4  \lambda_0 \times 2  \lambda_0$ in cross-range and range (depth), respectively. 

\begin{figure}[htbp]
\begin{center}
\begin{tabular}{cccc}
\includegraphics[scale=0.15]{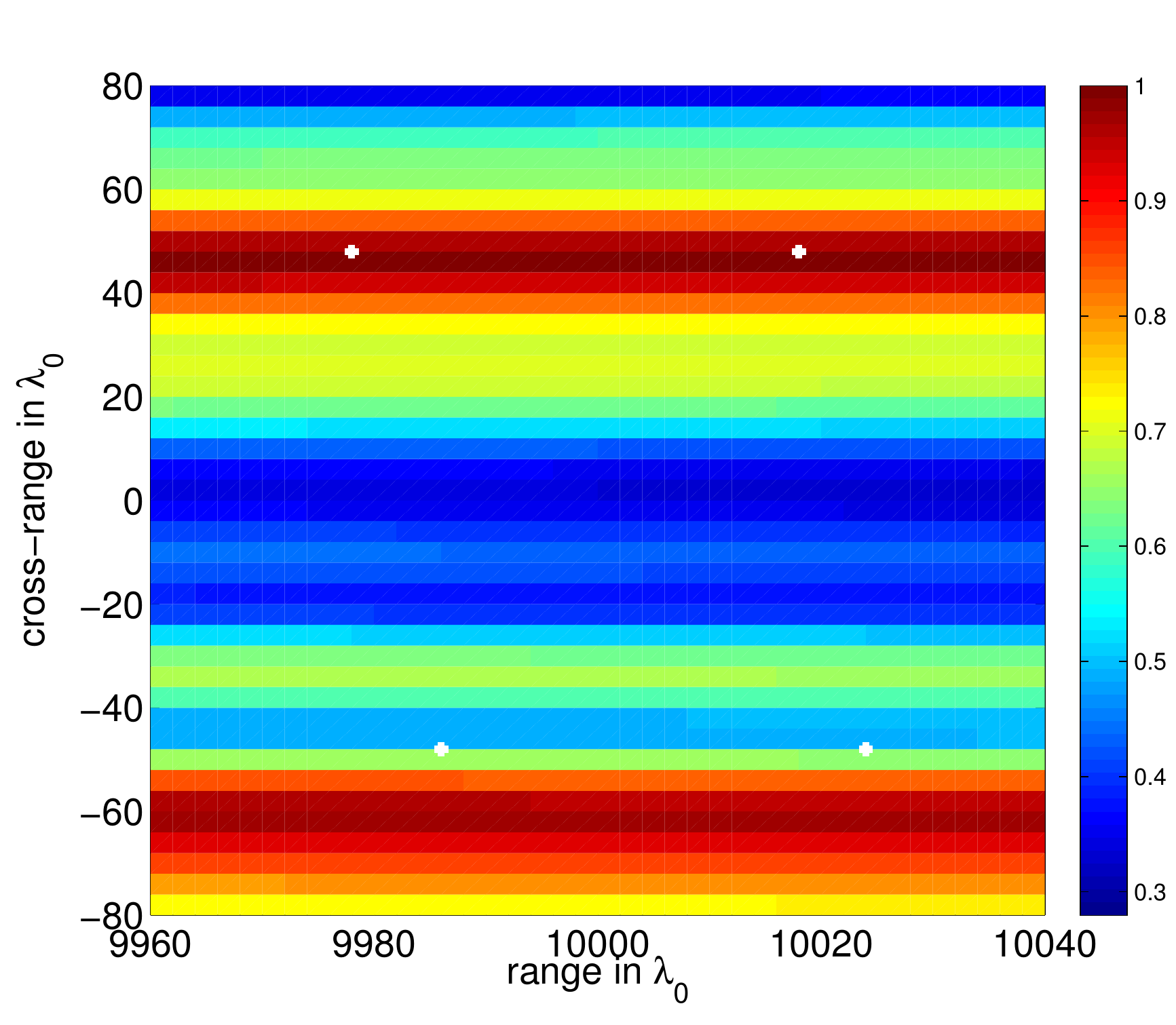}&
\includegraphics[scale=0.15]{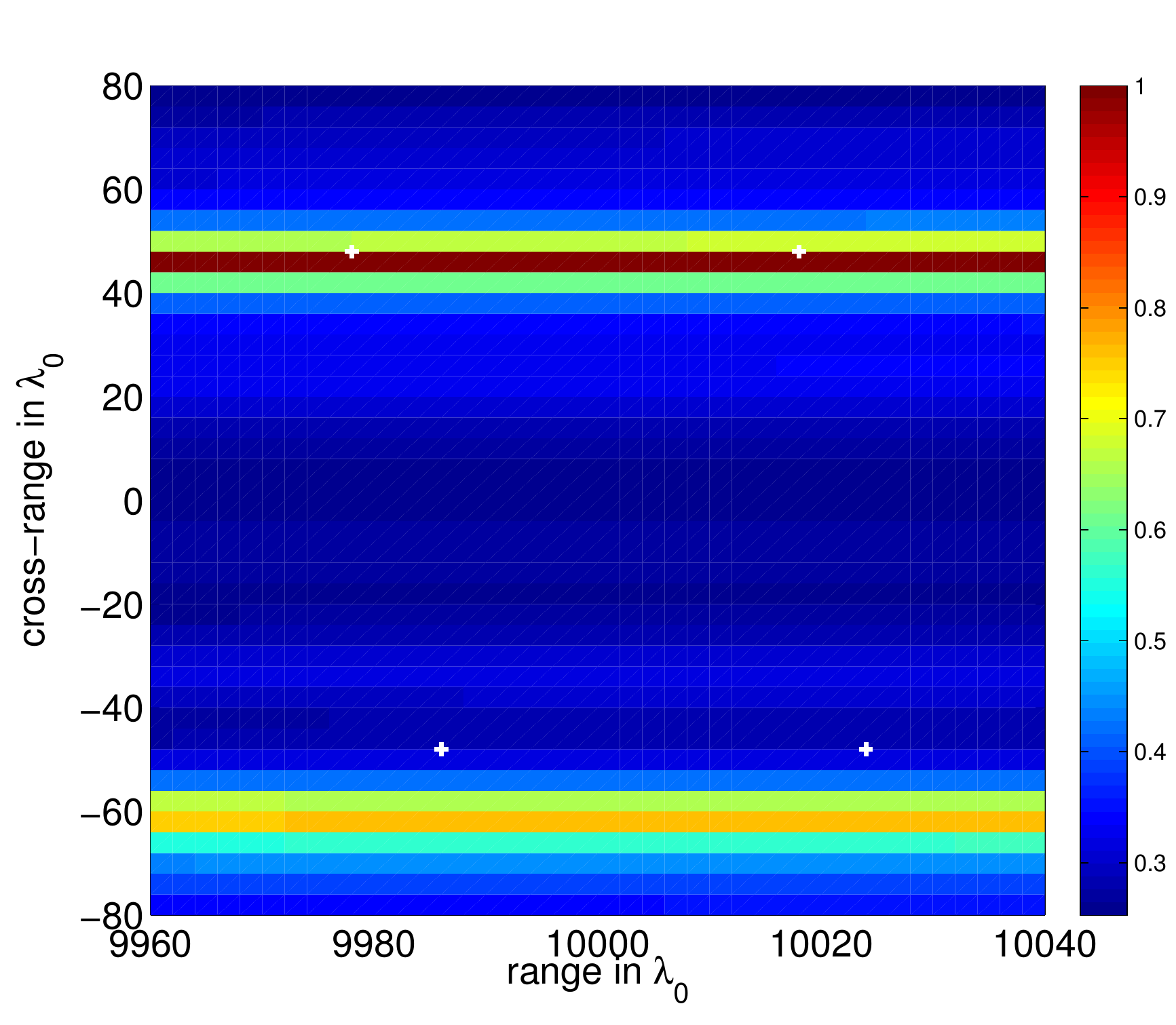}&
\includegraphics[scale=0.15]{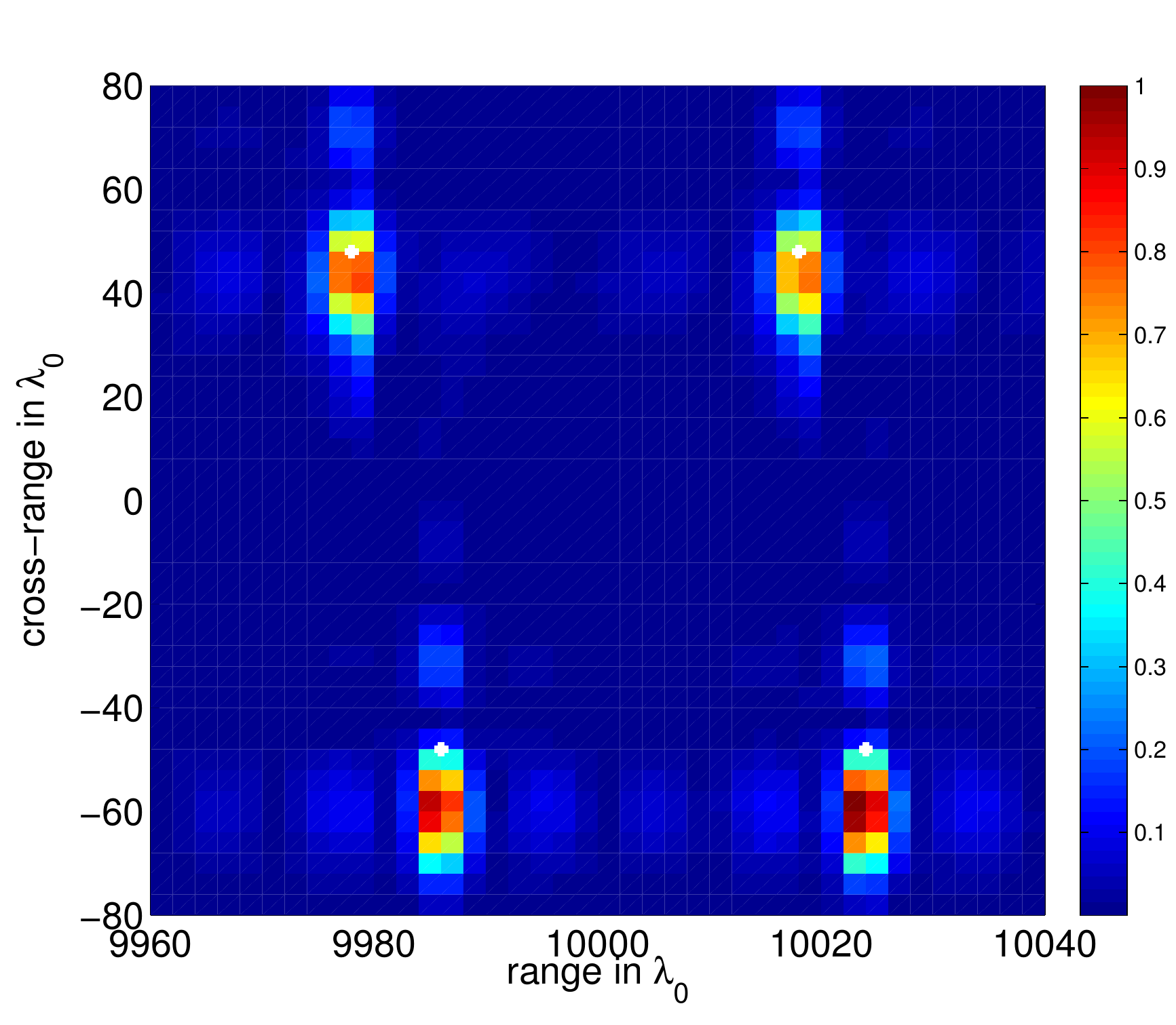}&
\includegraphics[scale=0.15]{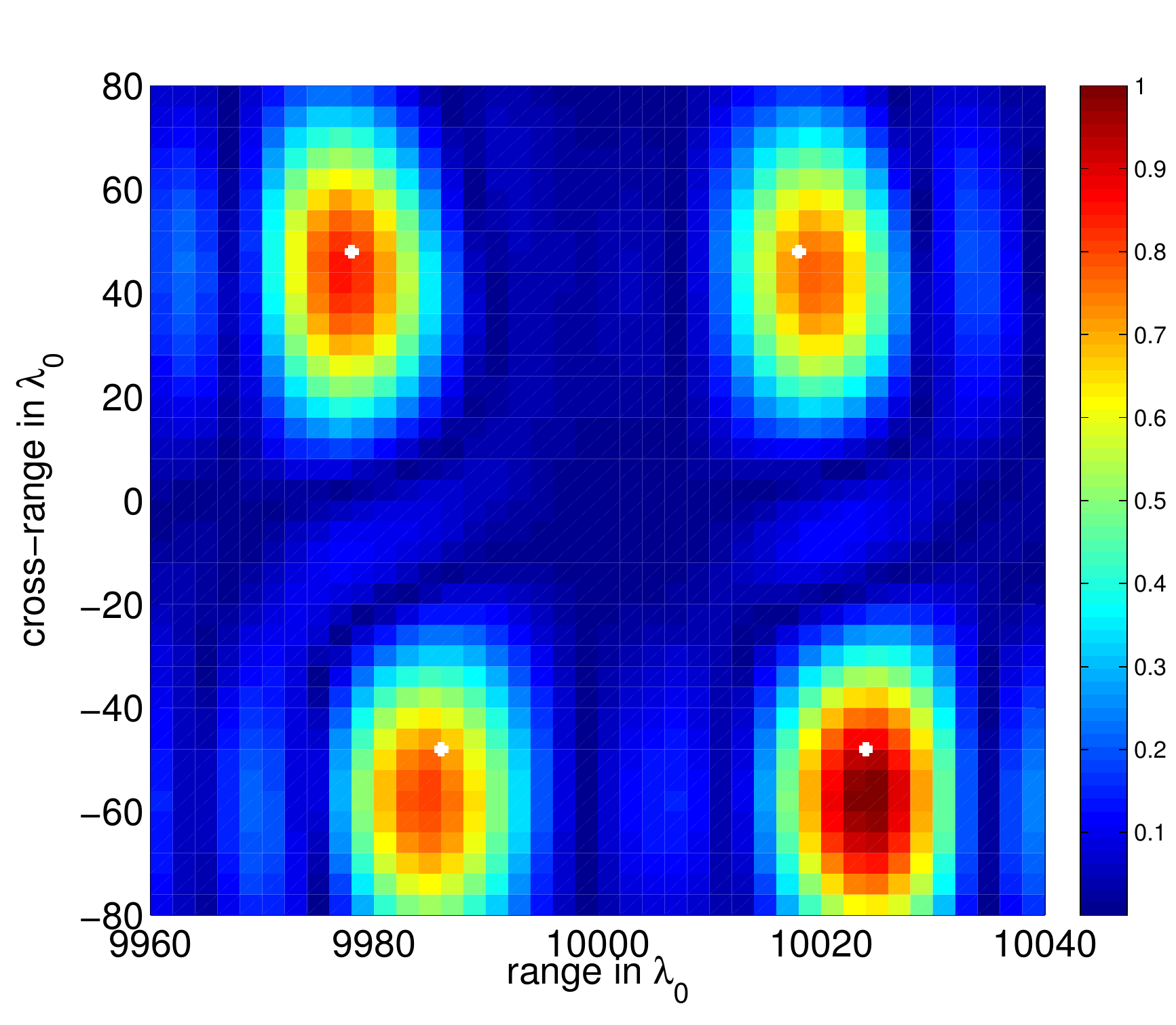}\\
\includegraphics[scale=0.15]{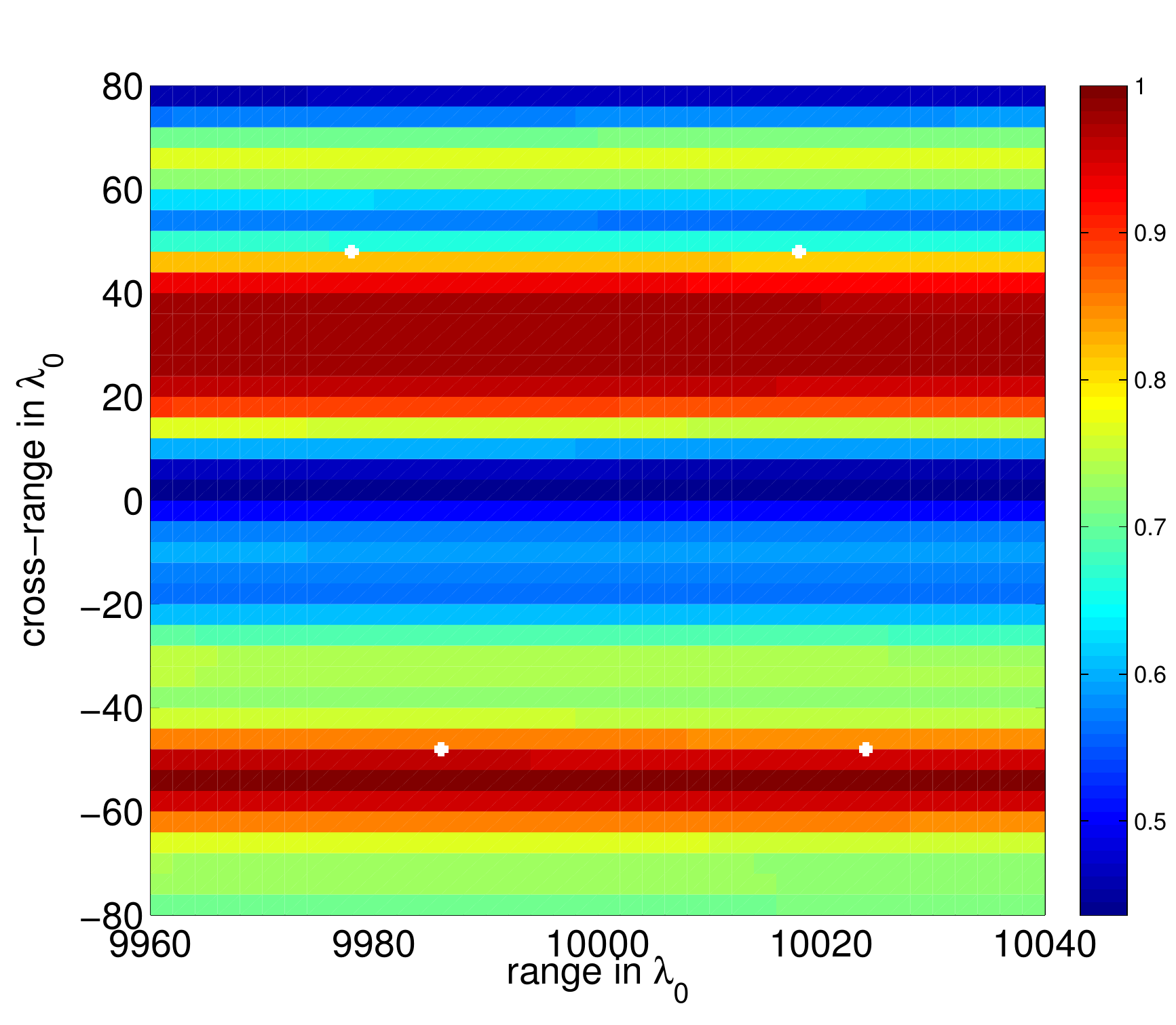}&
\includegraphics[scale=0.15]{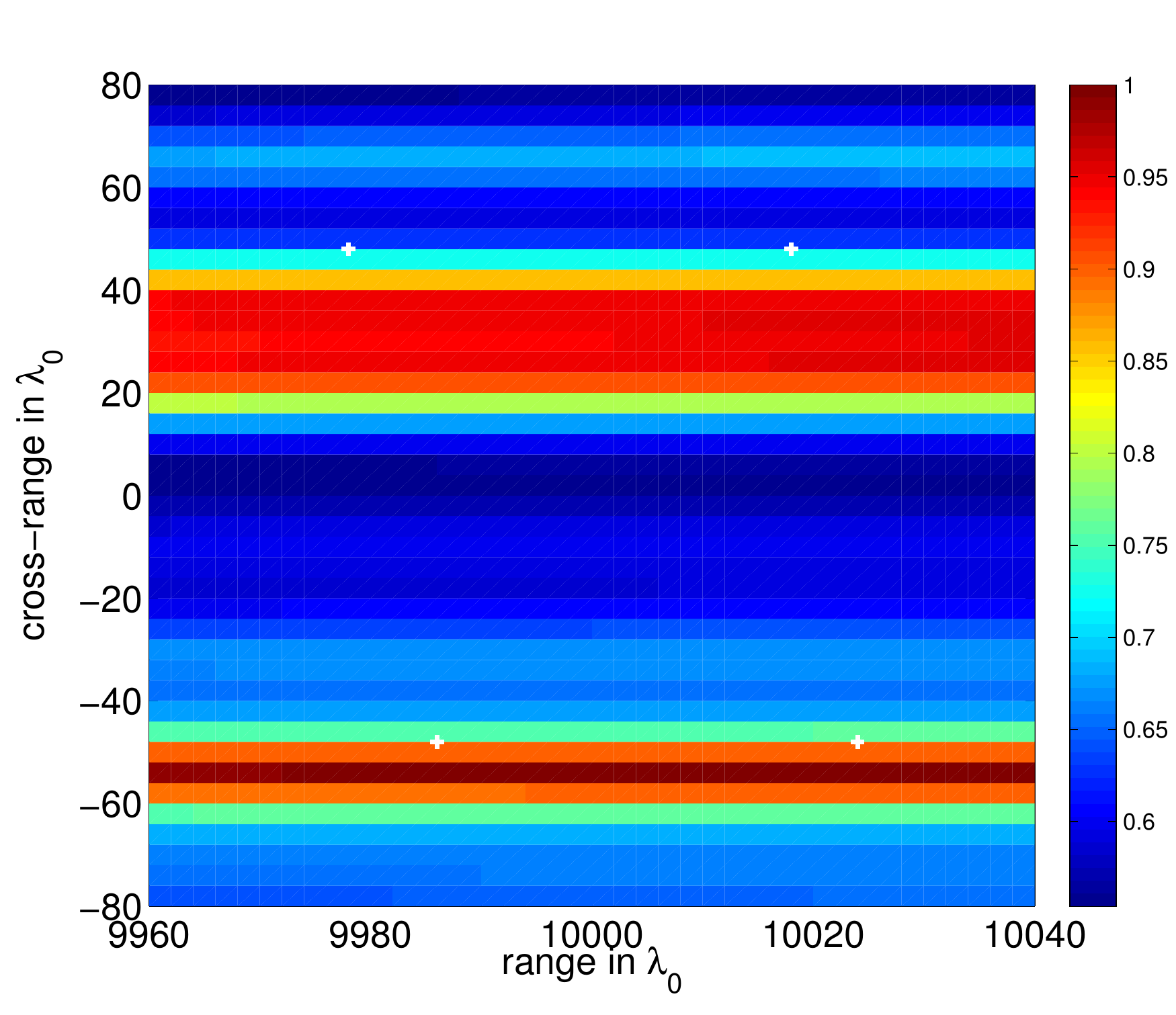}&
\includegraphics[scale=0.15]{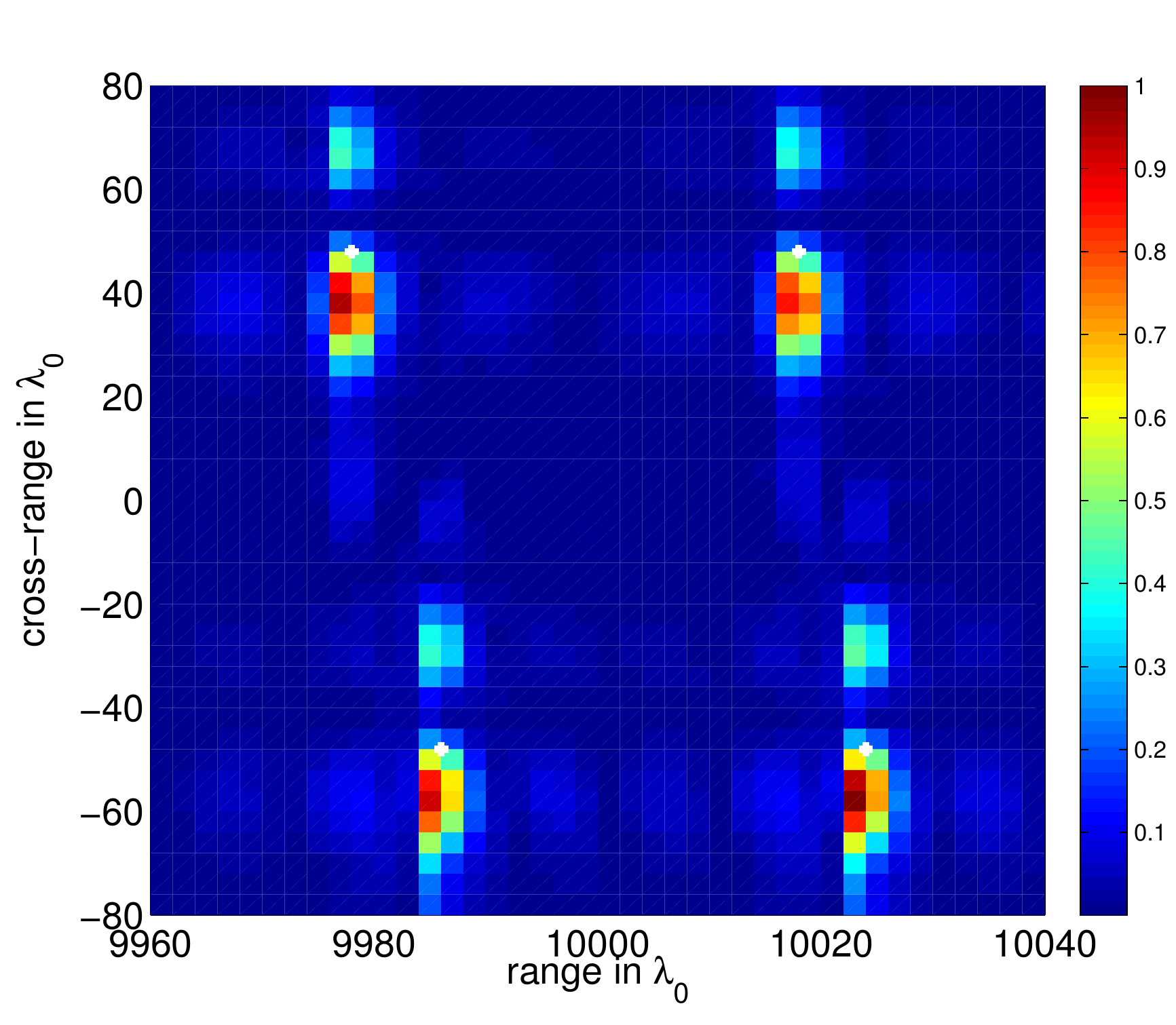}&
\includegraphics[scale=0.15]{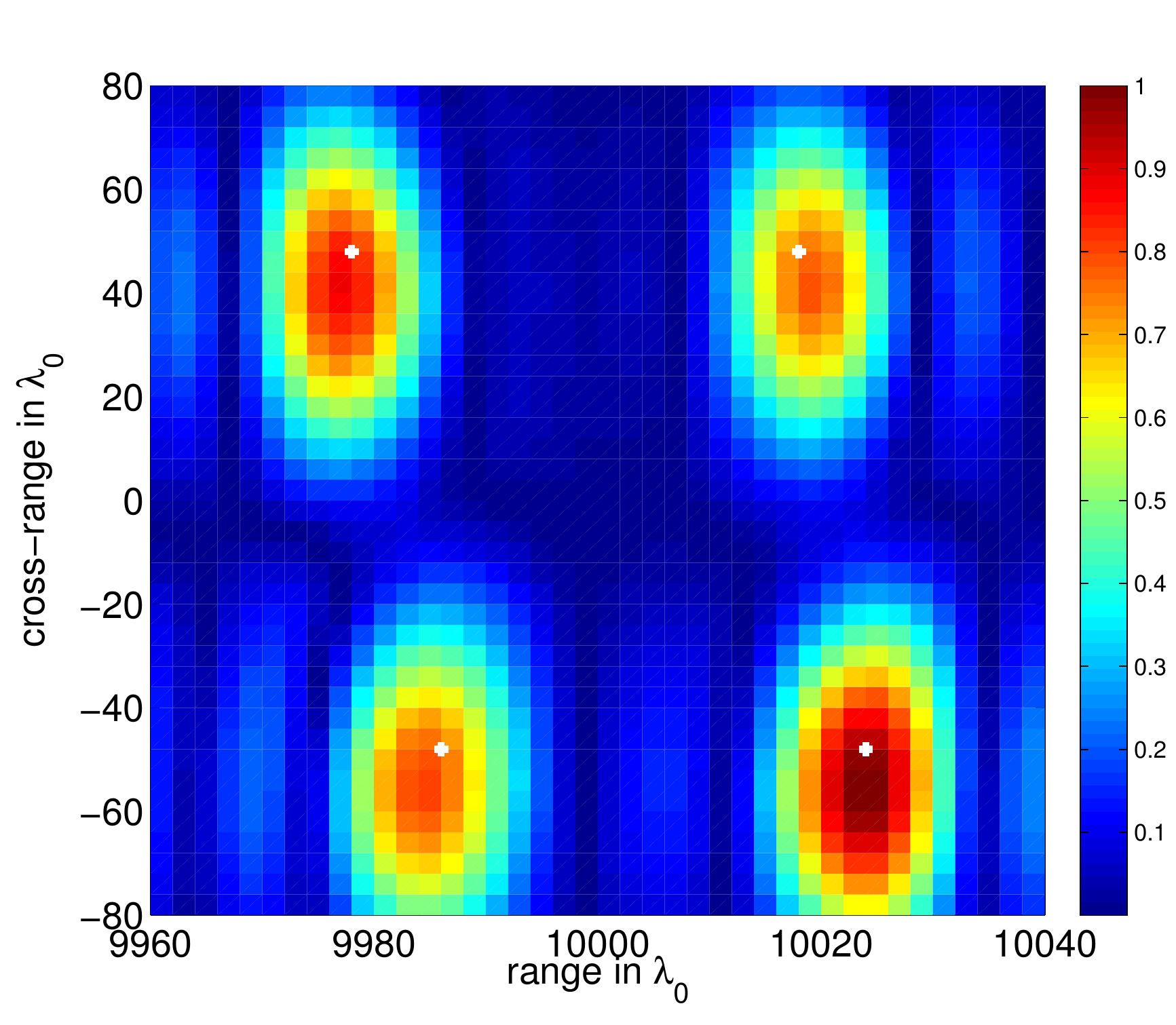}\\
\includegraphics[scale=0.15]{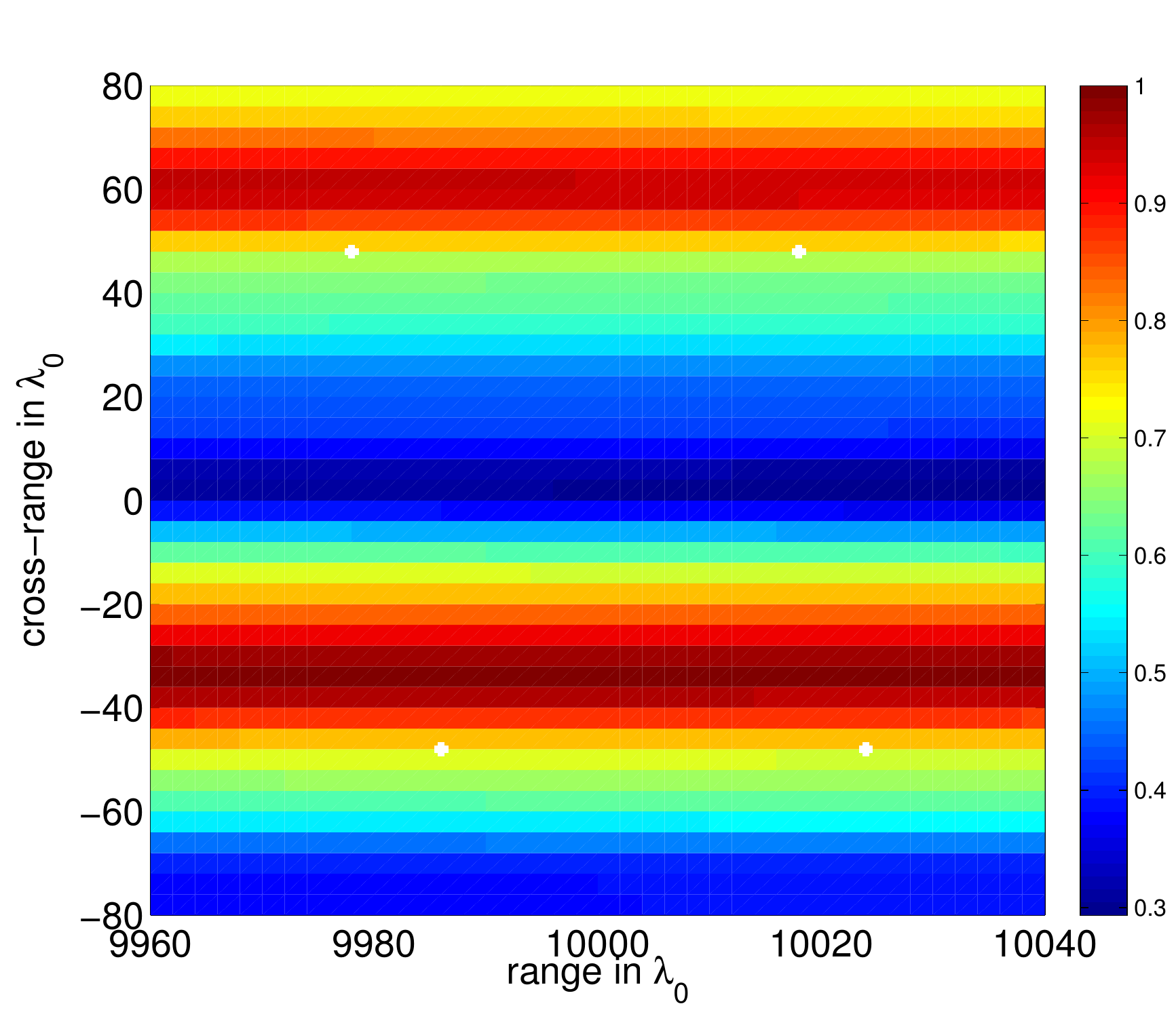}&
\includegraphics[scale=0.15]{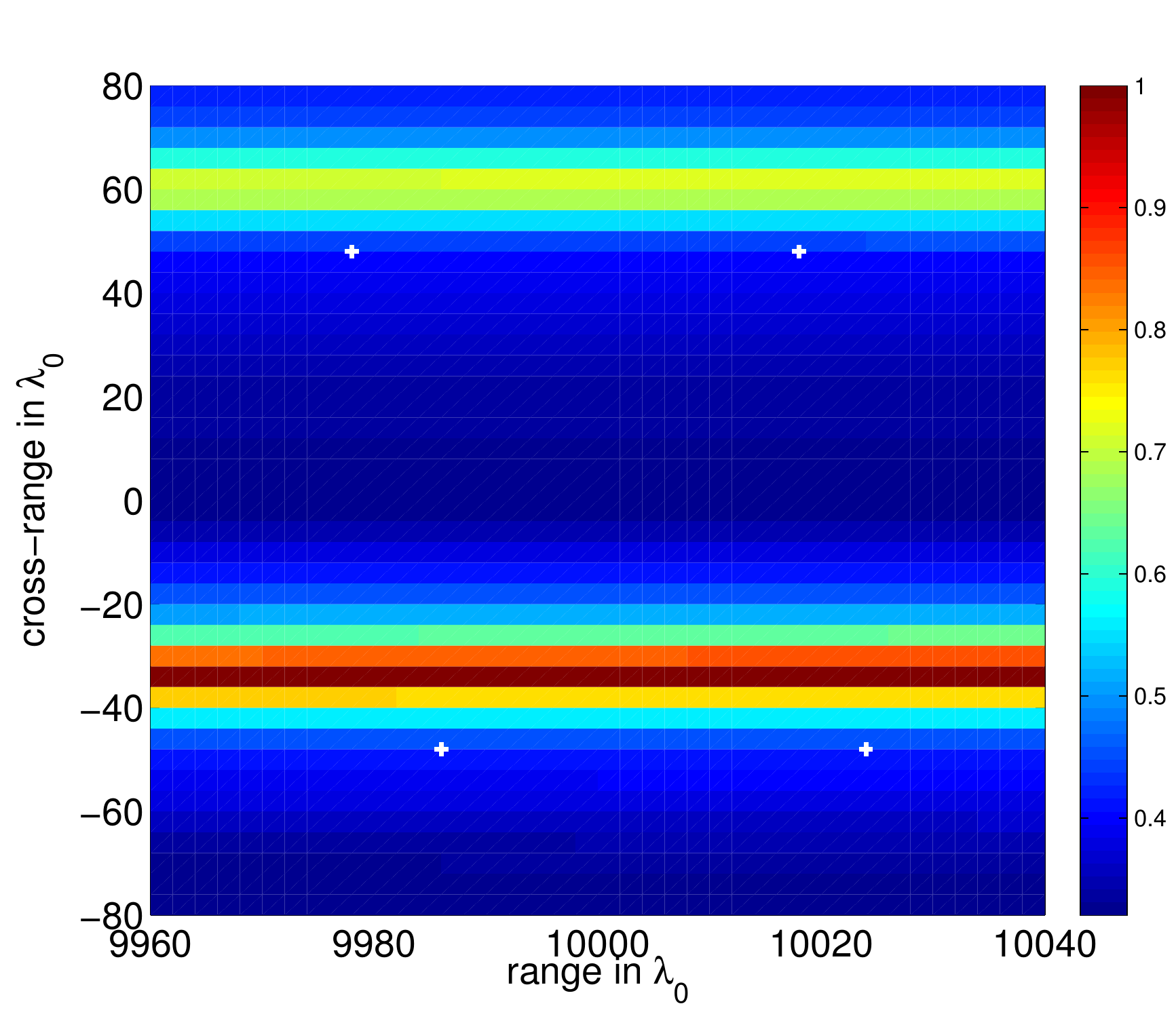}&
\includegraphics[scale=0.15]{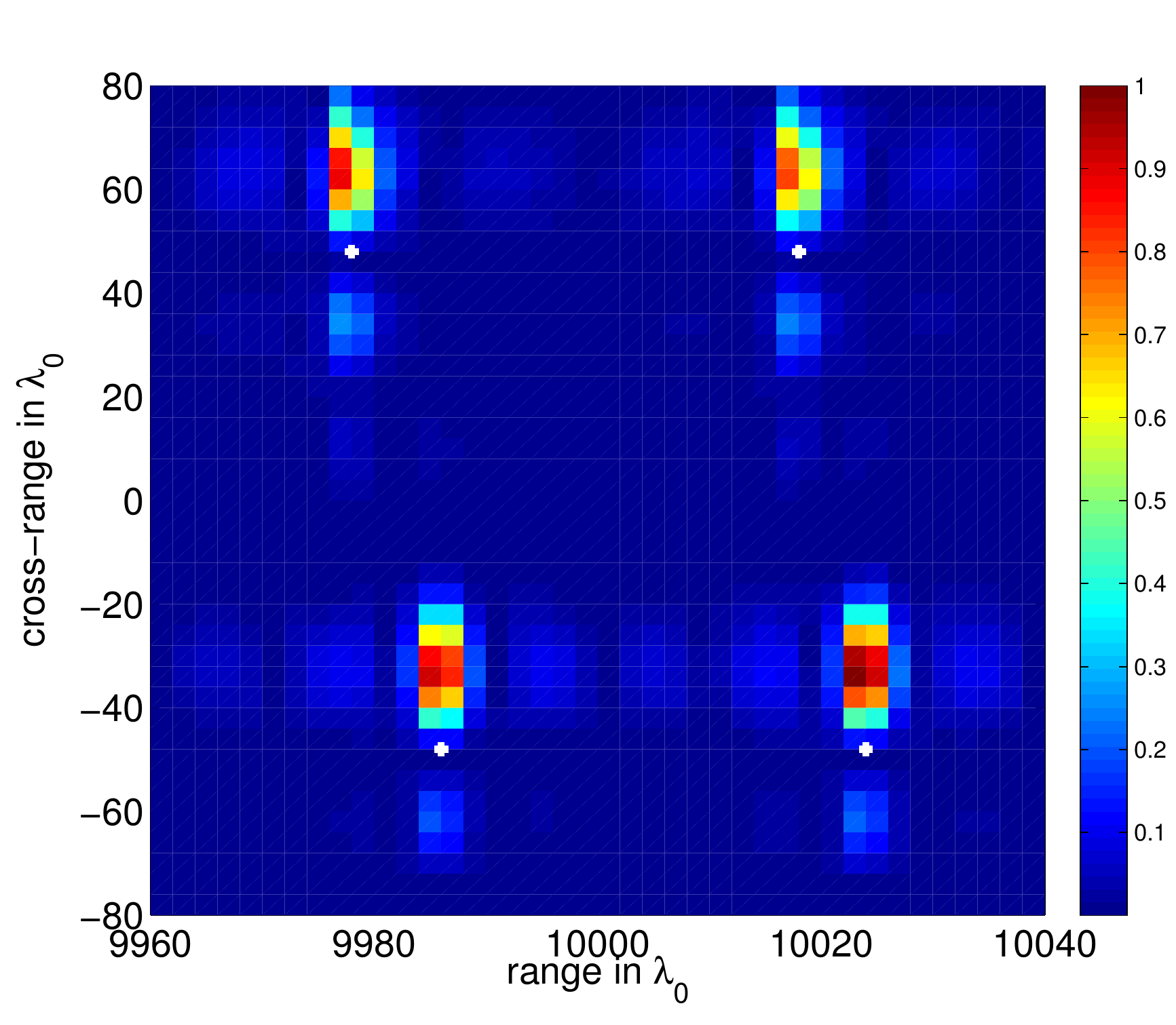}&
\includegraphics[scale=0.15]{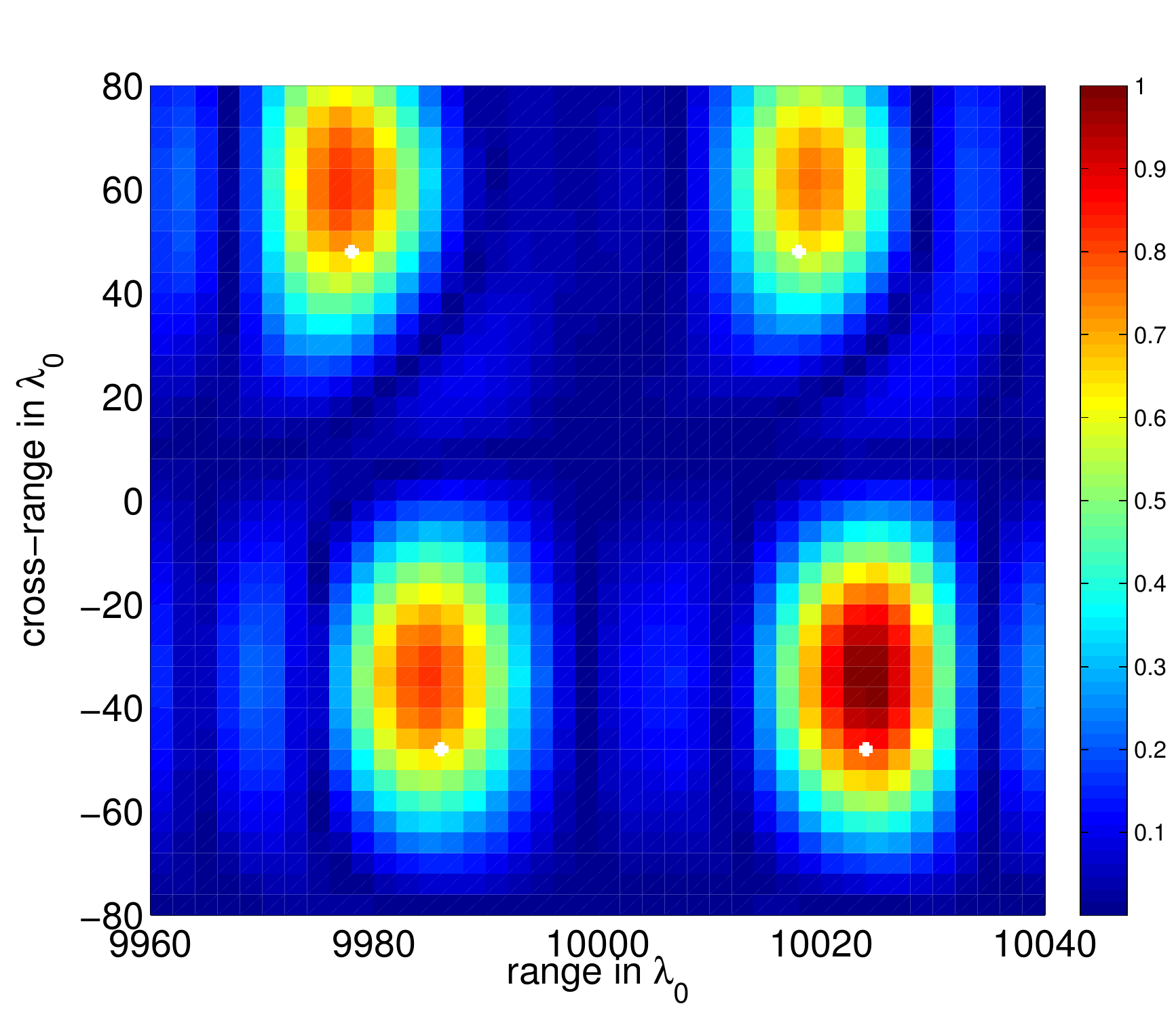} 
\end{tabular}
\end{center}
\caption{{\bf Multifrequency interferometric data}. We use 46 frequencies equally spaced in the bandwidth $B=[540,660]$THz. {\bf Random medium}. The scatterers are off the grid. From left to right: 
 $\mathcal{I}^{SIGNAL}$ as defined on \eqref{SIGNAL}, $\mathcal{I}^{MUSIC}$ as defined in \eqref{MUSIC}, $\mathcal{I}^{Interf}$ as defined in \eqref{eq:interfR},  
and  $\mathcal{I}^{SRINT}$ as defined in \eqref{eq:AMA}  using a mask with $\Omega_d=0.12B$ and $X_d=0.25 a$.  From top to bottom three realizations of the random medium. 
Note that for  $\mathcal{I}^{SIGNAL}$ and $\mathcal{I}^{MUSIC}$ we add incoherently over the multiple frequencies while for $\mathcal{I}^{Interf}$ and $\mathcal{I}^{SRINT}$ the image is constructed by adding 
coherently over the multifrequency data. Note also that to construct the images in the first two columns we need to recover the full matrix $\Mm$ while for the last two columns only the single receiver element matrix $\Mm_r$ is used.}
\label{fig:stability}
\end{figure}

We implement the interferometric imaging approach using masks as described in Section \ref{sec:masks} (see Eq. \eqref{eq:AMA}) with 
$\Omega_d=0.12B$ and $X_d=0.25 a$. By reducing the 
distance between the sensors and the frequencies used for imaging,
we gain stability but we lose some resolution. 
The parameters $\Omega_d$ and $X_d$ can be obtained by an optimization procedure \cite{BPT-ADA}. Here, however,
we experimented with different values of $\Omega_d$ and $X_d$, and we picked the ones that provide a good compromise between stability 
gain and resolution loss. 

\begin{figure}[htbp]
\begin{center}
\begin{tabular}{ccc}
\includegraphics[scale=0.32]{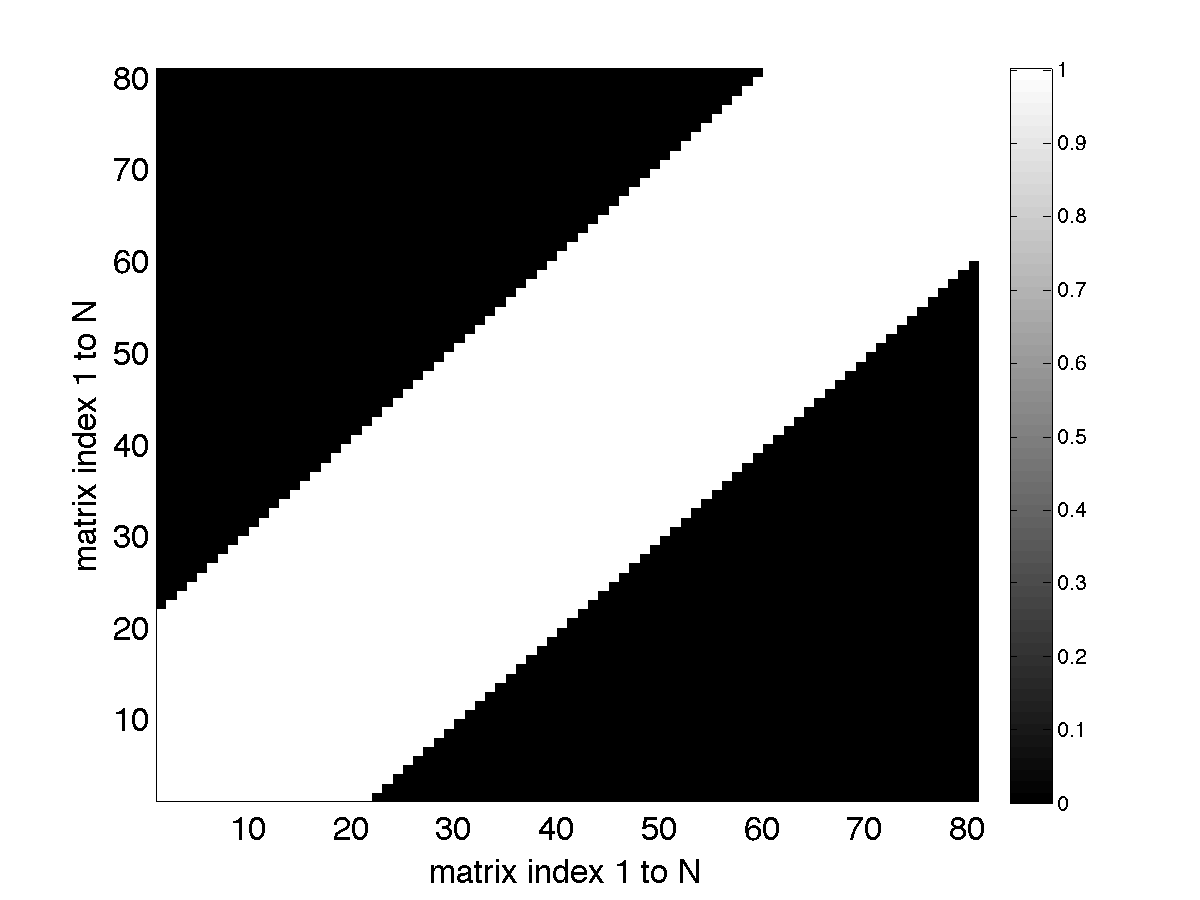}&
\includegraphics[scale=0.32]{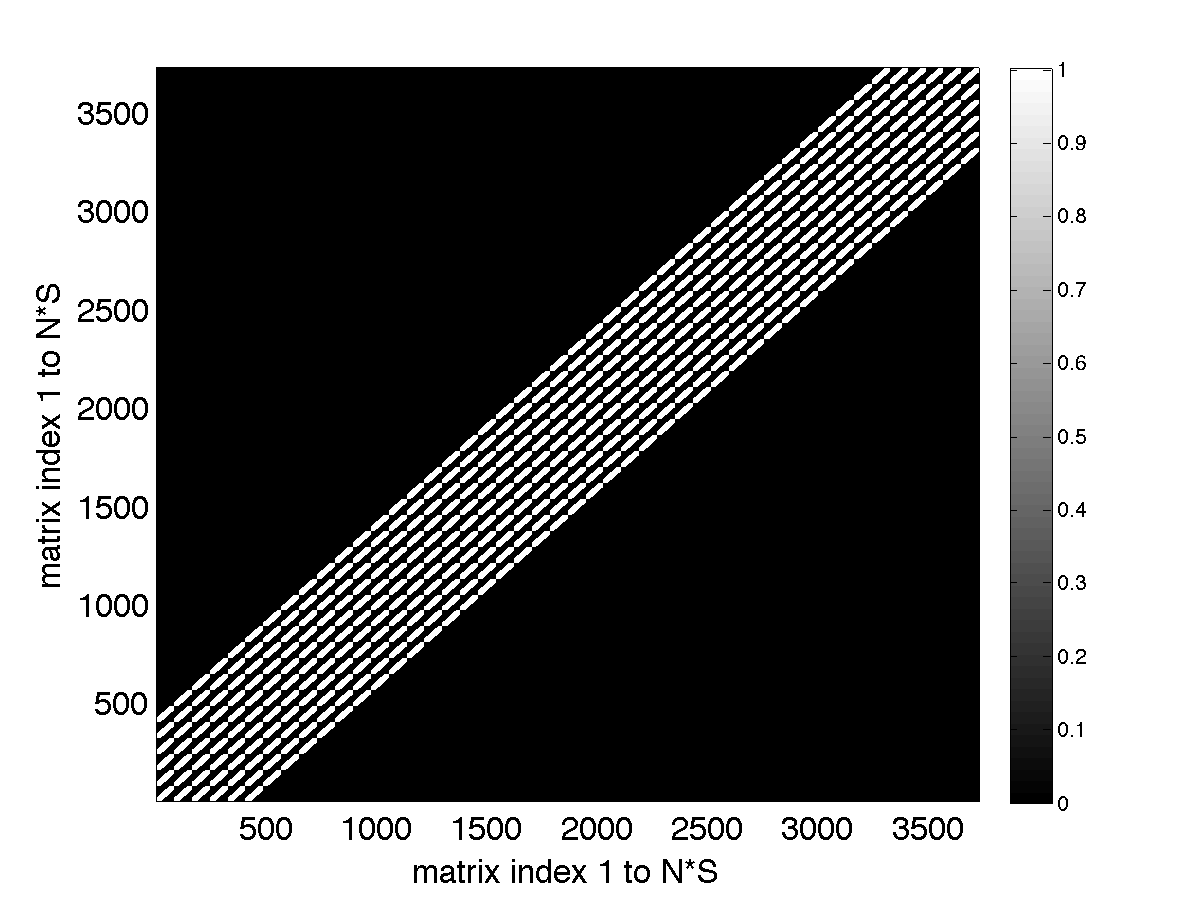}\\
\end{tabular}
\end{center}
\caption{Example of masks $\vect{\Zc}$. On the left is the mask corresponding to single frequency data and limits the interferometric data only in space. On the right is the mask used in the numerical simulations that limits interferometric data to nearby frequencies and nearby sources. The values of the matrix are zeros and ones, with one plotted in white and zero in black.}
\label{fig:masks}
\end{figure}

Figure~\ref{fig:stability} shows the images in three different realizations of a random medium. From left to right we show 
the images obtained with $\mathcal{I}^{SIGNAL}$ as defined on \eqref{SIGNAL}, $\mathcal{I}^{MUSIC}$ as defined in \eqref{MUSIC}, 
$\mathcal{I}^{Interf}$ as defined in \eqref{eq:interfR}, and $\mathcal{I}^{SRINT}$ as defined in \eqref{eq:AMA}. 
The mask used in $\mathcal{I}^{SRINT}$ is displayed in the right plot of Figure \ref{fig:masks}. Both imaging functionals $\mathcal{I}^{Interf}$ and 
$\mathcal{I}^{SRINT}$ use the multifrequency interferometric data matrix obtained from the intensities gathered at the receiver located at the center 
of the array. Similar results, not shown here, are obtained for other receiver locations. The receiver location does not really affect the imaging results. 
This is so because the array is small and is located at a large distance from the image window (IW). 

In all the cases shown in Figure~\ref{fig:stability}, the scatterers do not lie on the grid and hence, as expected, neither $\mathcal{I}^{SIGNAL}$ 
nor  $\mathcal{I}^{MUSIC}$ (first and second columns) are able to locate the scatterers because the range resolution is lost, 
as was the case in a homogeneous medium (see Figure~\ref{fig:h3}).
In addition, we now observe that when the medium is random, the estimated cross-range varies from one realization to another. 
This phenomenon is also noticeable in the  ${\cal I}^{Interf}$ images shown in the third column. 
Although the resolution of the ${\cal I}^{Interf}$ images is far better than 
the one given by  
$\mathcal{I}^{SIGNAL}$ and  $\mathcal{I}^{MUSIC}$, the peaks obtained in the ${\cal I}^{Interf}$ images dance around the true locations of the 
scatterers, meaning that ${\cal I}^{Interf}$ is not stable and, therefore, it is not useful for imaging in regimes with significant wavefront distortions. 
However, $\mathcal{I}^{SRINT}$ (fourth column), that uses masked multifrequency interferometric data, 
gives statistically stable results in these media. We also observe a significant loss in resolution due to the use of the masks.
Because $\Omega_d=0.12B$ and $X_d=0.25 a$, we only use a small part of the available bandwidth and of the array aperture and, therefore, the 
resolution decreases to $ \lambda_0 L/X_d$  and  $C_0/\Omega_d$  in cross-range and range, respectively (cf. \cite{BPT-05}).

In Figure \ref{fig:masks} we display the mask $\vect{\Zc}$ used to produce the SRINT images in Figure~\ref{fig:stability}. 
The left plot is the monochromatic version of the mask which illustrates windowing only in the spatial direction. 
We observe that $\vect{\Zc}$ is a band matrix whose bandwidth limits the correlations used in $\mathcal{I}^{SRINT}$ so that $|\bx-\bx'|\ll X_d$. 
The right plot is the mask used in our simulations for $S=46$ frequencies and $N=81$ array elements. 
Recall that the index of $\Mm_r$, and therefore of the mask as well, is defined as $i= s+ (l-1)\cdot N$ for source location $s=1,\ldots,N$ and frequency 
$l=1,\ldots,S$.

\begin{figure}[htbp]
\begin{center}
\begin{tabular}{ccc}
\includegraphics[scale=0.22]{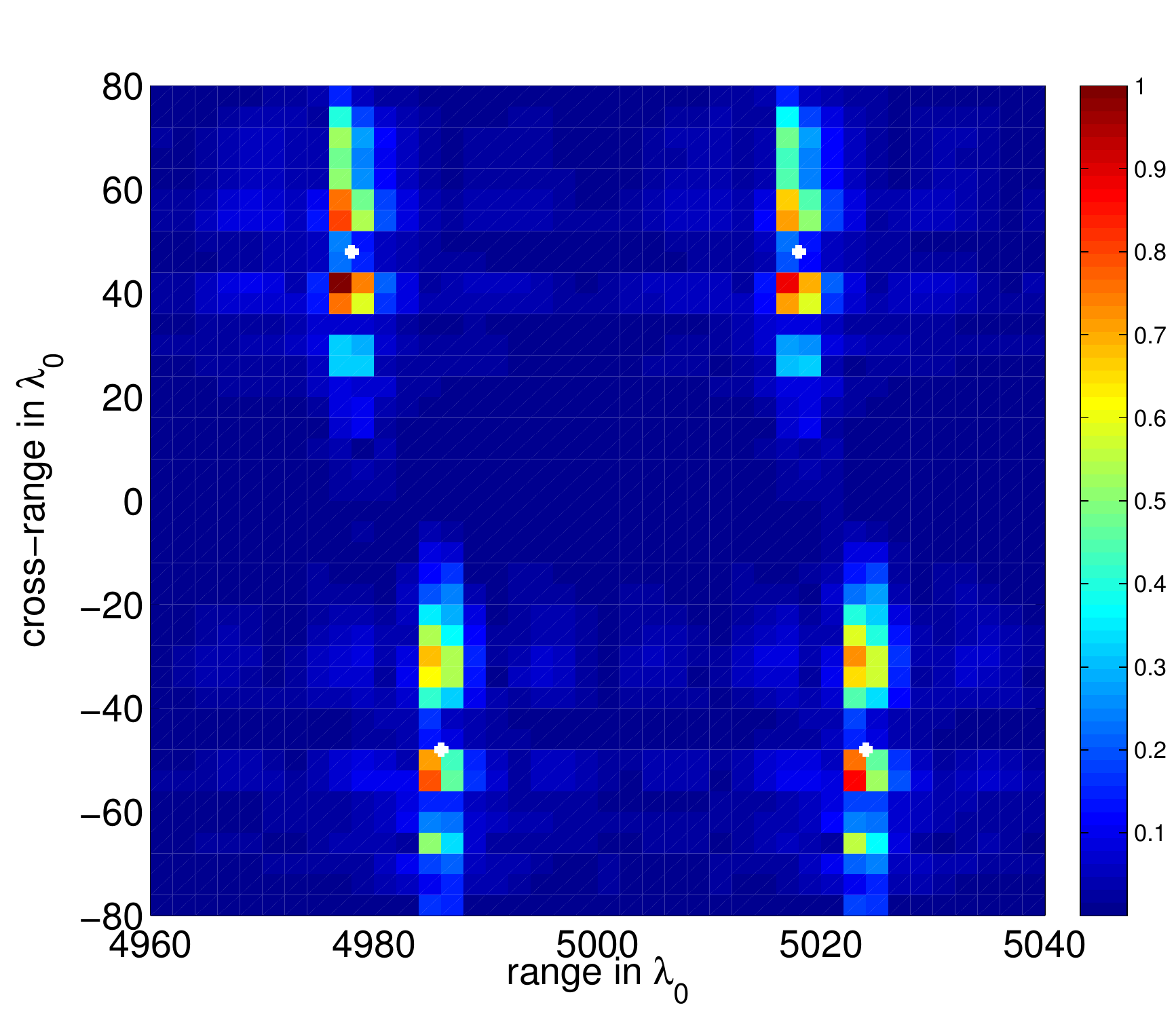}&
\includegraphics[scale=0.22]{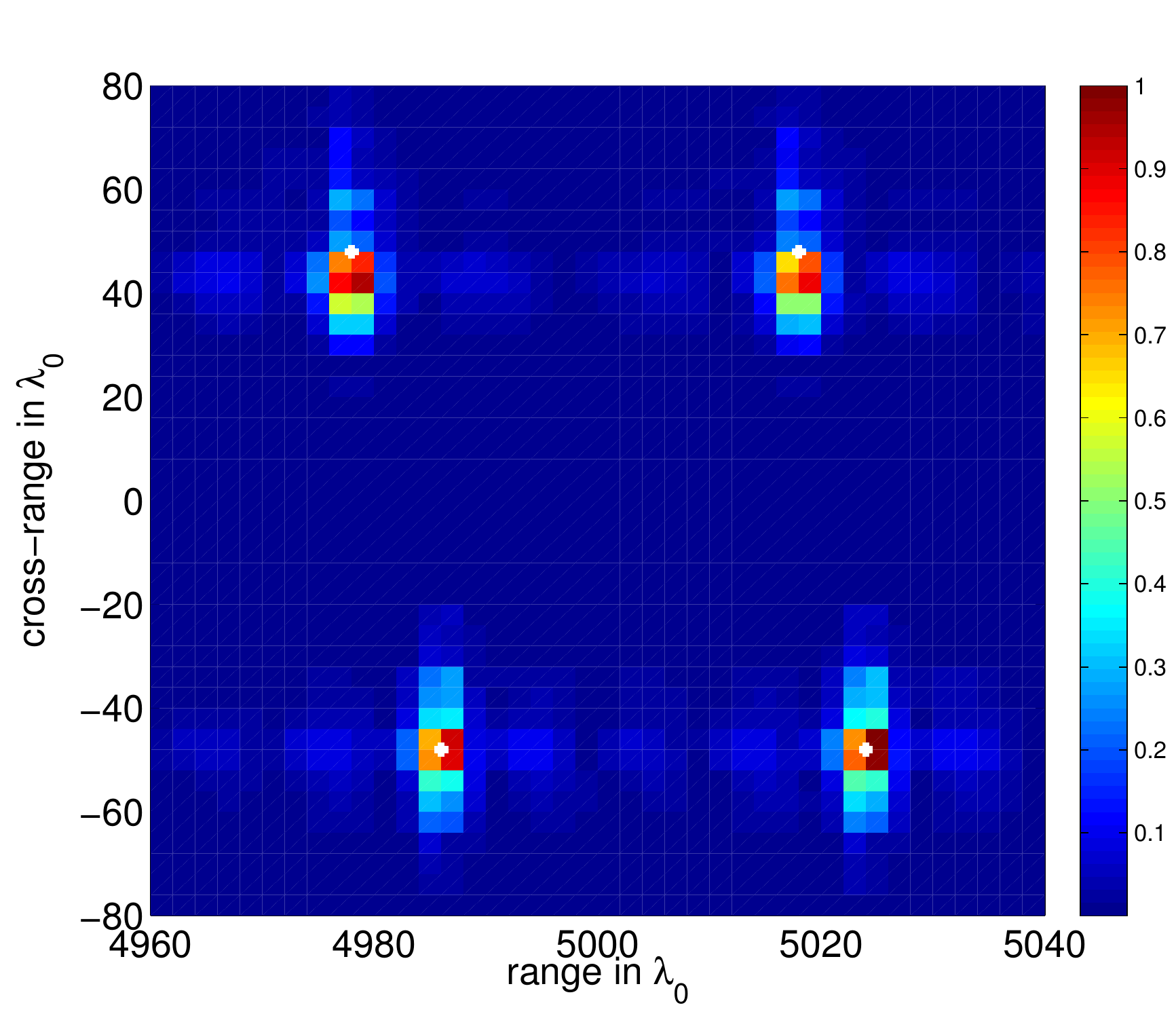}&
\includegraphics[scale=0.22]{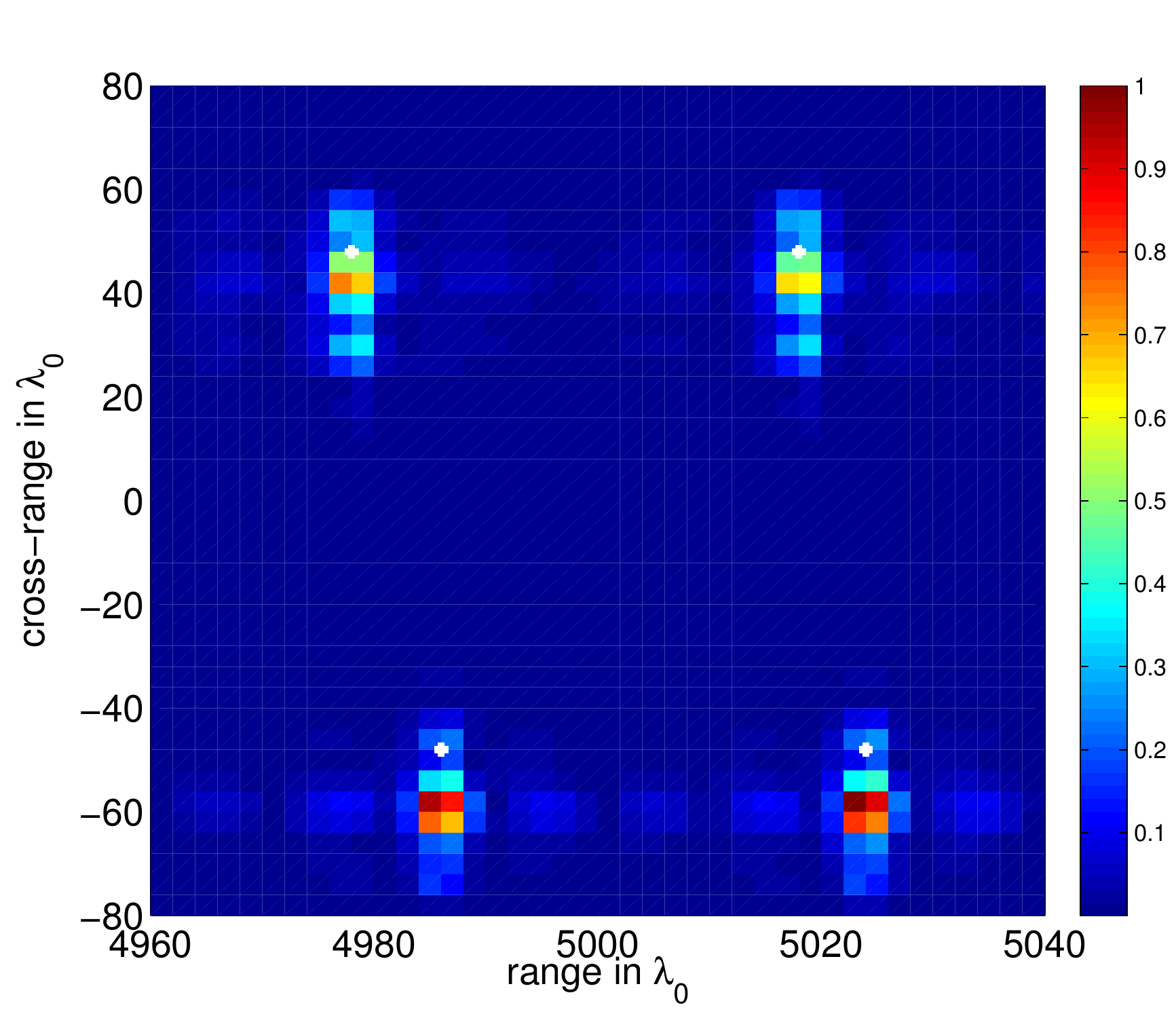} \\
\includegraphics[scale=0.22]{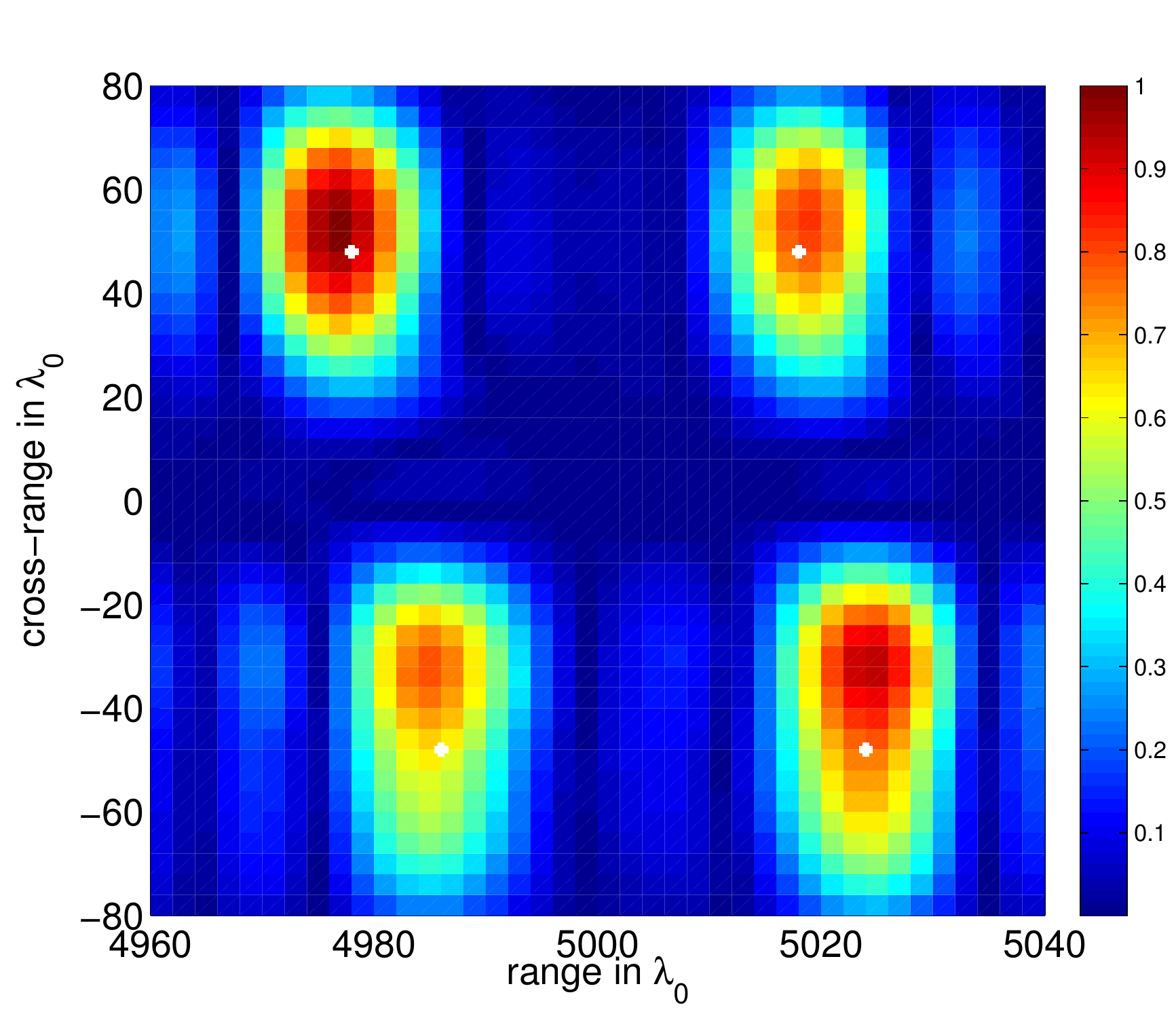}&
\includegraphics[scale=0.22]{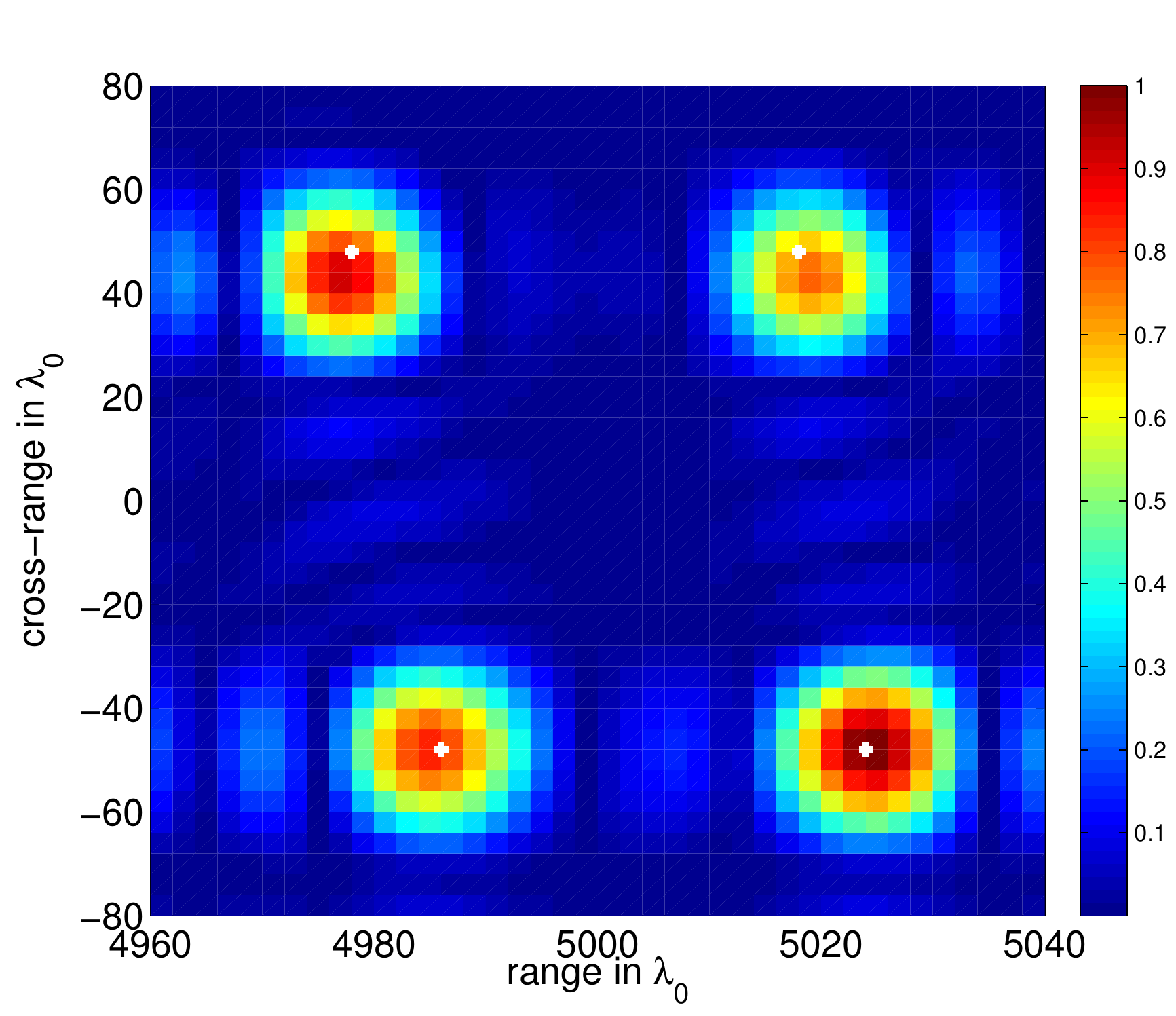}&
\includegraphics[scale=0.22]{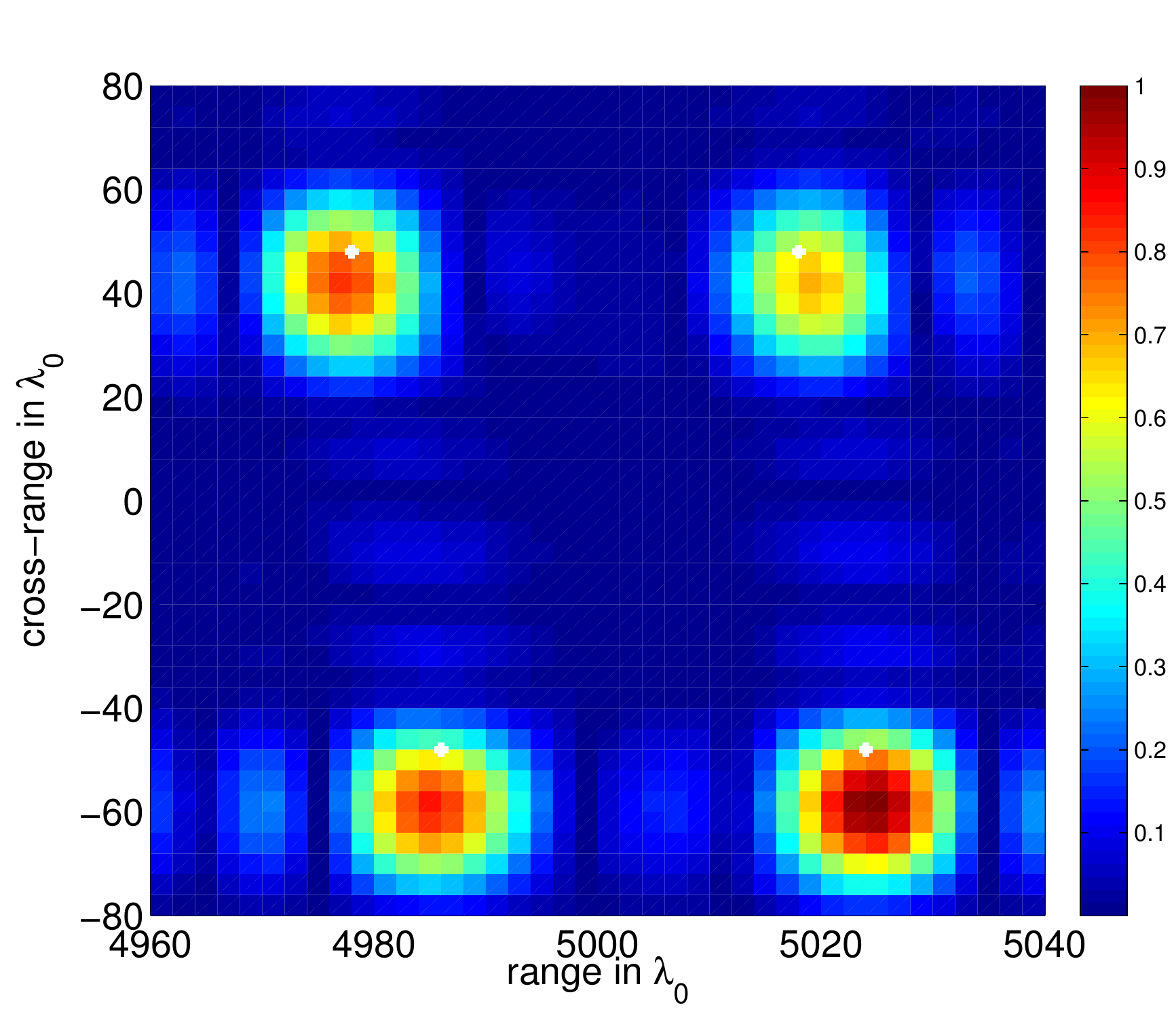}\\
\end{tabular}
\end{center}
\caption{{\bf Single receiver multifrequency interferometric data}. We use 46 frequencies equally spaced in the bandwidth $B=[540,660]$THz. {\bf Random medium}. The scatterers are off the grid. We change here the propagation distance to $L/2$ and at the same time we increase the strength of the fluctuations by a factor $\sqrt{2}$ (to keep $\sigma \sqrt{l L}/ \lambda_0$ constant). In the top row we show the ${\cal I}^{Interf}$ images and in the bottom the ${\cal I}^{SRINT}$ images that use the masked data. From left to right we illustrate results for three realizations of the random medium.  The same mask with $\Omega_d=0.12B$ and $X_d=0.25 a$ as before is used. We do not show the MUSIC images because they are just as bad as before.}
\label{fig:error1}
\end{figure}

\begin{figure}[htbp]
\begin{center}
\begin{tabular}{ccc}
\includegraphics[scale=0.22]{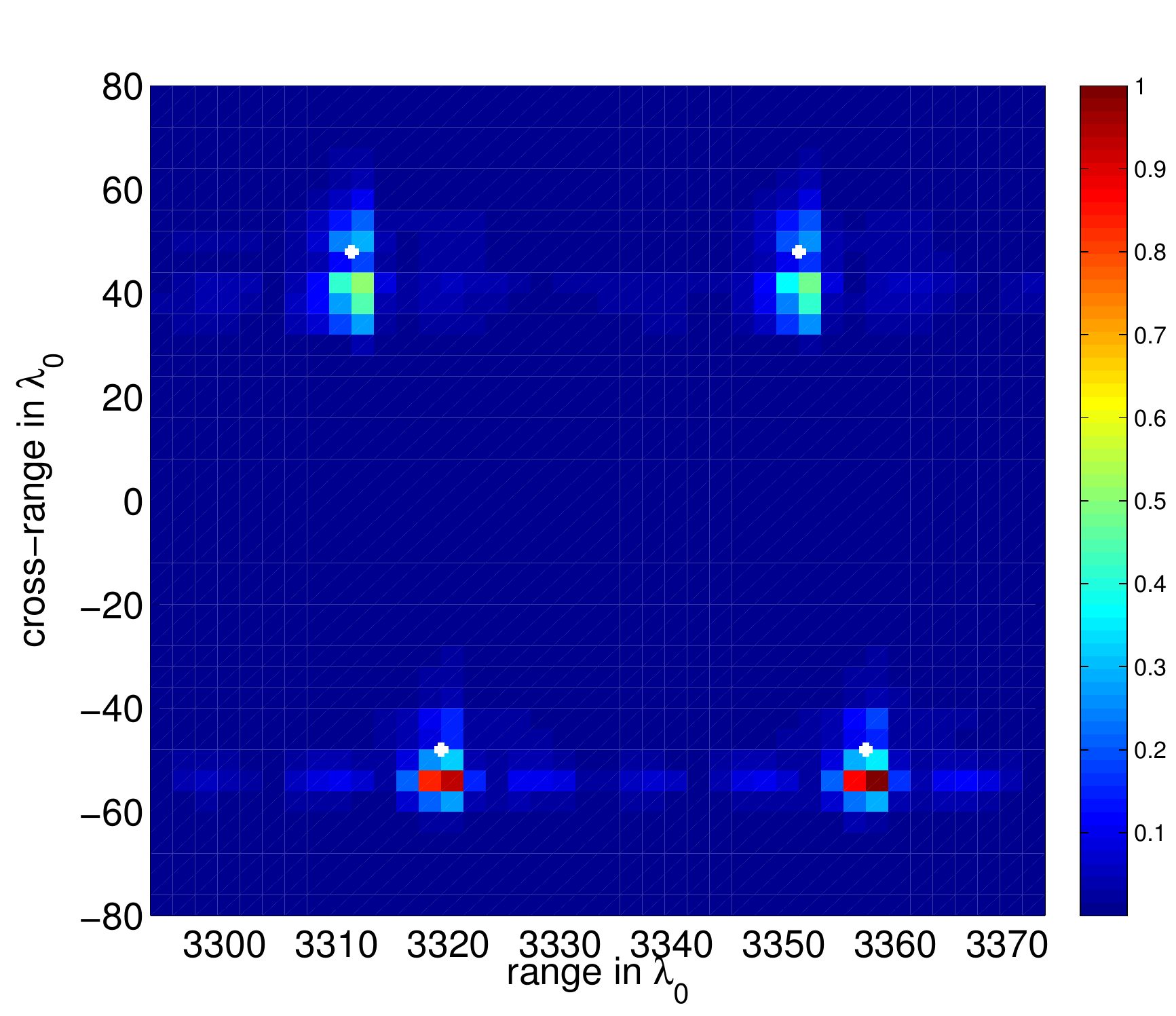}&
\includegraphics[scale=0.22]{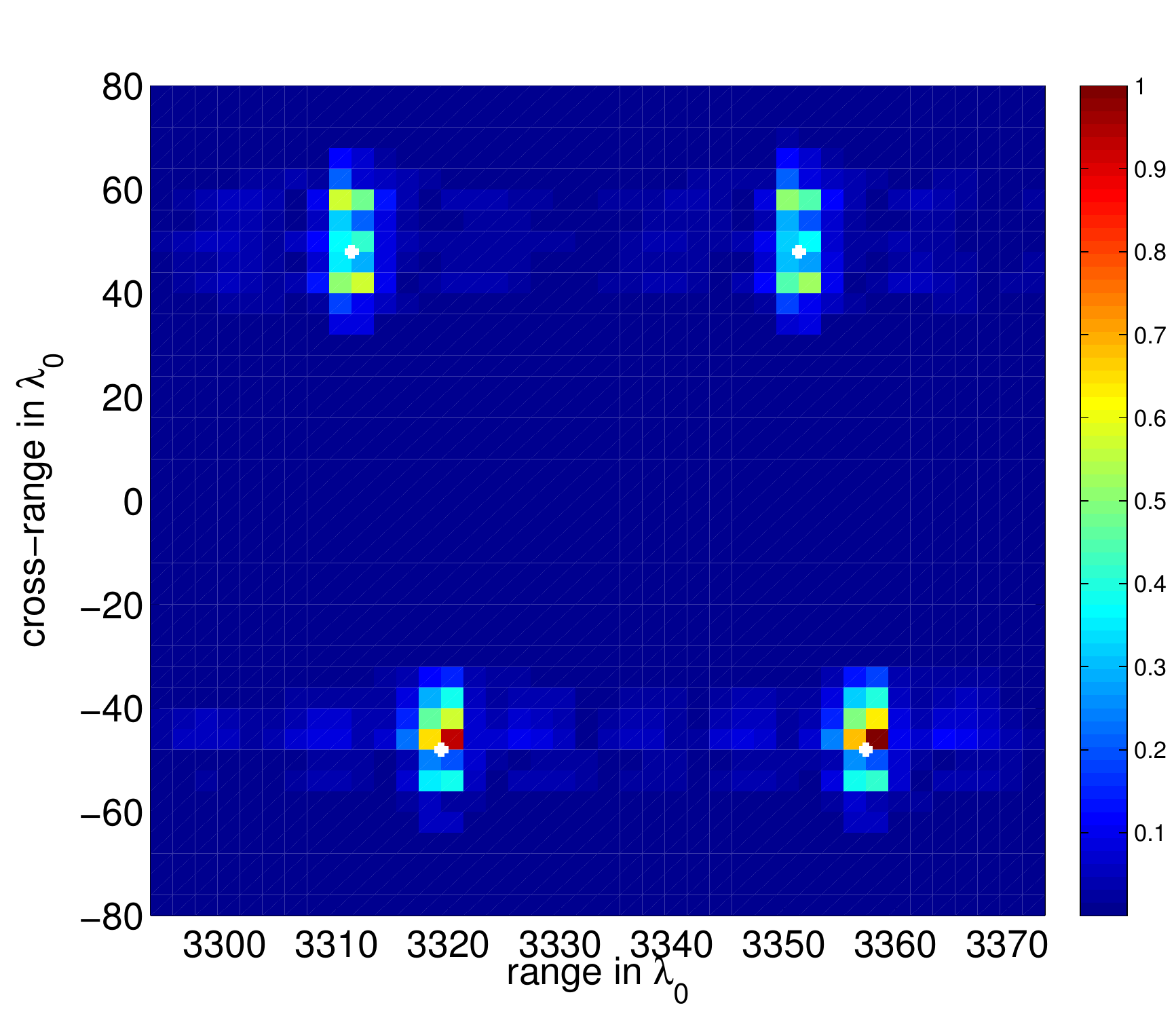}&
\includegraphics[scale=0.22]{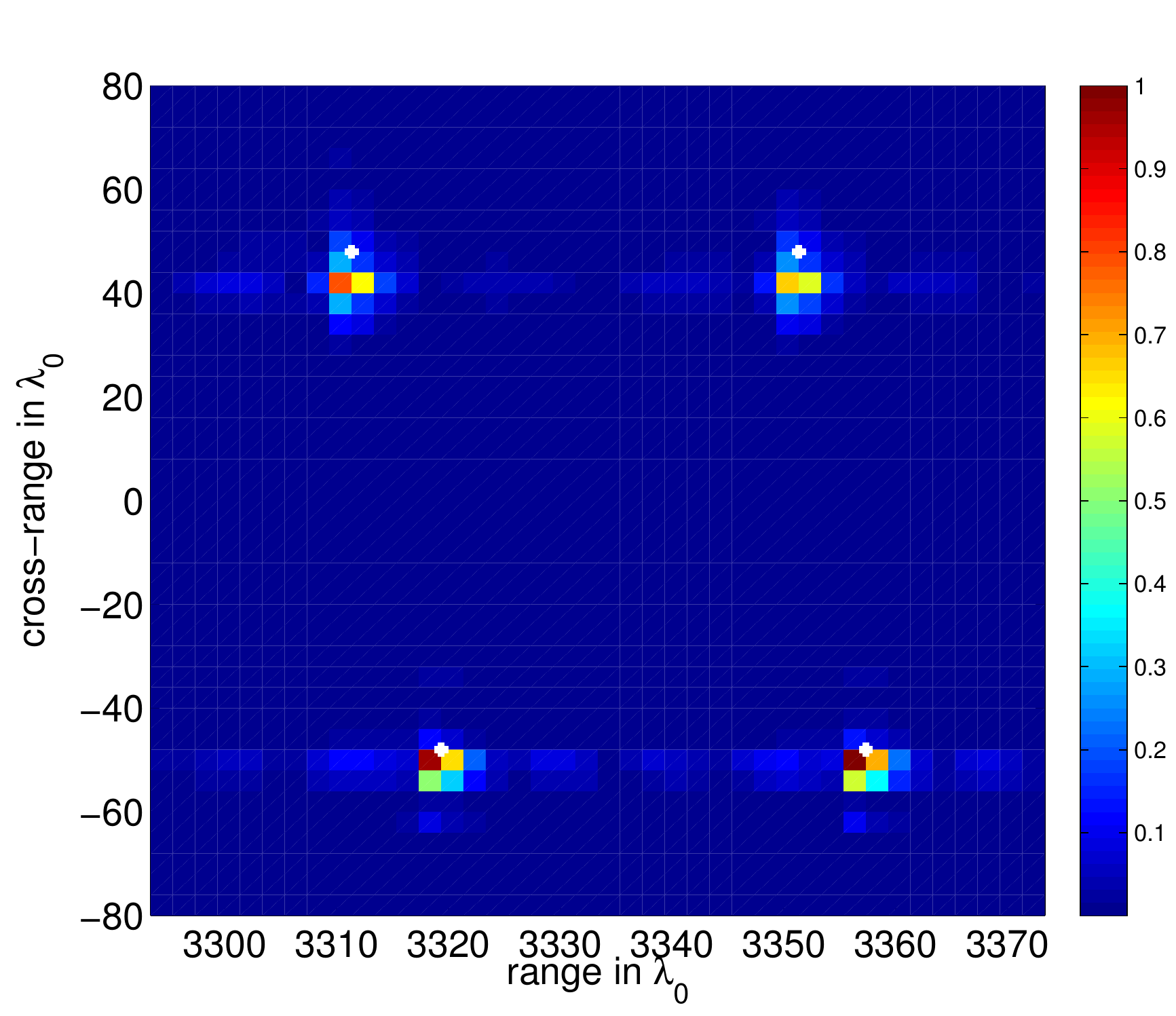} \\
\includegraphics[scale=0.22]{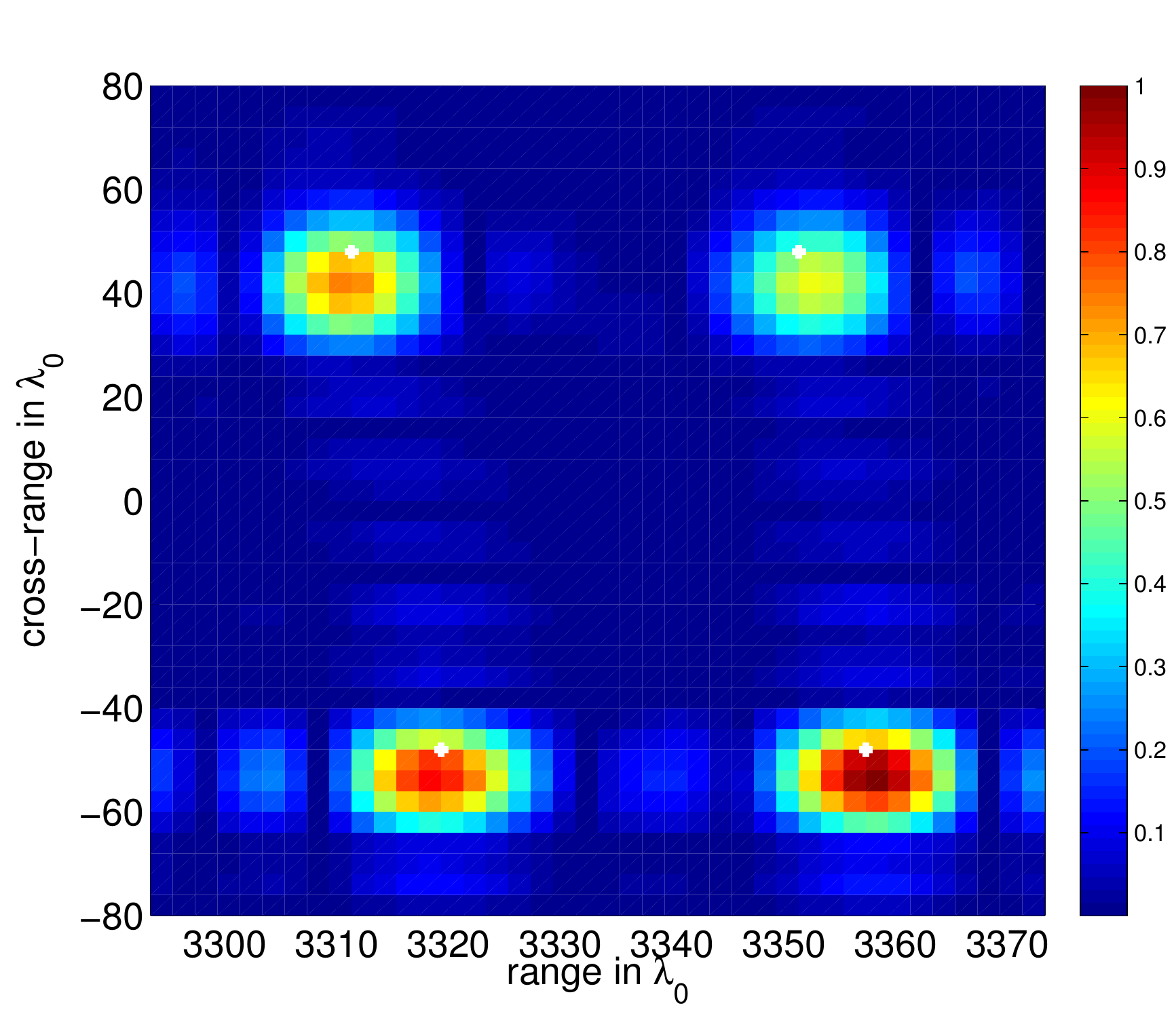}&
\includegraphics[scale=0.22]{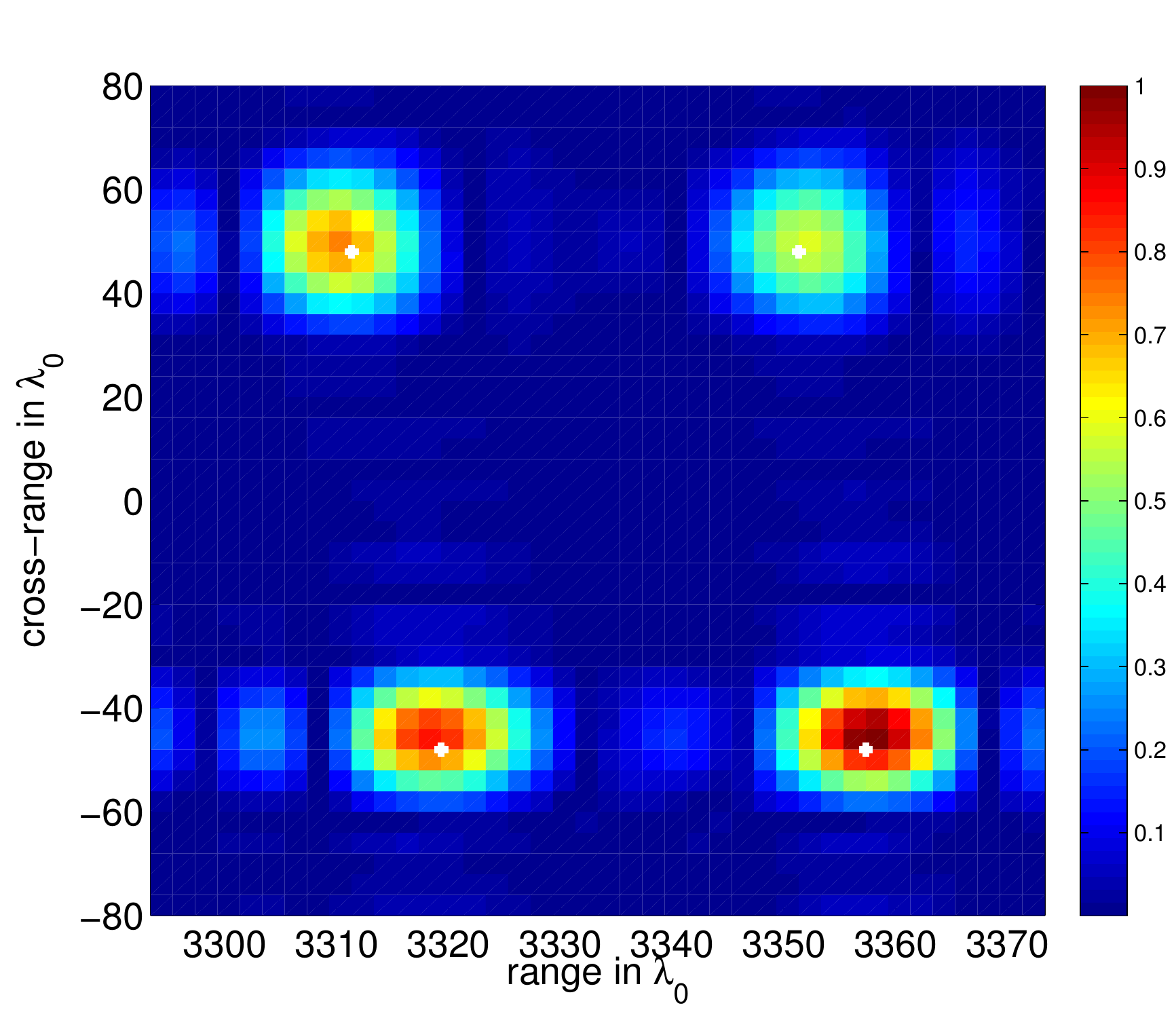}&
\includegraphics[scale=0.22]{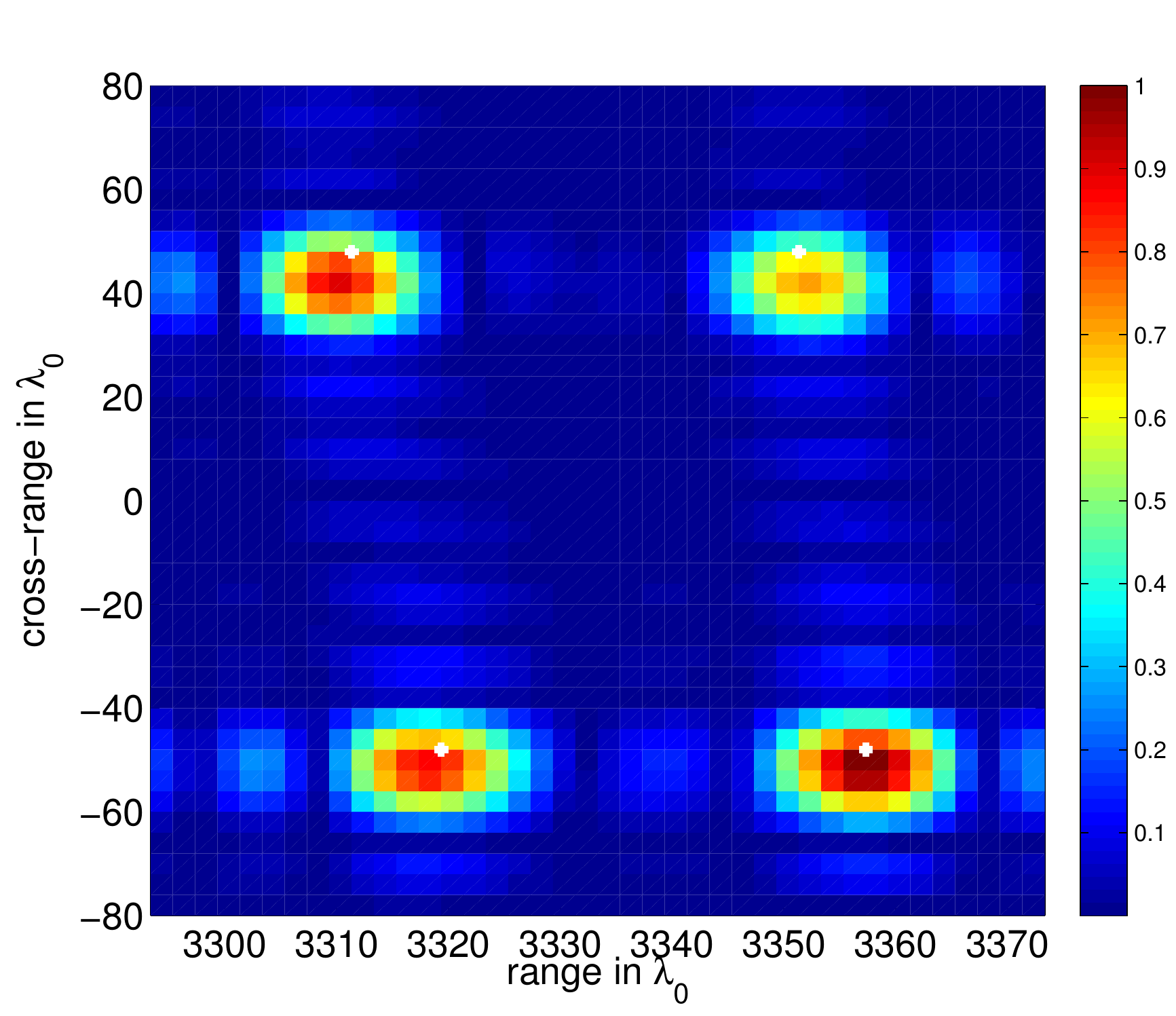}\\
\end{tabular}
\end{center}
\caption{This is the same as Figure \ref{fig:error1} but now we change the propagation distance to $L/3$ and multiply the strength of the fluctuations by a factor $\sqrt{3}$ to keep $\sigma \sqrt{l L}/ \lambda_0$ constant. }
\label{fig:error2}
\end{figure}

In order for SRINT to produce reliable results in inhomogeneous media there must be some coherence in the recovered interferometric data matrix $\Mm_r$. This happens 
when the fluctuations in the phases induced  by the random phase model are of order one (or less). From the discussion in Section \ref{sec:imagingrandom}, it
is seen that the standard deviation of the phases recorded at the array is of order $\sigma \sqrt{l L}/ \lambda_0= \sigma/\sigma_0=\epsilon$.
The numerical results shown in Figures~\ref{fig:error1} and~\ref{fig:error2} confirm this. We observe good images when we change the 
propagation distance to $L/2$ (Figure~\ref{fig:error1}) and $L/3$ ( Figure~\ref{fig:error2}) while multiplying the strength of the fluctuations by a factor 
$\sqrt{2}$ and $\sqrt{3}$, respectively, to keep the parameter $\epsilon$ fixed. In our numerical simulations we have $\epsilon =0.2$.
Note that as the propagation distance is reduced, the resolution in cross-range improves, as seen by 
comparing Figure \ref{fig:stability} with Figures \ref{fig:error1} and \ref{fig:error2}.

\section{Imaging in inhomogeneous background media}
\label{sec:imagingrandom}

Imaging in inhomogeneous random media is fundamentally different from imaging in homogeneous or
smoothly varying media. This is so because when the medium is inhomogeneous we know, at best, the large
scale, but we cannot know its detailed, small scale structure which is impossible to be determined. Hence,
all imaging methods use a homogeneous medium (or smoothly varying one) as a reference medium, even when the 
collected data are affected by the medium inhomogeneities. Indeed, in randomly inhomogeneous media
the data inherit the uncertainty of the fluctuations of the media, resulting in 
wave distortions that lead to space-correlated multiplicative noise, which is very different from additive uncorrelated 
noise usually taken into account in imaging. In fact, many of the
usual imaging methods used in homogeneous (or smoothly varying) media fail, e.g. Kirchhoff migration,  
because the images change unpredictably with the detailed features of the medium and, thus, they become unstable.

Single receiver interferometric imaging deals with wave distortions by restricting, or thresholding, the data to remove 
excessive incoherence as described in the Introduction (see \eqref{eq:SRINT}) and in Section \ref{sec:imaginghomo} (see \eqref{eq:AMA}). 
This stabilizes the images. Note that the thresholding, captured by the parameters $\Omega_d$ and $X_d$, depends 
on the properties of the inhomogeneous, random medium. Using a random phase model, we briefly describe  next how these parameters emerge 
from the properties of the random medium. This model, frequently used to account for weak phase distortions, is analyzed in detail in 
\cite{BGPT-rtt}, which we follow here.

The random phase model characterizes wave propagation in
the high-frequency regime in random media with weak fluctuations and small correlation
lengths $l$ compared to the wavelength $\lambda_0$. It provides an analytical approximation for the Green’s function 
between two points $\vec{\vect x}$ and $\vec{\vect y}$ at a distance $L$ from each other such that $L\gg l \gg \lambda_0$.
This approximation is given by
\begin{equation}\label{eq:random_green_func}
G(\vec{\vect x},\vec{\vect y}; \omega)=G_0(\vec{\vect x},\vec{\vect y}; \omega)
\exp{\left[i \omega \nu (\vx,\vy)\right]}.
\end{equation}
Here, $ G_0$ denotes the Green's function in a homogeneous medium (\ref{homo_green}), and
$\nu (\vx,\vy)$ is the random function 
\begin{eqnarray}
\nu (\vx,\vy) &=& \frac{\sigma |\vx-\vy|}{2c_0} \int_0^1 ds \, \mu \left(
\frac{\vy}{l} + (\vx-\vy)\frac{s}{l} \right)  \, ,
\label{eq:TT120}
\end{eqnarray}
which accounts for the phase distortions induced by the random fluctuations of the wave speed
modeled as 
\begin{equation}\label{eq:random_wave_speed}
\frac{1}{c^2(\vec{\vect x})}=\frac{1}{c_0^2}\bigg(1+\sigma\mu(\frac{\vec{\vect x}}{l})\bigg)\, .
\end{equation}
In \eqref{eq:random_wave_speed}, $c_0$ denotes the average speed, $\sigma$ denotes the strength of the flucutations
with correlation length $l$,  and $\mu(\cdot)$ is a stationary random process of a dimensionless argument with
zero mean and normalized autocorrelation function $R(|\vec{\vect x}-\vec{\vect x'}|)=
\mE(\mu(\vec{\vect x})\mu(\vec{\vect x'}))$, so that $R(0)=1$, and $\int_0^\infty R(r) r^2 dr < \infty$. 
Here, we only consider weak fluctuations such that $\sigma\ll1$. 
The propagation distance $L$ is large, though, so the cumulative scattering effects are $\mathcal{O}(1)$, 
but not too large so the phases of the collected signals still maintain certain degree of coherence. Hence, 
there is some degree of correlation between signals coming from  sources whose locations and frequencies are not very far away.

The model~\eqref{eq:random_green_func}-\eqref{eq:random_wave_speed} is valid when (i) the wavelength $\lambda_0$ is much smaller than the
correlation length $l$ so the geometric optics approximation holds, (ii) the correlation
length $l$ is much smaller than the propagation distance $L$ so the statistics of the phase
become Gaussian, and (iii) the strength of the fluctuations $\sigma$ is small, so the amplitude of
the wave is kept approximately constant, but large enough to ensure that the perturbations of the phases
are not negligible. The last condition holds when 
\begin{equation}
\label{cond}
\frac{\sigma^2L^3}{l^3}\ll\frac{\lambda^2}{\sigma^2lL} \sim 1\, ,
\end{equation}
as discussed in detail in~\cite{Tatarski61,Rytov89,BGPT-rtt}. Note that although we take weak fluctuations, 
the distortion of the wavefronts caused by the inhomogeneities of the medium is observable
because the waves travel long distances. Comparing \eqref{eq:random_green_func} to the homogeneous Green's function \eqref{homo_green} we see that, in this
regime, only the phases of the waves are perturbed by the random medium, while the amplitudes remain unchanged. 


\subsection*{The moments of the random function $\nu$}Assume that $\vy=(0,L)$ is in the IW and $\vx=(\bx,0)$ is at the array, 
then the distance between them is of order $L$, $L \gg l$. 
Suppose  $\mu(\cdot)$ is statistically homogeneous and Gaussian, then 
one can show (see Lemma 3.1 in~\cite{BGPT-rtt}, or Appendix \ref{appendix:moments_nu} of this paper) 
that
the random process 
\begin{equation}
\label{app:defnu}
\nu(\bx) := \nu \big( \vx=(\bx,0),\vy\big), 
\end{equation}
has Gaussian statistics with mean zero and
covariance function
\begin{equation}
\mE\left\{ \nu(\bx)\nu(\bx') \right\} \approx \tau_c^2
\cC\left(\frac{|\bx-\bx'|}{l} \right),
\label{app:lem.1p}
\end{equation}
where the variance is 
\begin{equation}
\tau_c^2 = \frac{\sqrt{2 \pi} \sigma^2 l L}{4 c_0^2},
\label{app:lem.1pp}
\end{equation}
 and 
 $$\cC(r) = \frac{1}{r} \int_0^r du \, e^{-u^2/2}.$$ 

Note that $\tau_c$ has dimensions of time. For a better understanding of the  parameters that are meaningful in the random phase model, 
we use adimensional variables. We scale all length variables by the central wavelength $\lambda_0$ 
\begin{eqnarray}
{\tilde \bx} = \frac{\bx}{ \lambda_0},\quad  {\tilde a} = \frac{a}{ \lambda_0},\quad  {\tilde l} = \frac{l }{ \lambda_0}, \quad  {\tilde L} = \frac{L }{ \lambda_0}, 
\label{app:ad1}
\end{eqnarray}
and the frequency (respectively time) by the central frequency $\omega_0=2\pi c_0/ \lambda_0$  (respectively $1/\omega_0$) 
\begin{eqnarray}
\quad  {\tilde \omega} = \frac{\omega}{\omega_0}\, ,  
\quad  {\tilde \tau_c} =  \tau_c \omega_0\, .
\label{app:ad2}
\end{eqnarray}
We also introduce the dimensionless parameter 
\begin{eqnarray}
\label{eq:sigma0}
\sigma_0= \lambda_0/\sqrt{lL} =1/\sqrt{\tilde l \tilde L},
\end{eqnarray}
which is a characteristic strength of the fluctuations of the inhomogeneities
for which the standard deviation of the random phase fluctuations of the collected signals is $\mathcal{O}(1)$. Note that it is $\sigma_0^2/\sigma^2$ that appears in \eqref{cond} and should be close to one. 
We therefore define the strength of the fluctuations $\sigma$ in terms of $\sigma_0$ and, thus, define the dimensionless
parameter 
\begin{equation}
 \epsilon = \frac{\sigma}{\sigma_0}\,. 
\label{app:epsilon}
\end{equation}
Using these dimensionless variables and parameters, we have 
\begin{equation}
 {\tilde \tau_c}^2 = \pi^2 \sqrt{2\pi} \epsilon^2\, ,
\end{equation}
which shows that the variance of the fluctuations of the random phases at the array only depends on the standard deviation
$\epsilon$ of the fluctuations of the wave speed measured in units of $\sigma_0$. This means in turn that $\tilde \tau_c$ 
only depends on the product of the dimensionless parameters ${\tilde l}$ and ${\tilde L}$.

\subsection*{The threshold parameters and moment formulas}
\label{sec:statstability}

The stability analysis of the imaging functionals (\ref{eq:SRINT}) or (\ref{eq:AMA}) relies on computations involving statistical moments of the 
Green's function \eqref{eq:random_green_func}. The detailed stability analysis is in \cite{BGPT-rtt} and will not be repeated here. We will
only show the first two moment formulas of the Green's function where the threshold parameters appear. Higher order moments are computed using the Gaussian property of the random function $\nu$. From the moment formulas we can see
that thresholding in  (\ref{eq:SRINT}) or (\ref{eq:AMA}) is a form of denoising. Another way to interpret thresholding is as the removal of relatively
incoherent data that will not contribute to stable image formation.

Assume, as before, that $\vy=(0,L)$ is in the IW, and $\vx=(\bx,0)$, $\vx'=(\bx',0)$ are the positions of two array elements.
Then, we can show that (see Lemma 3.2 in~\cite{BGPT-rtt}, or 
Appendix \ref{appendix:moments_expnu}) 
\begin{eqnarray}
\label{eq:meanC10} 
\mE \left\{e^{ i \om \nu(\bx)}  \right\} 
&=&
 \mbox{exp} \left\{ -
\frac{\om^2  \tau_c^2}{2} \right\}, 
\end{eqnarray}
and
\begin{eqnarray}
\mE \left\{e^{ i \om \nu(\bx) - i \om' \nu(\bx')}\right\} 
&=& \mbox{exp} \left\{ -
\frac{(\om-\om')^2 \tau_c^2}{2} - \om \om' \tau_c^2 \left[ 1-
  \cC\left(\frac{|\bx-\bx'|}{l}\right) \right] \right\}.
\label{eq:meanC1} 
\end{eqnarray}
If, in addition, the array elements are nearby so that $|\bx-\bx'|\ll{l}$ (and, thus, 
we can expand the covariance function around zero so $\cC(r)=1-r^2/6 +O(r^4)$), and the bandwidth is relatively small so that 
$\om \om' \approx \om_0^2$, we get that 
\begin{equation}
\mE \left\{e^{ i \om \nu(\bx) - i \om' \nu(\bx')}\right\} \approx  \mbox{exp} \left\{ -
  \frac{(\om-\om')^2 }{2 \Omega_d^2} - \frac{|\bx-\bx'|^2}{2 X_d^2} \right\},
\label{eq:meanC2} 
\end{equation}
with
\begin{equation}
\Omega_d = \frac{1}{\tau_c} 
, \quad \quad 
X_d = \frac{\sqrt{3} l}{\om_o \tau_c}  \,.
\label{eq:meanC3} 
\end{equation}
In dimensionless units, ${\tilde \Omega}_d = 1/{\tilde \tau_c}\approx 0.2/ \epsilon$ and 
${\tilde X_d} = \sqrt{3}{\tilde l}/{\tilde \tau_c} \approx 0.35\,{\tilde l}/\epsilon$. 
We deduce that the dimensionless decoherence frequency ${\tilde \Omega}_d$ only depends on $\epsilon$, while the dimensionless decoherence length 
${\tilde X}_d$ also depends on the dimensionless correlation length ${\tilde l}$. 

We note the following concerning the Green's function between a point $\vx=(\bx,0)$ on the array and a point $\vy=(0,L)$ 
in the IW: the Green's function $G(\vx,\vy;\om)$ that models wave propagation between these two points is a random process
with mean
\begin{equation}
\mE \left\{G(\om,\vx,\vy)\right\} =  G_0(\om,\vx,\vy) 
\mE \left\{e^{i \om \nu(\bx)}\right\} =
G_0(\om,\vx,\vy) e^{-\frac{\om^2 \tau_c^2}{2}},
\label{eq:meanG}
\end{equation}
and  variance 
\begin{equation}
Var \left\{ G(\om,\vx,\vy)\right\} =
\left|  G_0(\om,\vx,\vy)  \right|^2 
\left( 1- e^{-\om^2 \tau_c^2} \right).
\label{eq:VarG}
\end{equation}
 
Thus, according to \eqref{eq:meanG} and \eqref{eq:VarG}, the mean of the Green's function goes to zero when $\om \tau_c$ is large, while 
the variance remains bounded. This instability will be inherited by any imaging functional that backpropagates these data, as recorded on the array,
in a homogeneous background medium. 
    
If, instead,  we backpropagate 
interferometric data, we have to consider the random process $G(\vx,\vy;\om) \overline{G(\vx',\vy;\om')}$. The mean and the  variance 
of this random process is given by
 \begin{equation}
\begin{array}{r}
\mE \left\{  G(\vx,\vy;\om) \overline{G(\vx',\vy;\om')} \right\} \dsp =  G_0(\vx,\vy;\om) \overline{G_0(\vx',\vy;\om)} 
\mE \left\{e^{i \om \nu(\bx) - i \om' \nu(\bx')} \right\} \\[12pt]
 =\dsp  G_0(\vx,\vy;\om) \overline{G_0(\vx',\vy;\om)} 
\mbox{exp} \left\{ -
  \frac{(\om-\om')^2 }{2 \Omega_d^2} - \frac{|\bx-\bx'|^2}{2 X_d^2} \right\}
\end{array}
\label{eq:meanGG}
\end{equation}
and 
\begin{equation}
\begin{array}{r}
Var\left\{ G(\vx,\vy;\om) \overline{G(\vx',\vy;\om')} \right\} = \hspace*{6cm}\\
\dsp \left| G_0(\vx,\vy;\om)  \overline{G_0(\vx',\vy;\om)}  \right|^2 
\left( 1-  \mbox{exp} \left\{ -
  \frac{(\om-\om')^2 }{\Omega_d^2} - \frac{|\bx-\bx'|^2}{ X_d^2} \right\}\right),
\label{eq:VarGG}
\end{array}
\end{equation}
respectively. Hence, 
we observe that the variance  goes to zero as $| \bx -\bx'| \to 0$  and $|\om-\om'| \to 0$.
Moreover, as $| \bx -\bx'| \to 0$  and $|\om-\om'| \to 0$ the expected value  reduces to the corresponding quantity in 
 a homogeneous medium. This means, that the restricted multifrequency interferometric data converge to the corresponding deterministic quantities. 
Hence, if we use SRINT with masks, so that $| \bx -\bx'| \ll X_d$ and $|\om -\om'|\ll\Omega_d$,
the  multifrequency interferometric data 
used for imaging are restricted appropriately to data that lead to statistically stable images. 
The fact that we gain stability by thresholding but we also lose some resolution is
discussed in \cite{BGPT-rtt}.

\section{Concluding remarks}
\label{sec:conclusions}

In this paper, we introduce a holography-based approach for imaging with intensity-only measurements. We show that by controlling the illuminations, 
we can recover the multifrequency interferometric data matrix (\ref{eq:Mr}) from intensities recorded  at a single receiver. The recovered data matrix $\Mm_r$ can then be used for imaging, which can be done by back-propagation. This allows us to reconstruct a full three-dimensional image, including depth information, from intensity-only measurements. Moreover, we show that in homogeneous media there is no resolution loss compared to imaging with full data, including phases. 
We also consider inhomogeneous media where wavefront distortions can arise. To image in such media we restrict the interferometric data to nearby 
frequencies $|\omega-\omega'|\le \Omega_d$ and nearby receivers $|\bx -\bx'|\le X_d$. An efficient matrix implementation of this thresholding operation is introduced, using a mask that sets to zero all the entries in the multifrequency interferometric data matrix
that do not satisfy the required coherence constraints. 
The robustness of the interferometric approach with respect to the image window discretization as well as in the case of wavefront distortions is explored with numerical simulations carried out in an optical (digital) microscopy imaging regime.

\section*{Acknowledgment} 

Miguel Moscoso's work was partially supported by  AFOSR FA9550-14-1-0275 and the Spanish MICINN grant
FIS2013-41802-R. The work of George Papanicolaou and Chrysoula Tsogka was partially supported by
AFOSR grant FA9550-14-1-0275. Alexei Novikov's work was partially supported by  AFOSR FA9550-14-1-0275 and NSF DMS-1515187.

\appendix

\section{Moments of the random process $\nu(\bx)$} \label{appendix:moments_nu}
Because $\mu$ is a zero-mean stationary random process, $\mE\{\nu(\bx)\}=0$. Next, we compute the second order moment
of the random process $\nu(\bx)=\nu((\bx,0),(0,L))$. Assuming that $ | \bx | \ll L$,  
\begin{eqnarray}
\mE\left\{ \nu(\bx)\nu(\bx') \right\} &=& \frac{\sigma^2 L^2}{4 c_0^2}\int_{0}^{1} ds \int_{0}^{1} ds'
\mE \left\{ \mu(\frac{\vy + s(\vx-\vy)}{l}) \mu(\frac{\vy + s'(\vx'-\vy)}{l})\right\} \nonumber \\
&=& \frac{\sigma^2 L^2}{4 c_0^2}\int_{0}^{1} ds \int_{0}^{1} ds' R\left( \frac{| s(\vx-\vy) - s'(\vx'-\vy) |}{l}  \right) \nonumber \\
&=& \frac{\sigma^2 L^2}{4 c_0^2}\int_{0}^{1} ds \int_{0}^{1} ds' \exp\{-\frac{| s(\vx-\vy) - s'(\vx'-\vy) |^2}{2\l^2}\} \nonumber \\
&=& \frac{\sigma^2 L^2}{4 c_0^2}\int_{0}^{1} ds \int_{0}^{1} ds' \exp\{-\frac{| s \bx - s' \bx' |^2}{2\l^2} - \frac{(s-s')^2L^2}{2l^2}\}  \nonumber \\
&=& \frac{\sigma^2 L^2}{4 c_0^2}\int_{0}^{1} ds \int_{-s}^{1-s} ds'' \exp\{-\frac{|s\bx-(s''+s)\bx'|^2}{2\l^2} - \frac{s''^2L^2}{2l^2}\} \, .
\end{eqnarray}
If $L\gg l$, the $s''$ integral is approximately the integral of a Gaussian that can be extended to the real line, and whose value
is $\sqrt{2\pi}\,l/L$. Hence,
\begin{eqnarray}
\mE\left\{ \nu(\bx)\nu(\bx') \right\} &\approx& 
\frac{\sqrt{2\pi}\sigma^2 l L}{4 c_0^2}\int_{0}^{1} ds \exp\{-\frac{s^2|\bx-\bx'|^2}{2\l^2}\} \nonumber \\
&=& \frac{\sqrt{2\pi}\sigma^2 l L}{4 c_0^2}\frac{l}{|\bx-\bx'|}\int_{0}^{|\bx-\bx'|/l} du \exp\{-u^2/2\}  \nonumber \\
&=& \tau_c^2 \cC(|\bx-\bx'|/l)\, ,
\label{app:A2}
\end{eqnarray}
with $\tau_c^2 = \frac{\sqrt{2 \pi} \sigma^2 l L}{4 c_0^2}$, and $\cC(r) = \frac{1}{r} \int_0^r du \, e^{-u^2/2}$. 

\section{Moments of the random process $e^{i\omega\nu({\bf x})}$} \label{appendix:moments_expnu}
To compute the first moment of the random process $e^{i\omega\nu({\bf x})}$ we use that the expectation of an exponential function
is $\mE\left\{ e^{aX}\right\}=e^{a^2/2}$ if $X\sim N(0,1)$. Hence,
\begin{eqnarray}
\mE\left\{ e^{i\omega\nu({\bx})} \right\} = e^{-\omega^2\mE\left\{[\nu({\bx})]^2\right\}/2} = e^{-\omega^2\tau_c^2/2} \, ,
\end{eqnarray}
where we have used \eqref{app:A2} with $\bx=\bx'$. To obtain the second order moment 
\begin{eqnarray}
\mE\left\{ e^{i\omega\nu({\bx}) -i\omega'\nu({\bx'})} \right\} = 
e^{-\mE\left\{[\omega\nu({\bx}) -\omega'\nu({\bx'})]^2\right\}/2}
\label{app:B2}
\end{eqnarray}
we compute the expectation
\begin{eqnarray}
\mE\left\{ [\omega \nu(\bx)- \omega'\nu(\bx')]^2 \right\} &=& 
\omega^2\mE\left\{ [\nu(\bx) ]^2 \right\} + \omega'^2\mE\left\{ [\nu(\bx') ]^2 \right\} - 2 \omega\omega'\mE\left\{ [\nu(\bx)\nu(\bx')] \right\} \nonumber \\
&=&(\om-\om')^2  \mE\left\{ [\nu(\bx) ]^2 \right\} + 2\om\om'\mE\left\{ [\nu(\bx) ]^2 \right\}- 2 \omega\omega'\mE\left\{ [\nu(\bx)\nu(\bx')] \right\} \nonumber \\
&=&(\om-\om')^2 \tau_c^2 + 2\om\om'\tau_c^2 - 2 \omega\omega'\tau_c^2\cC(|\bx-\bx'|/l)\, ,
\label{app:B3}
\end{eqnarray}
where we have used \eqref{app:A2}. Inserting \eqref{app:B3} into \eqref{app:B3} we obtain
\begin{eqnarray}
\mE\left\{ e^{i\omega\nu({\bx}) -i\omega'\nu({\bx'})} \right\} = 
e^{-(\om-\om')^2 \tau_c^2/2 - \omega\omega'\tau_c^2(1-\cC(|\bx-\bx'|/l))}\,.
\label{app:B4}
\end{eqnarray}
If, in addition, $|\bx-\bx'|/l \ll 1$ so $\cC(|\bx-\bx'|/l))\approx 1 - |\bx-\bx'|^2/6l^2 + \dots$, and assuming $\om\om'\approx\om_0^2$, we find that
\begin{eqnarray}
\mE\left\{ e^{i\omega\nu({\bx}) -i\omega'\nu({\bx'})} \right\} = 
e^{-(\om-\om')^2 /2\Omega_d^2 - |\bx-\bx'|^2/2X_d^2}\,,
\label{app:B5}
\end{eqnarray}
with $\Omega_d = 1/\tau_c$ and $X_d = \sqrt{3} l/\om_o \tau_c$.

\end{document}

%% file: mathmacros.tex
\newcommand{\om}{\omega}

\newcommand{\vy}{{\vec{\mathbf{y}}}}



%% file: setup2.pdf_t
\begin{picture}(0,0)%
\includegraphics{setup2.pdf}%
\end{picture}%
\setlength{\unitlength}{1579sp}%
\begingroup\makeatletter\ifx\SetFigFont\undefined%
\gdef\SetFigFont#1#2#3#4#5{%
  \reset@font\fontsize{#1}{#2pt}%
  \fontfamily{#3}\fontseries{#4}\fontshape{#5}%
  \selectfont}%
\fi\endgroup%
\begin{picture}(11024,6695)(399,-6908)
\put(9451,-4411){\makebox(0,0)[lb]{\smash{{\SetFigFont{9}{10.8}{\familydefault}{\mddefault}{\updefault}{\color[rgb]{0,0,0}IW}%
}}}}
\put(1651,-1486){\makebox(0,0)[lb]{\smash{{\SetFigFont{9}{10.8}{\familydefault}{\mddefault}{\updefault}{\color[rgb]{0,0,0}${\vect x}_r$}%
}}}}
\put(1951,-4111){\makebox(0,0)[lb]{\smash{{\SetFigFont{9}{10.8}{\familydefault}{\mddefault}{\updefault}{\color[rgb]{0,0,0}$\lambda$}%
}}}}
\put(1501,-4786){\makebox(0,0)[lb]{\smash{{\SetFigFont{9}{10.8}{\familydefault}{\mddefault}{\updefault}{\color[rgb]{0,0,0}${\vect x}_s$}%
}}}}
\put(4126,-5086){\makebox(0,0)[lb]{\smash{{\SetFigFont{9}{10.8}{\familydefault}{\mddefault}{\updefault}{\color[rgb]{0,0,0}$L$}%
}}}}
\put(526,-4036){\makebox(0,0)[lb]{\smash{{\SetFigFont{9}{10.8}{\familydefault}{\mddefault}{\updefault}{\color[rgb]{0,0,0}$a$}%
}}}}
\put(8701,-3061){\makebox(0,0)[lb]{\smash{{\SetFigFont{9}{10.8}{\familydefault}{\mddefault}{\updefault}{\color[rgb]{0,0,0}${\vect y}_j$}%
}}}}
\put(931,-4951){\makebox(0,0)[lb]{\smash{{\SetFigFont{9}{10.8}{\familydefault}{\mddefault}{\updefault}{\color[rgb]{0,0,0}$h$}%
}}}}
\end{picture}%